\documentclass{amsart}

\usepackage[utf8]{inputenc}

\usepackage{amsmath,amssymb,amsfonts,mathtools}
\usepackage{bm}
\usepackage{accents}

\usepackage{color}
\usepackage{graphicx}
\usepackage{epsf}

\usepackage[margin=1.5in,top=1in,bottom=1in]{geometry}

\usepackage{enumerate}
\usepackage[noautoscale]{youngtab}

\usepackage[cmtip,arrow]{xy}
\xyoption{all}
\usepackage{pb-diagram,pb-xy}

\usepackage{subcaption}

\usepackage{tikz}
\usepackage{tikz-cd}

\usepackage{float}
\usepackage{algorithm}
\usepackage{algpseudocode}
\usepackage{etoolbox}

\ifx\pdfpageheight\undefined
    \usepackage[dvips,colorlinks=true,linkcolor=blue,citecolor=red,urlcolor=blue]{hyperref}
    \makeatletter
    \edef\Gin@extensions{\Gin@extensions,.mps}
    \makeatother
\else
    \usepackage[
        pdftex=true,
        bookmarksopen=false,
        breaklinks=true,
        backref=page,
        pagebackref=true,
        plainpages=false,
        hyperindex=true,
        pdfstartview=FitH,
        colorlinks=true,
        linkcolor=blue,
        citecolor=red,
        urlcolor=blue,
        pdfpagelabels=true,
        hypertexnames=false
    ]{hyperref}
\fi

\makeatletter
\tikzset{
    column sep/.code=\def\pgfmatrixcolumnsep{\pgf@matrix@xscale*(#1)},
    row sep/.code=\def\pgfmatrixrowsep{\pgf@matrix@yscale*(#1)},
    matrix xscale/.code=\pgfmathsetmacro\pgf@matrix@xscale{\pgf@matrix@xscale*(#1)},
    matrix yscale/.code=\pgfmathsetmacro\pgf@matrix@yscale{\pgf@matrix@yscale*(#1)},
    matrix scale/.style={/tikz/matrix xscale={#1},/tikz/matrix yscale={#1}}
}
\def\pgf@matrix@xscale{1}
\def\pgf@matrix@yscale{1}
\makeatother

\newtheorem{theorem}{Theorem}
\newtheorem{lemma}{Lemma}[section]
\newtheorem{corollary}{Corollary}
\newtheorem{proposition}{Proposition}[section]

\newtheorem*{claim*}{Claim}

\newtheorem{assumption}{Assumption}

\newtheorem*{theorem*}{Theorem}
\newtheorem*{corollary*}{Corollary}

\theoremstyle{definition}
\newtheorem{definition}{Definition}[section]
\newtheorem{example}{Example}

\newtheorem{notation}{Notation}

\algnewcommand\algorithmicinput{\textbf{Input:}}
\algnewcommand\INPUT{\item[\algorithmicinput]}

\algnewcommand\algorithmicoutput{\textbf{Output:}}
\algnewcommand\OUTPUT{\item[\algorithmicoutput]}

\algnewcommand\algorithmicproc{\textbf{Procedure:}}
\algnewcommand\PROCEDURE{\item[\algorithmicproc]}

\algnewcommand\algorithmiccomplexity{\textbf{Complexity:}}
\algnewcommand\COMPLEXITY{\item[\algorithmiccomplexity]}

\newlength{\continueindent}
\makeatletter
\newcommand*{\ALG@customparshape}{%
    \parshape 2
    \leftmargin \linewidth
    \dimexpr\ALG@tlm+\continueindent\relax
    \dimexpr\linewidth+\leftmargin-\ALG@tlm-\continueindent\relax
}
\apptocmd{\ALG@beginblock}{\ALG@customparshape}{}{\errmessage{failed to patch}}
\makeatother

\newcommand{\SDO}{$\texttt{SDO}$}

\theoremstyle{remark}
\newtheorem{remark}{Remark}

\theoremstyle{observation}

\definecolor{DarkBlue}{rgb}{0,0.1,0.55}

\numberwithin{equation}{section}

\newcommand {\hide}[1]{}

\newcommand{\zero}{\mbox{\bf zero}}

\newcommand {\junk}[1]{}

\newcommand {\Der} {{\rm Der}}

\newcommand {\la}   {{\langle}}
\newcommand {\ra}   {{\rangle}}

\def\addots{\mathinner{\mkern1mu
		\raise1pt\vbox{\kern7pt\hbox{.}}
		\mkern2mu\raise4pt\hbox{.}\mkern2mu
		\raise7pt\hbox{.}\mkern1mu}}

\DeclareMathOperator{\trace}{Tr}

\DeclareMathOperator{\svectorize} {svec}

\DeclareMathOperator{\ri} {ri}

\newcommand{\x}{\mathbf{x}}
\newcommand{\X}{\mathbf{X}}

\newcommand{\y}{\mathbf{y}}

\newcommand{\LO}{$\texttt{LO}$}

\newcommand{\LCP}{$\texttt{LCP}$}
\newcommand{\AHO}{$\texttt{AHO}$}
\newcommand{\IPM}{$\texttt{IPM}$}

	\title[Ramification of Weighted Central Paths and Super-linear IPMs for SDO]
	{
		Ramification of Weighted Central Paths and Super-linear Interior-Point Methods for Semi-definite Optimization
	}

    	\author[Mohammad-Nezhad]{Ali Mohammad-Nezhad}
	\address{Department of Statistics and OR,
            University of North Carolina at Chapel Hill, Chapel Hill, NC 27599, U.S.A.}
	\email{alimn@unc.edu}

        \author[Sampourmahani]{Pouya Sampourmahani}
        \address{Department of Industrial and Systems Engineering, Lehigh University, Bethlehem, PA 18015, U.S.A.}
	\email{pos220@lehigh.edu}
        \author[Terlaky]  {Tam\'as Terlaky}
        \address{Department of Industrial and Systems Engineering, Lehigh University, Bethlehem, PA 18015, U.S.A.}
	\email{terlaky@lehigh.edu}

\begin{document}
\begin{abstract}
The local convergence of primal-dual interior-point methods for semi-definite optimization is closely tied to the regularity of the central path at its limiting point. Without strict complementarity, the central path can exhibit fractional power asymptotics and need not be differentiable in the ordinary barrier parameter.

We show that this obstruction can be removed simultaneously for a
family of weighted \AHO~central paths. More precisely, there exists a
common ramification exponent $N$ such that all weighted paths in a
fixed neighborhood of the identity become analytic at the limiting
point after the reparametrization $\mu\mapsto\mu^N$. Using
higher-order derivatives of these reparametrized paths, we construct
a feasible $p$th-order predictor-corrector interior-point method.

Under a uniform predictor regularity condition, and without assuming
strict complementarity, we prove that the normalized complementarity gaps converge
with $Q$-order at least $(p+1)/\rho$, where $\rho$ is the reparametrization exponent used by the method. We also establish super-linear convergence of the simultaneously vanishing primal-dual blocks, derive an explicit quantitative reduction bound throughout the predictor regularity region, and obtain a doubly logarithmic outer iteration bound once the iterates enter the local super-linear regime. Finally, we illustrate the predicted
higher-order behavior on a collection of degenerate and
non-degenerate semi-definite optimization problems.
\end{abstract}
\subjclass[2020]{Primary 90C22, 90C51; Secondary 14P10}

\keywords{Semi-definite optimization, interior-point methods, weighted central paths, semi-algebraic functions, Puiseux expansion, ramification index}

\maketitle
\tableofcontents

\section{Introduction}\label{intro}
Semi-definite optimization (\SDO) is the minimization of a linear objective over the intersection of an affine subspace and the cone of positive semi-definite matrices. Let $\mathbb{S}^n$ be the vector space of symmetric matrices endowed with the inner product $\langle X, S\rangle:=\trace(XS)$. A pair of primal-dual \SDO~problems is defined as 
\begin{align*}
&(\mathrm{P}) \qquad z^*_p:=\inf_{X \in \mathbb{S}^{n}} \left\{\langle C, X \rangle \mid \langle A_i, X \rangle =b_i, \ i=1,\ldots,m, \ X \succeq 0 \right\},\\
&(\mathrm{D}) \qquad z^*_d:=\sup_{(y,S) \in \mathbb{R}^m \times \mathbb{S}^n} \left \{b^Ty  \mid  \sum_{i=1}^m y_i A_{i}+S=C, \ S \succeq 0, \ y \in \mathbb{R}^m \right \},
\end{align*}
where $\succeq 0$ means positive semi-definite.
\begin{notation}\label{primal_dual_optimal_sets}
We use the notation $(\cdot,\cdot,\cdots,\cdot)$ to represent vectors of $\mathbb{R}^n$ or an arbitrary vector space. A primal-dual vector of \SDO, identified by $(X,y,S) \in \mathbb{S}^n \times \mathbb{R}^m \times \mathbb{S}^n$, is called an \textit{optimal solution} if $X \in \mathcal{P}^*$ and $(y,S) \in   \mathcal{D}^*$, where 
\begin{align*}
&\mathcal{P}^*:=\left\{X \mid  \langle A_{i}, X \rangle=b_i, \quad  i=1,\ldots, m, \ \ X \succeq 0, \ \ \langle C,  X \rangle = z^*_p  \right\}, \\
&\mathcal{D}^*:=\left\{(y,S) \mid \sum_{i=1}^m y_i A_{i}+S=C, \ \ S \succeq 0,  \ b^Ty=z^*_d \right\}.
\end{align*}
For a primal-dual feasible vector $(X,y,S)$, we define the
\emph{normalized complementarity gap} by
\begin{align*}
\mu(X,S):=\frac{\langle X,S\rangle}{n}.
\end{align*}
\end{notation}
\noindent
To avoid trivialities, we make the following assumptions from this point on, which in turn guarantee that $\mathcal{P}^*,\mathcal{D}^*$ are non-empty and compact, and $z^*_p=z^*_d$ (see e.g.,~\cite[Corollary~4.2]{T01}).
\begin{assumption}\label{IPM_Cond}
The matrices $A_1,\ldots,A_m$ are linearly independent in $\mathbb{S}^n$, and there exists a primal-dual feasible solution $(X,y,S)$ with $X,S\succ 0$.
\end{assumption}
\noindent
\SDO~lies at the intersection of convex
optimization, computational mathematics, and real algebraic geometry, with applications ranging from combinatorial optimization and control
to polynomial optimization and sums of squares~\cite{L15}. In this paper, our
focus is on the interaction between the semi-algebraicity of weighted central paths and the local convergence of interior-point methods (\IPM s).

\subsection{\IPM s and central path}
Analogous to linear optimization (\LO), primal-dual path-following \IPM s~\cite{AG03, NN94} are among the most efficient algorithms to solve \SDO~numerically. The core idea, originating from~\cite{FM90}, is to follow the \textit{central path}. In this paper, we consider the following formulation of the central path (see also Remark~\ref{other_formulations}): the central path of \SDO~is the map
\begin{align*}
\xi:(0,\infty)\longrightarrow
\mathbb{S}^n\times\mathbb{R}^m\times\mathbb{S}^n,
\qquad
\xi(\mu):=(X(\mu),y(\mu),S(\mu)),
\end{align*}
where for each fixed $\mu>0$, $(X(\mu),y(\mu),S(\mu))$ uniquely solves
\begin{equation}
\begin{aligned}\label{central_path_description}
\langle A_i, X \rangle &=b_i, \quad i=1,\ldots,m, \\ \sum_{i=1}^m y_i A_i + S -C&=0,\\
X S + SX - 2\mu I_n&=0, \ X,S \succ 0,
\end{aligned}
\end{equation}
in which the third equation is referred to as the centrality condition. It is easy to see that along the central path, the normalized complementarity gap coincides with the path parameter $\mu$. Under Assumption~\ref{IPM_Cond}, the central path exists and is real analytic; see~\cite[Theorems~3.1 and~3.3]{Kl02}. Further, $\xi$ is a \textit{semi-algebraic} function (see Section~\ref{semi-algebraic_def}) and converges to an optimal solution (see Proposition~\ref{weighted_central_path_convergence} with $W = I_n$ or~\cite{BM22,BM24}).
 
\vspace{5px}
\noindent
The machinery of primal-dual path-following \IPM s involves generating a trajectory of solutions that approximately solve~\eqref{central_path_description}, while taking $\mu \downarrow 0$, towards an approximate optimal solution of \SDO. The local behavior of primal-dual \IPM s is closely
tied to the regularity of the central path as $\mu\downarrow0$. Under \textit{strict complementarity} (see Definition~\ref{maximal_strict_comp}), the central path admits an analytic extension to its limiting point, and this regularity underlies classical super-linear and quadratic convergence results.
Without strict complementarity, however, the geometry changes qualitatively.  The central path need not be differentiable at $\mu=0$, and its local expansions may involve fractional powers of $\mu$. Thus, the ordinary barrier parameter does not in general
provide the higher-order derivatives needed to construct higher-order predictor directions.

The central idea of this paper is to overcome this obstruction by ramification. We prove that an entire family of \textit{weighted \AHO~central
paths} admits a common analytic reparametrization near the limiting point. The resulting higher-order derivatives are then used directly
to construct and analyze a higher-order predictor-corrector \IPM.

\section{Main results}\label{main_results}
Our results connect three components: the ramified regularity of weighted central paths, a higher-order predictor-corrector construction, and its convergence without strict complementarity.

\subsection{Ramification of weighted central paths}\label{semi-algebraic_weighted_central_paths} 
The core of our higher-order \IPM~is a semi-algebraic description of weighted central paths, which allows for their analytic reparametrization. For the purpose of our algorithmic development, we consider a \textit{feasible} parameterized family of \textit{weighted central paths} defined by functions 
\begin{align*}
\xi(\mu,W):=(X(\mu,W),y(\mu,W),S(\mu,W)):(0,\infty) \times \mathbb{S}^n_{++}  \to \mathbb{S}^{n} \times \mathbb{R}^m \times \mathbb{S}^{n},
\end{align*}
where for each fixed $\mu > 0$ and $W \succ 0$, $(X(\mu,W),y(\mu,W),S(\mu,W))$ uniquely solves 
\begin{equation}\label{Weighted-Central_Path}
\begin{aligned}
\langle A_i, X \rangle &=b_i, \quad i=1,\ldots,m, \\ \sum_{i=1}^m y_i A_i + S -C&=0,\\
X S + SX - 2\mu W&=0, \ X,S \succ 0,
\end{aligned}
\end{equation}
where $\mathbb{S}^n_{++}$ denotes the interior of the positive semi-definite cone. For a fixed $W \succ 0$, we denote the corresponding weighted central path by
\begin{align*}
\xi_W(\mu):=\xi(\mu,W), \quad X_W(\mu):=X(\mu,W), \quad y_W(\mu):=y(\mu,W), \quad S_W(\mu):=S(\mu,W),
\end{align*}
\noindent
and with a fixed $\mu > 0$ $\big(X_W(\mu),y_W(\mu),S_W(\mu)\big)$ is called a \emph{weighted central solution}. 

\vspace{5px}
\noindent
When $W=I_n$, a weighted central path is equivalent to the central path defined in~\eqref{central_path_description}. Under Assumption~\ref{IPM_Cond}, when $W \succ 0$, $\xi(\mu,W)$ exists and is analytic on $(0,\infty) \times \mathbb{S}^n_{++}$; see Proposition~\ref{prop:weighted_path_existence}. Further, $\xi(\mu,W)$ is jointly semi-algebraic, and for a fixed $W \succ 0$, it converges to an optimal solution as $\mu \to 0$; see Proposition~\ref{weighted_central_path_convergence}. More importantly, for the purposes of this paper, we establish the existence of a minimal positive integer $N$ such that the reparametrized path $\xi_W(\mu^{N})$ is analytic at $\mu = 0$ over a bounded subset 
\begin{align*}
\mathcal W_\beta:=
\{W\succ0\mid\|W-I_n\|_F\leq\beta\}, 
\end{align*}
where $0 < \beta < 1$ is a fixed parameter. This analytic reparametrization arises naturally from the \textit{Puiseux expansion} of weighted central paths, which we formally introduce in Section~\ref{Puiseux_Def}.
\begin{notation}
Throughout the paper, $\|\cdot\|_F$ denotes the Euclidean norm for vectors and the Frobenius norm for matrices.
\end{notation}
\begin{theorem}\label{global_ramification}
For every fixed $0<\beta<1$, there exists a minimal positive integer $N=N(\beta)$ such that, for every $W\in\mathcal W_\beta$, the reparametrized path $\xi_W(\mu^N)$ is analytic at $\mu=0$.
\end{theorem}
\noindent
Theorem~\ref{global_ramification} provides a common ramified reparametrization, but does not
assert joint analyticity in $(\mu,W)$ or uniform bounds on derivatives as $W$ varies. The latter requirement enters separately through the predictor regularity condition used in the convergence analysis of our higher-order \IPM.
\begin{remark}[Reparametrization exponent]
\label{rem:admissible_rho}
The integer $N$ in Theorem~\ref{global_ramification} is the smallest positive integer that provides a common analytic reparametrization of the weighted central paths.
We call a positive integer $\rho$ \textit{admissible} on $\mathcal W_\beta$ if $\xi_W(\mu^\rho)$ is analytic at $\mu=0$ for every $W\in\mathcal W_\beta$. In particular, $\rho=N$ is admissible, and every positive integer multiple of $N$ is also admissible. In what follows, we use $\rho$ to denote the reparametrization exponent employed by the higher-order \IPM.
\end{remark}

\begin{definition}
We refer to $N$ as the \textit{global ramification index} of weighted central paths on $\mathcal{W}_{\beta}$. 
\end{definition}

\begin{remark}
The proof of Theorem~\ref{global_ramification} can be modified to extend the result to all $W\succ0$. However, doing so would require a substantially longer and more cumbersome statement, which is unnecessary for the algorithmic development in this paper.
\end{remark}

\subsection{A higher-order \IPM}\label{main_results:HIPM}
We develop a higher-order feasible predictor-corrector \IPM~based on the \AHO~(Alizadeh-Haeberly-Overton~\cite{AHO98})~search directions. Our method exploits the analytic reparametrization $\mu \mapsto \mu^{\rho}$ of weighted central paths to construct higher-order predictor directions. Specifically, given an initial strictly feasible solution, the predictor directions are obtained from the $p$th-order Taylor approximation of the reparametrized weighted central path. We forward-track this path with respect to an outer neighborhood of the reparametrized central path. Upon reaching the boundary of the neighborhood, we apply centering to return to an inner neighborhood of the reparametrized central path.

\begin{remark}
\label{rem:nonadmissible_algorithm}
Our higher-order \IPM~does not require the use of an admissible reparametrization exponent $\rho$. Indeed, any positive integer $\rho$ for which the reparametrized weighted central paths possess derivatives through order $p$ can be used to construct the $p$th-order predictor. Theorem~\ref{global_ramification} guarantees
the existence of such a choice, since one may always take $\rho=N$, or any positive integer multiple of $N$.

Admissibility becomes important in our asymptotic convergence analysis.
In particular, Theorems~\ref{Algorithm_Convergence} and~\ref{convergence_of_vanishing_blocks} require an admissible reparametrization exponent together with suitable uniform predictor regularity of the
reparametrized weighted central paths.
\end{remark}
\noindent
We establish the super-linear convergence of the complementarity gaps $\la X_k,S_k \ra/n$ of Algorithm~\ref{alg:higher-order-IPM-imp}, without the strict complementarity condition. 

\begin{notation}
For the $k$th iterate of Algorithm~\ref{alg:higher-order-IPM-imp}, we denote the normalized complementarity gap by 
\begin{align*}
\mu_k:=\frac{\langle X_k,S_k\rangle}{n}.
\end{align*}
\end{notation}
\begin{theorem}[Super-linear convergence] \label{Algorithm_Convergence}
Let $0 < \beta < \frac{1}{4}$, $\rho$ be an admissible reparametrization exponent, and let
$p\geq \rho$. Suppose, in addition, that the family of reparametrized weighted
central paths is uniformly $(p,\rho)$-predictor regular on $\mathcal{W}_{\beta}$ (see Definition~\ref{predictor_regularity}). Then the normalized complementarity gaps
satisfy
\begin{align*}
\mu_{k+1}
=
O\!\left(
\mu_k^{(p+1)/\rho}
\right).
\end{align*}
Consequently, the sequence $\{\mu_k\}$ converges with Q-order at least
\begin{align*}
\frac{p+1}{\rho}>1.
\end{align*}
In particular, Algorithm~\ref{alg:higher-order-IPM-imp} converges
$Q$-super-linearly in the normalized complementarity gap.
\end{theorem}
\noindent
We also prove that the simultaneously vanishing blocks of $X_k$ and $S_k$, in the sense of the \textit{optimal partition} of \SDO~(see Definition~\ref{optimal_partition}), converge super-linearly, with a smaller order of convergence.
\begin{theorem}\label{convergence_of_vanishing_blocks}
Under the assumptions of Theorem~\ref{Algorithm_Convergence}, let $Q:=(Q_{\mathcal{B}},Q_{\mathcal{T}}, Q_{\mathcal{N}})$ be an orthonormal basis partitioned according to $\mathcal{B}$, $\mathcal{T}$, and $\mathcal{N}$. Then the components $Q^T_{\mathcal{N}} X_kQ_{\mathcal{N}}$ and $Q^T_{\mathcal{B}} S_kQ_{\mathcal{B}}$ converge super-linearly with order $(p+1)/\rho$. Further, if $p$ is large enough such that
\begin{align*}
p > \rho (2^{n-1}-1)-1,
\end{align*}
then the components  $Q^T_{\mathcal{T}} X_kQ_{\mathcal{T}}$ and $Q^T_{\mathcal{T}} S_kQ_{\mathcal{T}}$ converge super-linearly with order 
\begin{align*}
(p+1)/(\rho(2^{n-1}-1)).
\end{align*}
\end{theorem}
\noindent
Finally, we quantify the local iteration complexity of Algorithm~\ref{alg:higher-order-IPM-imp}.
Once the iterates enter the super-linear regime of Theorem~\ref{Algorithm_Convergence}, the remaining number of outer iterations required to attain a prescribed accuracy is doubly logarithmic in the target tolerance.

\begin{theorem}[Outer iteration complexity]
\label{iteration_complexity_result}
Let $0 < \beta < \frac{1}{4}$, $\rho$ be an admissible
reparametrization exponent, and let $p\geq\rho$. Suppose, in addition, that the family of reparametrized weighted central paths is uniformly $(p,\rho)$-predictor regular on $\mathcal{W}_{\beta}$ with constants $\tilde t>0$ and $c>0$. Assume that $0<t_0\leq\tilde t$, $\mu_0=t_0^\rho$, and set
\begin{align*}
q_0:=
\frac{1}{
1+
\left(
\frac{2\beta}
     {c\,t_0^d}
\right)^{1/\rho}
}
<1,
\end{align*}
where $d:=p+1-\rho \ge 1$. Then for every $\eta\in(0,\mu_0)$, Algorithm~\ref{alg:higher-order-IPM-imp} reaches
an iterate satisfying $\mu_k\leq\eta$ after at most
\begin{align*}
    K_g(\eta)
    :=
    \left\lceil
        \frac{\log(\mu_0/\eta)}
             {\rho\log(1/q_0)}
    \right\rceil
\end{align*}
outer iterations. Further, there exists $\mu_{\mathrm{loc}}>0$ such that, once $\mu_k\le\mu_{\mathrm{loc}}$, Algorithm~\ref{alg:higher-order-IPM-imp} computes a $\zeta$-optimal solution
in
\begin{align*}
O\!\left(
\log\log\frac{1}{\zeta}
\right)
\end{align*}
additional outer iterations. Consequently, taking $\eta=\mu_{\mathrm{loc}}$, for
$0<\zeta<\mu_{\mathrm{loc}}$,
Algorithm~\ref{alg:higher-order-IPM-imp} computes a $\zeta$-optimal solution in
\begin{align*}
   K_g(\mu_{\mathrm{loc}})
    +
    O\!\left(\log\log\frac{1}{\zeta}\right)
\end{align*}
outer iterations. 
\end{theorem}

\begin{remark}[Higher-order acceleration of a polynomial-time \IPM]
\label{rem:hybrid-ipm}
The higher-order method developed in this paper can also be used as a local acceleration phase for a standard polynomial-time predictor-corrector \IPM~for \SDO. For example, one may first apply a
polynomial-time method based on the \AHO, or more generally the Monteiro-Zhang~\cite{M98}, search directions until the normalized complementarity gap reaches a threshold $\bar\mu$ for which the local assumptions of Theorem~\ref{Algorithm_Convergence} hold, and then switch to Algorithm~\ref{alg:higher-order-IPM-imp}. If the standard phase reaches such a threshold in
\begin{align*}
O\!\left(
\sqrt n\log\frac{\mu_0}{\bar\mu}
\right)
\end{align*}
outer iterations, then the resulting two-phase method requires
\begin{align*}
O\!\left(
\sqrt n\log\frac{\mu_0}{\bar\mu}
+
\log\log\frac{1}{\zeta}
\right)
\end{align*}
outer iterations to attain a $\zeta$-optimal solution. This observation does not require modifying Algorithm~\ref{alg:higher-order-IPM-imp} itself. The higher-order \IPM~is used only after the iterates have entered its local super-linear regime.
\end{remark}
\noindent
The rest of this paper is organized as follows. In Section~\ref{main_techniques}, we briefly review the techniques used in the proofs. In Section~\ref{survey}, we review prior results on the analytic and algebraic properties of (weighted) central paths, super-linear convergence of \IPM s, and real algebraic geometry methods in optimization. In Section~\ref{AG_background}, we provide the preliminaries from real algebraic geometry, Puiseux series, and the central path that will be needed for the theoretical results of our higher-order \IPM. In Section~\ref{weighted_central_paths}, we extend the notion of analytic reparametrization to weighted central paths, and we also prove existence, convergence, analyticity, and semi-algebraicity of weighted central paths. We then prove Theorem~\ref{global_ramification}. In Section~\ref{higher-order_IPM}, we present our higher-order \IPM, and in Section~\ref{Proof_of_main_results}, we provide the proofs of the main theoretical results, i.e., Theorems~\ref{Algorithm_Convergence}-\ref{iteration_complexity_result}. In   Section~\ref{sec:numerical_demonstrations}, we summarize the elements of our implementation, and we present our numerical experiments. Finally, in Section~\ref{sec:conclusion}, we conclude the paper with directions for future research.

\subsection{Outline of the proofs}\label{main_techniques} In this section, we concisely review the ideas behind our proofs.

\vspace{5px}
\noindent
The proof of Theorem~\ref{global_ramification} utilizes tools from algorithmic real algebraic geometry~\cite{BPR06} and the theory of algebraic curves~\cite{W78}. More specifically, we apply parameterized bounded algebraic sampling~\cite[Algorithm~12.18]{BPR06}  to obtain a
parameterized real univariate representation containing the weighted central solutions. We then use cylindrical algebraic decomposition~\cite[Theorem~5.6]{BPR06} (see Definition~\ref{def:cad}) to partition
the parameter space so that the relevant real roots are delineable and their Thom encodings remain fixed on each cell. After specializing
the weight $W$, the Newton-Puiseux theorem (see Proposition~\ref{Puiseux_Newton_Thm}) shows that the corresponding coordinates admit Puiseux expansions whose ramification indices are bounded by the degree of the univariate representation. This yields
only finitely many possible ramification indices and hence a common reparametrization exponent for all $W\in\mathcal W_\beta$. 

\vspace{5px}
\noindent
For the proof of Theorem~\ref{Algorithm_Convergence}, we use the $p$th-order Taylor approximation
of the reparametrized weighted central paths (see~\eqref{Taylor_Approx_Central_Path} and~\eqref{HOSys}) together with the uniform $(p,\rho)$-predictor regularity condition to obtain a residual bound that is uniform over the relevant weights. When $\rho$ is admissible and $p\ge\rho$, this
bound yields a predictor parameter satisfying
\begin{align*}
t_+=O\!\left(t_k^{(p+1)/\rho}\right),
\end{align*}
and Lemma~\ref{corrector_gap} then yields $\mu_{k+1}=t_+^\rho$, which gives the $Q$-super-linear convergence asserted in Theorem~\ref{Algorithm_Convergence}.

\vspace{5px}
\noindent
For Theorem~\ref{convergence_of_vanishing_blocks}, we apply an error bound result for linear matrix inequalities from~\cite[Theorem~4.3]{MT19} to quantify the decay of the vanishing blocks associated with the optimal partition of \SDO. Combining these estimates with Theorem~\ref{Algorithm_Convergence} yields the stated block-wise super-linear
convergence rates.   

\vspace{5px}
\noindent
Finally, the proof of Theorem~\ref{iteration_complexity_result} has two parts. First, the quantitative predictor reduction of Proposition~\ref{prop:predictor-reduction} yields a uniform contraction
\begin{align*}
\mu_{k+1}\le q_0^\rho\mu_k
\end{align*}
throughout the predictor regularity region, and hence an explicit bound on the number of outer iterations required to reach any fixed local threshold. Once the iterates enter the super-linear regime,
Theorem~\ref{Algorithm_Convergence} gives
\begin{align*}
\mu_{k+1}\le C\mu_k^{(p+1)/\rho}.
\end{align*}
After rescaling the normalized complementarity gaps to absorb the multiplicative constant, this recurrence yields doubly exponential decay with respect to the local outer iteration counter, and therefore the bound
\begin{align*}
O\!\left(\log\log\frac{1}{\zeta}\right).
\end{align*}

\section{Background}\label{AG_background}
In this section, we recall some relevant notions from real algebraic geometry, the theory of algebraic curves, and the central path in accordance with~\cite{BM22,BM24, BPR06,BCR98, W78}. We use standard terminology from field theory and commutative algebra; see, e.g.,~\cite{Lang2002,Artin91}.  

\begin{notation}
In this paper, $\mathbb{R}$ and $\mathbb{C}$ denote the fields of real and complex numbers. Given a field $\mathbb{F}$, the ring of polynomials in $[X_1,\ldots,X_{k}]$ with coefficients in $\mathbb{F}$ is denoted by $\mathbb{F}[X_1,\ldots,X_{k}]$. Further, given a finite set $\mathcal{P}$ of polynomials in $\mathbb{F}[X_1,\ldots,X_k]$, the zero set of $\mathcal{P}$ in $\mathbb F^k$ is defined as
\begin{align*}
    \zero(\mathcal{P},\mathbb{F}^k):=\left\{x \in \mathbb{F}^k \mid \bigwedge_{P \in \mathcal{P}} P(x) = 0\right\}.
\end{align*}
\end{notation}

\subsection{Puiseux series}\label{Puiseux_Def}
In algebraic geometry, roots of a polynomial $F(X, Y) \in \mathbb{C}[X,Y]$ in terms of $X$ near $\mathbf{0}$ are expressed as Puiseux series. More specifically, in the proof of Theorem~\ref{global_ramification}, each coordinate of the central path or a weighted central path is described by a Puiseux series.

A \textit{Puiseux series} with coefficients in $\mathbb{C}$ (respectively  $\mathbb{R}$)
is a Laurent series with fractional exponents, i.e., a series of the form $\sum_{i = r}^{\infty} c_i \varepsilon^{i/q}$ where $c_i \in \mathbb{C}$ (resp. $c_i \in \mathbb{R}$), $i,r \in \mathbb{Z}$, and $q$ is a minimal positive integer, which is called the \textit{ramification index} of the Puiseux series. The set of Puiseux series in $\varepsilon$ with coefficients in $\mathbb{C}$ (resp. $\mathbb{R}$) is a field, which we denote by 
$\mathbb{C}\langle \langle \varepsilon \rangle\rangle$ (resp. $\mathbb{R}\langle \langle \varepsilon \rangle\rangle$). We let $o(\cdot)$ denote the order of a Puiseux series, and it is defined as $o(\sum_{i = r}^{\infty} c_i \varepsilon^{i/q})=r/q$ if $c_{r} \neq 0$, see~\cite[Section~2.6]{BPR06}. 

\vspace{5px}
\noindent
If $F(0,0)=0$ and $\partial F/\partial Y(0,0) \neq 0$, then the Implicit Function Theorem~\cite[Page~19]{GH94} provides a unique analytic function $Y=\psi(X)$ defined in a neighborhood of $X=0$, such that $\psi(0)=0$ and
\begin{align*}
F(X,\psi(X))=0.
\end{align*} 
If $\partial F/\partial Y(0,0)=0$, the Implicit Function Theorem no longer guarantees such an analytic parametrization, and fractional-power expansions may occur. In general, the Newton-Puiseux theorem~\cite[Theorem~3.1 of Chapter IV]{W78} shows that a Puiseux expansion of $F(X, Y) = 0$ always exists. 
\begin{proposition}[Theorems 3.2 and~4.1 and Section~4.2 of Chapter IV in~\cite{W78}]\label{Puiseux_Newton_Thm}
Let 
\begin{align*}
F(X,Y)=a_0 + a_1Y + \cdots + a_d Y^d \in \mathbb{C}[X][Y],
\end{align*}
where $a_d \neq 0$. There exist $d$ (not necessarily distinct) Puiseux series $\psi_i(X) \in \mathbb{C} \langle \langle X \rangle \rangle$ for $i=1,\ldots,d$ such that 
\begin{align*}
F(X,Y)=a_d \prod_{i=1}^d (Y - \psi_i).
\end{align*}  
\end{proposition} 
\noindent
A subfield of $\mathbb{R}\langle \langle \varepsilon \rangle\rangle$, denoted by $\mathbb{R}\langle \varepsilon\rangle$, can be identified with the field of germs at $0^+$ of continuous semi-algebraic functions~\cite[Theorem~3.14]{BPR06}. More precisely, the field $\mathbb{R}\langle \varepsilon\rangle$ consists of the elements of $\mathbb{R}\langle \langle \varepsilon \rangle\rangle$ that are algebraic over $\mathbb{R}(\varepsilon)$, the field of rational functions in the indeterminate $\varepsilon$ with coefficients in $\mathbb R$~\cite[Page~370]{Artin91}, and it is called the field of \emph{algebraic Puiseux series} with coefficients in $\mathbb{R}$. It is the real closure~\cite[Page~43]{BPR06} of the ordered field $\mathbb{R}(\varepsilon)$ in which $\varepsilon$ is positive and infinitesimal, i.e.,
\begin{align*}
0<\varepsilon<r
\qquad\text{for every }r\in\mathbb{R}_{>0}.
\end{align*}
Thus, every element of $\mathbb{R}\langle \varepsilon\rangle$ is represented, for sufficiently small $t_0>0$, by a continuous semi-algebraic function
\begin{align*}
(0,t_0)\longrightarrow\mathbb{R}.
\end{align*}
This interpretation makes $\mathbb{R}\langle \varepsilon\rangle$ particularly useful for studying germs of semi-algebraic functions, including the germs of weighted central paths at $\mu=0$.

\subsection{Semi-algebraic sets}\label{semi-algebraic_def}
Given the real field $\mathbb{R}$, a \textit{quantifier-free formula} $\mathrm{\Phi}(X_1,\ldots,X_{k})$ with coefficients in $\mathbb{R}$ is a Boolean combination of atoms $P > 0$, $P=0$, or $P < 0$ where $P \in \mathbb{R}[X_1,\ldots,X_{k}]$, and $\{X_1,\ldots,X_{k}\}$ are called the \emph{free variables} of $\mathrm{\Phi}$. 
A quantified formula is given by 
\begin{align*}
\Psi= (\mathrm{Q}_1 Y_1) \cdots (\mathrm{Q}_{\ell} Y_{\ell}) \ \Phi(X_1,\ldots, X_{k},Y_1,\ldots,Y_{\ell}),
\end{align*}
in which $\mathrm{Q}_i \in\{\forall,\exists\}$ are quantifiers and $\Phi$ is a quantifier-free formula with polynomials in
\begin{align*}
\mathbb{R}[X_1,\ldots,X_{k},Y_1,\ldots,Y_{\ell}].
\end{align*}
The set of all $(x_1,\ldots,x_{k}) \in \mathbb{R}^{k}$ satisfying $\Psi$ is called the $\mathbb{R}$-realization of $\Psi$. A \textit{semi-algebraic} subset of $\mathbb{R}^{k}$ is defined as the $\mathbb{R}$-realization of a quantifier-free formula. Alternatively, a semi-algebraic subset of $\mathbb{R}^{k}$ is a subset of the form
\begin{align*}
\bigcup_{i=1}^s \bigcap_{j=1}^{r_i}\left\{x \in \mathbb{R}^{k} \mid P_{ij}(x) \substack{<\\=}  \ 0 \right\},    
\end{align*}
where $P_{ij} \in \mathbb{R}[X_1,\ldots,X_{k}]$. Given a semi-algebraic set $S \subset \mathbb{R}^k$, a function $f:S \to \mathbb{R}$ is called semi-algebraic if its graph is a semi-algebraic subset of $\mathbb{R}^{k+1}$. 

\vspace{5px}
\noindent
The family of semi-algebraic subsets of $\mathbb{R}^{k}$ is closed under finite union, finite intersection, and complementation. Further, it is a classical result (see~\cite{Tarski51}) that every quantified formula is equivalent to a quantifier-free formula (see~\cite[Theorem~14.16]{BPR06} and~\cite[Algorithm~14.5]{BPR06}).

\subsection{Sign conditions, univariate representations, and Thom encodings}\label{sign_cond_Thom_encod}
In the proof of Theorem~\ref{global_ramification}, we describe weighted central points symbolically, using their \textit{real univariate representations}. We follow the representation used in~\cite{BPR06}. 

\vspace{5px}
\noindent
Given a finite family $\mathcal{P} \subset \mathbb{R}[X_1,\ldots,X_{k}]$, a \textit{sign condition} on $\mathcal{P}$ is an element of $\{-1,0,1\}^{\mathcal{P}}$, i.e., a mapping $\mathcal{P} \to \{-1,0,1\}$. Let $P \in \mathbb{R}[X]$ and 
\begin{align*}
\Der(P):= \big\{P, \frac{d P}{d X},\frac{d^2 P}{d X^2},\ldots,\frac{d^{(\deg(P))} P}{d X^{(\deg(P))}}\big\}
\end{align*}
denote a list of polynomials in which $d^i P/d X^i$ for $i > 0$ is the formal $i$th-order derivative of $P$ and $\deg(P)$ stands for the degree of $P$. A \textit{Thom encoding} $\sigma$ of real zero $t_{\sigma}$ of $P$ is a sign condition on $\Der(P)$ such that $\sigma(P)=0$. The following proposition~\cite[Proposition~2.27]{BPR06} indicates that for any
$P \in \mathbb{R}[X]$ and a root $x \in \mathbb{R}$ of $P$, its Thom encoding characterizes $x$. 

\begin{proposition} [Thom's Lemma] 
\label{prop:Thom}
Let $P \in \mathbb{R}[X]$ be a univariate polynomial and $\sigma \in \{-1,0,1\}^{\Der(P)}$. Then the realization of the sign condition $\sigma$ is either empty, a point, or an open interval.
\end{proposition} 
\noindent
In particular, when \(\sigma(P)=0\), its realization is either empty or a single point.

\vspace{5px}
\noindent
A $k$-univariate representation is a $(k+2)$-tuple of polynomials $u=\big(f, g_{0},\ldots,g_{k}\big) \in \mathbb{R}[T]^{k+2}$, where $f$ and $g_0$ are coprime. A real $k$-univariate representation of an $x \in \mathbb{R}^{k}$ is a pair $(u,\sigma)$ of a $k$-univariate representation $u$ and a Thom encoding $\sigma$ of a real root $t_{\sigma}$ of $f$ such that
\begin{align*}
x=\bigg(\frac{g_1(t_{\sigma})}{g_0(t_{\sigma})},\dots,\frac{g_{k}(t_{\sigma})}{g_0(t_{\sigma})}\bigg) \in \mathbb{R}^{k}.
\end{align*}
\noindent
Given a finite family $\mathcal{P} \subset \mathbb{R}[X_1,\ldots,X_k]$ with $\zero(\mathcal{P},\mathbb{R}^k)$ being non-empty and bounded, there exists an algorithm~\cite[Algorithm~12.16]{BPR06} to generate real $k$-univariate representations describing a set of points that meets every semi-algebraically connected component of $\zero(\mathcal{P},\mathbb{R}^k)$. A parametric version of this algorithm~\cite[Algorithm~12.18]{BPR06} will be used in the proof of Theorem~\ref{global_ramification}.

\subsection{Cylindrical algebraic decomposition}\label{cad}
The notion of cylindrical algebraic decomposition~\cite{BCR98,BPR06} plays a key role in semi-algebraic geometry, and will be needed later in Section~\ref{weighted_central_paths} for the proof of Theorem~\ref{global_ramification}. In what follows, we include the definition from~\cite[Definition~5.1]{BPR06}.
\begin{definition}[Cell decomposition]
\label{def:cad}
Fixing the standard basis of $\mathbb{R}^k$, we identify for each $i,1 \leq i \leq k$, $\mathbb{R}^i$ with the span of the first
$i$ basis vectors. A cylindrical algebraic decomposition (or simply a cell decomposition) of $\mathbb{R}$ is a $1$-tuple $(\mathcal{D}_1)$, where
$\mathcal{D}_1$ is a finite set of subsets of $\mathbb{R}$, each element being a point or an open interval, which together gives a partition of $\mathbb{R}$. 
A cell decomposition of $\mathbb{R}^k$ is a $k$-tuple $(\mathcal{D}_1,\ldots,\mathcal{D}_k)$,
where each $\mathcal{D}_i$ is a decomposition of $\mathbb{R}^i$,
$(\mathcal{D}_1,\ldots, \mathcal{D}_{k-1})$ is a cell decomposition of $\mathbb{R}^{k-1}$,
and $\mathcal{D}_k$ is a finite set of semi-algebraic subsets of $\mathbb{R}^k$ (called the cells of $\mathcal{D}_k$) giving a partition of $\mathbb{R}^k$ consisting 
of the following: for each $C \in \mathcal{D}_{k-1}$, there is a finite set of semi-algebraic continuous functions
$f_{C,1}, \ldots, f_{C,N_C}: C \rightarrow \mathbb{R}$ such that $f_{C_1} < \cdots < f_{C,N_C}$,
and each element of $\mathcal{D}_k$ is either the graph of a function $f_{C,i}$ 
or of the form 
\begin{enumerate}[(a)]
    \item $\{ (x,t) \;\mid\; x\in C, \ t < f_{C,1}(x) \}$, 
    
    \item  $\{ (x,t) \;\mid\; x\in C, \ f_{C,i}(x) < t < f_{C,i+1}(x) \}$, 
    
    \item  $\{ (x,t) \;\mid\; x\in C, \ f_{C,N_C}(x) < t  \}$,     
    \item  $\{ (x,t) \;\mid\; x\in C \}$
\end{enumerate}
(the last case arising is if the set of functions $\{f_{C,i} | 1 \leq i \leq N_C \}$ is empty).
\end{definition}
\noindent
We will say that the cell decomposition $(\mathcal{D}_1,\ldots,\mathcal{D}_k)$ is \textit{adapted} to 
a semi-algebraic subset $D$ of $\mathbb{R}^k$, if for each $C \in \mathcal{D}_k$, $C \cap D$ is either equal to $C$ or empty. It is well-known that for every semi-algebraic subset $D$ of $\mathbb{R}^k$ there exists a cell decomposition adapted to $D$~\cite[Theorem~5.6]{BPR06}. 

\vspace{5px}
\noindent
Given a finite family of polynomials
$\mathcal P\subset \mathbb{R}[X_1,\ldots,X_k]$, a cylindrical
algebraic decomposition is said to be \emph{adapted to $\mathcal P$}
if every $P\in\mathcal P$ has constant sign on each cell of the
decomposition. In particular, let
$P\in\mathbb R[X_1,\ldots,X_k,T]$ and define
\begin{align*}
\Der_T(P)
:=
\left\{
P,\frac{\partial P}{\partial T},\ldots,
\frac{\partial^{d}P}{\partial T^{d}}
\right\},
\qquad d=\deg_T(P).
\end{align*}
A cylindrical algebraic decomposition adapted to $\Der_T(P)$ can be chosen so that, over each cell of the induced decomposition of the parameter space
$\mathbb R^k$, the real roots of $P(x,T)$ are delineable.
Consequently, the number and ordering of the real roots, as well as
their Thom encodings, remain fixed as $x$ varies within a cell; see
\cite[Chapter~5]{BPR06}.

\subsection{Analytic reparametrization of the central path}\label{CP_SDO}
We first recall the strict complementarity condition.

\begin{definition}\label{maximal_strict_comp}
A primal-dual optimal solution
$(X^*,y^*,S^*)\in\mathcal P^*\times\mathcal D^*$ is called
\emph{strictly complementary} if
\begin{align*}
X^*+S^*\succ0.
\end{align*}
The strict complementarity condition holds if such an optimal solution
exists.
\end{definition}
\noindent
Since a central solution, for every positive $\mu$, is an isolated non-singular solution of~\eqref{central_path_description}~\cite[Theorem~3.3]{Kl02}, all higher-order derivatives of the central path exist, by the Implicit Function Theorem (see e.g.,~\cite[Corollary~2.9.8]{BCR98}). More explicitly, the  $j$th-order derivative $(X^{(j)}, y^{(j)}, S^{(j)})$ can be computed by solving 
\begin{align*}
\langle A_i, X^{(j)} \rangle &=0,\\[-1\jot]
\sum_{i=1}^m y_i^{(j)} A_i  + S^{(j)} &=0,\\[-1\jot]
X^{(j)} S(\mu)+X(\mu) S^{(j)} + S^{(j)} X(\mu) + S(\mu) X^{(j)}&=\\
\begin{cases} 2I_n \ \ &j=1,\\
-\sum_{k=1}^{j-1} {j \choose k}\big(X^{(k)} S^{(j-k)} + S^{(k)} X^{(j-k)}\big), \ \ &j>1.
\end{cases}
\end{align*} 
Under the strict complementarity condition, the central path admits an
analytic extension to $\mu=0$~\cite[Theorem~1]{H02}. Consequently, its
derivatives of every order have finite limits as $\mu\downarrow 0$. Without strict complementarity, this endpoint regularity may fail; in particular, the first-order derivative of the central path can become unbounded as $\mu\downarrow 0$~\cite{GS98}. The following example, adopted from~\cite[Example~4]{BM24}, illustrates this phenomenon.

\begin{example}\label{derivatives_diverge}
Consider the minimization of a linear objective function over the \emph{3-elliptope}, see Figure~\ref{3elliptope_alg_curve}, as follows:
\begin{align*}
\inf_{(x,y,z)}\bigg\{4x-4y-2z \mid 
\begingroup 
\setlength{\tabcolsep}{.75pt} 
\renewcommand{\arraystretch}{.8}
\begin{pmatrix} 1 & \ x & \ y\\x & \ 1 & \ z\\y & \ z & \ 1 \end{pmatrix}
\endgroup \succeq 0\bigg\},
\end{align*} 
\noindent
which has a unique non-strictly complementary optimal solution 
\begin{align*}
\setlength{\arraycolsep}{1pt}
X^*=\begin{pmatrix} \ \ 1 & \ -1 & \ \ \ 1 \\ -1 & \ \ \ 1 & \ -1\\ \ \ 1 & \ -1 & \ \ \ 1 \end{pmatrix}, \ \
y^*=(-4, \ -1, \ -1)^T,\ \
S^*=\begin{pmatrix} \ \ 4 & \ \ \ 2 & \ -2 \\ \ \ 2 & \ \ \ 1 & \ -1\\ -2 & \ -1 & \ \ \ 1 \end{pmatrix}.
\end{align*}
The central path equations can be equivalently written as 
\begin{align*}
\underbrace{\begin{pmatrix} 1  & X_{12} & X_{13}\\ X_{12}  & 1 & X_{23}\\X_{13} & X_{23} & 1 \end{pmatrix}}_{X(\mu)} \underbrace{\begin{pmatrix} -y_1 & 2 & -2 \\2 & -y_2 & -1\\ -2 & -1 & -y_3 \end{pmatrix}}_{S(\mu)} =\begin{pmatrix} \mu & 0 & 0 \\0 & \mu & 0\\ 0 & 0 & \mu \end{pmatrix},
\end{align*}
\noindent
which, after eliminating $y_1$, $y_2$, and $y_3$, yields
\begin{equation}\label{central_path_after_elimination}
\begin{aligned}
 \begin{cases} F(\mu,X_{12}):=2X_{12}^3+(2-\frac{\mu}{2})X_{12}^2-(\mu+2)X_{12}-2=0,\\X_{13} = -X_{12},\\X_{23} = -2X^2_{12}+\frac{\mu}{2} X_{12}+1.   
\end{cases}
\end{aligned} 
\end{equation}
\noindent
By Proposition~\ref{Puiseux_Newton_Thm}, the first equation in~\eqref{central_path_after_elimination} has a root corresponding to $X_{12}(\mu)$ that can be expressed as the Puiseux series
\begin{align}\label{coordinate_Puiseux_expansion}
X_{12}(\mu)&=-1+\frac{\sqrt{8}}{8} \mu^{\frac12}+\frac{1}{32} \mu -\frac{11\sqrt{8}}{2048} \mu^{\frac32}-\frac{3}{512} \mu^{2}-\frac{121\sqrt{8}}{1048576} \mu^{\frac52}\\
&\quad+\frac{15}{32768} \mu^3+\frac{19405\sqrt{8}}{268435456}\mu^{\frac72}+ \cdots \nonumber
\end{align}
with ramification index $2$. All this implies the non-Lipschitzian convergence 
\begin{align*}
\|X(\mu)-X^{*}\| = O(\sqrt{\mu}) \ \ \text{and} \ \ \|S(\mu)-S^{*}\| = O(\sqrt{\mu}).
\end{align*}
It is clear from~\eqref{coordinate_Puiseux_expansion} that the first-order derivative of $X_{12}(\mu)$ does not exist at $\mu = 0$.
\begin{figure}
\center
\begin{minipage}[c]{0.5\textwidth}
\includegraphics[height=2.2in]{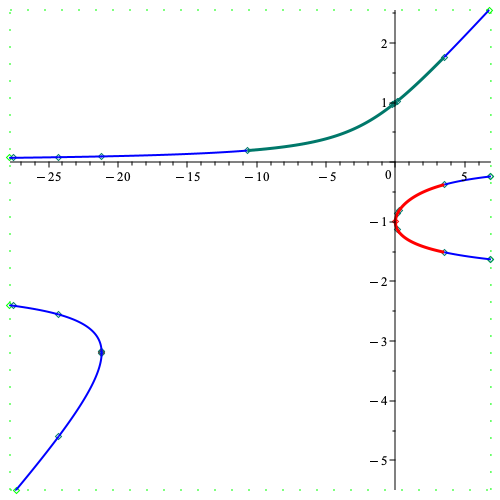}
\end{minipage}
\begin{minipage}[c]{0.4\textwidth}
\includegraphics[height=1.5in]{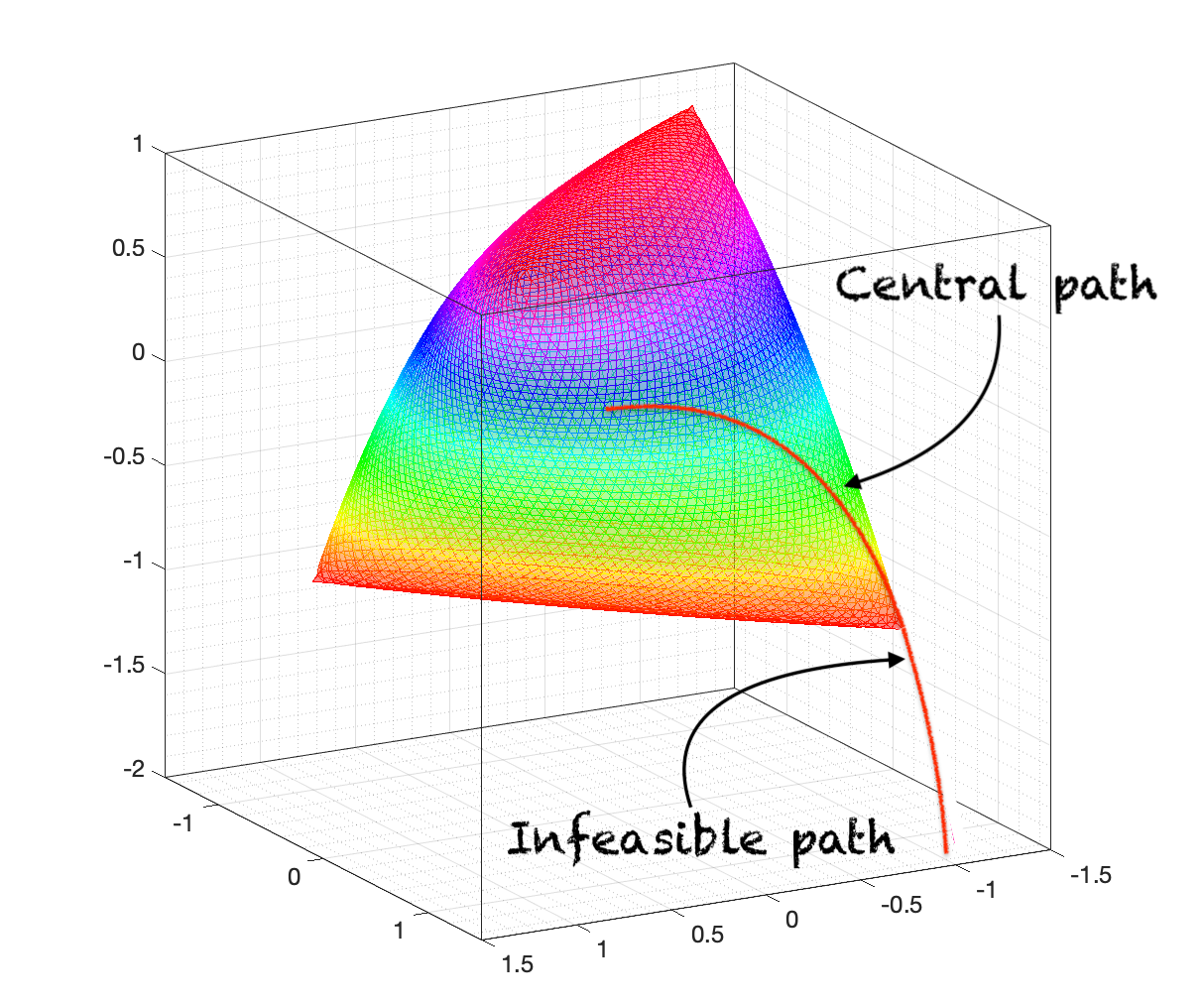}
\end{minipage}
\caption{The central path is not analytic at $\mu = 0$.}
\label{3elliptope_alg_curve}
\end{figure}
However, the ramification index of $X_{12}(\mu)$ suggests the local parametrization $(\mu^{2},\Phi(\mu))$ of the curve $F = 0$ around $\mu = 0$, where 
\begin{align*}
\Phi(\mu)=-1+\frac{\sqrt{8}}{8} \mu+\frac{1}{32} \mu^2 -\frac{11\sqrt{8}}{2048} \mu^{3}-\frac{3}{512} \mu^{4}-\frac{121\sqrt{8}}{1048576} \mu^{5}+\frac{15}{32768} \mu^6+ \ldots, 
\end{align*}
see Figure~\ref{3elliptope_alg_curve}. Thus, a reparametrization $\mu \mapsto \mu^{\rho}$ recovers the analyticity at $\mu = 0$, where $\rho$ is a positive integer multiple of the ramification index of the Puiseux expansion. 
\end{example}

\vspace{5px}
\noindent
The preceding example illustrates the general mechanism behind analytic reparametrization of the central path. This mechanism follows from the semi-algebraicity of the central path and the Newton-Puiseux theorem, as we briefly summarize below.

\begin{notation}\label{matrix_identification}
We identify a symmetric matrix $X$ by a vector using the linear transformation 
\begin{align*}
\svectorize:\mathbb{S}^{n} \to \mathbb{R}^{t(n)}, \quad t(n):={n+1 \choose 2},
\end{align*}
which multiplies the off-diagonal entries of a symmetric matrix by $\sqrt{2}$, and stacks the upper triangular part into a vector, i.e.,
\begin{align}
\svectorize(X):=\big(X_{11}, \sqrt{2}X_{12},\ldots, \sqrt{2}X_{1n}, X_{22}, \sqrt{2} X_{23},\ldots,\sqrt{2}X_{2n},\ldots, X_{nn}\big)^T. \label{linear_transformation_svec}
\end{align}
\end{notation}

\begin{notation}\label{vector_form}
We use the notation $(\cdot;\cdot;\cdots;\cdot)$ for the concatenation of column vectors. Accordingly, using the identification~\eqref{linear_transformation_svec}, we identify the central path by $v(\mu):=(x(\mu);y(\mu);s(\mu))$. Furthermore, the coordinates of the central path are denoted by $v_i(\mu)$ for $i=1,\ldots,\bar{n}$, where 
\begin{align*}
\bar{n}:=m+2t(n).
\end{align*}
This notation will also be adopted for weighted central solutions.
\end{notation}

\begin{notation}\label{Intermediates_Notation}
Using the identification of Notation~\ref{vector_form}, the indeterminates of the polynomials defining~\eqref{central_path_description} are denoted by $V_1,\ldots, V_{\bar{n}}$.
\end{notation}

\subsubsection{Puiseux reparametrization of the central path.}
The central path is a semi-algebraic function (see Proposition~\ref{weighted_central_path_convergence}).
Hence, for each coordinate $v_i(\mu)$, the germ of $v_i$ at $0^+$
defines an element of the field $\mathbb{R}\langle \mu\rangle$ of
algebraic Puiseux series. In particular, $v_i$ is algebraic over
$\mathbb{R}(\mu)$. Therefore, there exists a non-zero polynomial
\begin{align*}
P_i\in\mathbb{R}[\mu,V_i]
\end{align*}
such that
\begin{align*}
P_i(\mu,v_i(\mu))=0,
\qquad i=1,\ldots,\bar n,
\end{align*}
for all sufficiently small $\mu>0$. By Proposition~\ref{Puiseux_Newton_Thm}, the coordinate $v_i(\mu)$ admits a convergent
Puiseux expansion of the form
\begin{align*}
v_i(\mu)=\sum_{j=r_i}^{\infty} c_{ij}\mu^{j/q_i},
\end{align*}
where $q_i$ is its ramification index. Since the central path converges
as $\mu\downarrow 0$ (see Proposition~\ref{weighted_central_path_convergence}), each $v_i(\mu)$ has a finite limit at $\mu=0$,
and therefore $r_i\geq 0$. Thus,
\begin{align*}
v_i(\mu)=\sum_{j=0}^{\infty} c_{ij}\mu^{j/q_i}.
\end{align*}

Let $q$ be the least common multiple of the ramification indices
$q_1,\ldots,q_{\bar n}$. Then, for every positive integer $\rho$
divisible by $q$, each coordinate
\begin{align*}
v_i(\mu^\rho)
\end{align*}
is represented by a convergent power series in $\mu$ in a neighborhood
of $\mu=0$. Hence the reparametrized central path
$v(\mu^\rho)$ admits an analytic extension to $\mu=0$. In particular, under strict complementarity no ramification is needed: the analyticity of the central path at $\mu=0$ implies that all coordinate ramification indices are equal to one, and hence $q=1$.

\section{Prior and related work}\label{survey}
\subsection{Central path}
The endpoint regularity discussed in Section~\ref{CP_SDO} has been studied
extensively in the literature. Under strict complementarity, analyticity of the central path of \SDO~at $\mu=0$ was established in~\cite{GS98, H02, PS04}. This yields, in particular, the Lipschitzian estimates
\begin{align*}
\|X(\mu)-X^{**}\|_F=O(\mu),
\qquad
\|S(\mu)-S^{**}\|_F=O(\mu),
\end{align*}
where $(X^{**},y^{**},S^{**})$ denotes the limit point of the central
path. Without strict complementarity, substantially more complicated
asymptotic behavior may occur, including fractional powers of $\mu$.

\vspace{5px}
\noindent
The failure of analyticity has been extensively studied for \LO~and linear complementarity problems (\LCP). Reparametrization provides a way to recover analyticity in this
degenerate setting. For \LCP~problems, analytic
extension can be recovered by the reparametrization $\mu \mapsto \mu^2$~\cite{SW98, SW99, SWM98}. This behavior is largely explained by the fact that variables along the central path have magnitudes of order $O(1)$, $O(\mu)$, or $O(\sqrt{\mu})$~\cite{IPRT00}. For \SDO, such a classification holds only under the strict complementarity condition. In its absence, however, solutions along the central path may exhibit blocks of order $O(\mu^d)$, where the exponent $d$ can be exponentially large (see~\cite{MT19}). In~\cite{H02}, the author mentions that
\textit{``... some examples of SDPs show that the variables along the central path for SDP display much more
complex behavior than those for LCP. A question now arises as to how the results of this paper may be used to construct higher-order methods for SDP ... the question that still remains open is how to define paths for SDP, so that their analyticity would serve for the construction of
algorithms similarly as the analyticity of the weighted paths serves for LCP."}

For general \SDO, Basu and Mohammad-Nezhad~\cite{BM24} showed that the central path admits an analytic reparametrization $\mu\mapsto\mu^\rho$ at $\mu=0$. They characterized the minimal such exponent as the least common multiple of the ramification indices of the Puiseux expansions of the coordinate functions of the central path, and established the bound
\begin{align*}
\rho \leq 2^{O(m^2+n^2m+n^4)}.
\end{align*}
In particular, the required ramification exponent may grow exponentially with the problem dimension.

 \subsection{Weighted central paths and regularity}
\label{sec:weighted-paths-literature}
Weighted central paths had already played an important role in the analysis of \LO~and \LCP, where their structure is considerably simpler (see e.g.,~\cite{G94,SWM98,MT96}).  Their extension to \SDO~is more subtle because the notion of weighted centrality depends on the chosen symmetrization of the matrix complementarity equation. Instead, weighted paths in \SDO~are typically defined directly through a symmetrized complementarity equation
\begin{align*}
\Phi(X,S)=\mu W, \qquad W\succ0,
\end{align*}
where the choice of the interior-point map $\Phi$ determines the corresponding family of weighted paths~\cite{MP98,MZ00}. In particular, Monteiro and Pang~\cite{MP98} established existence and uniqueness properties for weighted trajectories associated with the \AHO~map
\begin{align*}
        \Phi_{\text{\AHO}}(X,S)=\frac{XS+SX}{2},
\end{align*}
and Monteiro and Zanj\'acomo~\cite{MZ00} subsequently developed a general framework for weighted paths associated with a broad class of interior-point maps. Alternative notions of weighted centers for \SDO~have also been studied; see, e.g., Sturm and Zhang~\cite{SZ00}.

\begin{remark}[Comparison with previous weighted path regularity results]
The limiting regularity of weighted central paths has been studied in
several settings. Under strict complementarity, Prei{\ss} and Stoer~\cite{PS04} proved, for \AHO~weighted paths of monotone
semi-definite \LCP s, analytic dependence jointly on the path parameters, including the barrier parameter and the weight matrix. For \SDO, Lu and
Monteiro~\cite{LM07} established analyticity at $\mu=0$, in the original barrier parameter, for each fixed $W \succ 0$, together with results on the limiting behavior and derivatives of
the corresponding \AHO~weighted path. Related regularity results for other symmetrizations have also been obtained; see, e.g.,~\cite{LM05}.

Particularly relevant to Theorem~\ref{global_ramification} is
\cite[Proposition~3.5]{LM07}. Without assuming strict complementarity, Lu and Monteiro state that, for each fixed $W \succ 0$, there exists an analytic scalar reparametrization $\nu=\nu(t)$, with $\nu(0)=0$, such that the
corresponding \AHO~weighted path
\begin{align*}
t\longmapsto
\bigl(X(\nu(t),W),y(\nu(t),W),S(\nu(t),W)\bigr)
\end{align*}
admits an analytic extension to $t=0$. The reparametrization is thus obtained separately for each fixed $W$ and is not further
characterized.

Theorem~\ref{global_ramification} proceeds in a different direction. Without assuming strict
complementarity, we show that the reparametrization can be chosen in
the explicit power form $\mu\mapsto\mu^N$, where a single positive integer $N$ works simultaneously for every $W\in\mathcal W_\beta$.
Thus, in contrast with the fixed-$W$ reparametrization above, the same
ramification exponent applies to an entire family of weighted \AHO~paths. On the other hand, Theorem~\ref{global_ramification} does not assert joint analyticity in $(\mu,W)$ or uniform bounds on the derivatives as $W$ varies. The
latter uniformity is imposed separately through the predictor regularity condition in Definition~\ref{predictor_regularity}.
\end{remark}

\vspace{5px}
\noindent
The regularity of weighted central paths is closely connected with the
local convergence analysis of primal-dual \IPM s; in particular, the role
played by the weight and by additional asymptotic centering conditions depends on the chosen symmetrization. We discuss these issues next.

\subsection{Local convergence of primal-dual \IPM s}\label{local_convergence_IPMs}
Classical local convergence analyses of primal-dual \IPM s are based on first-order search directions computed in neighborhoods of the central path. Under suitable regularity conditions, super-linear or quadratic convergence has been established for broad classes of primal-dual \IPM s, including \SDO
~\cite{AHO98, JPS99, JPS99b, KSS98, LSZ98, PS98, PR98, PS96}, semi-definite \LCP s~\cite{KSS98, KSS99, LM05b, S11}, and \LCP s~\cite{JPH95, JPTZ91, Mi96, MW94, P01, P04, PS98b, PS97, PS10, S99}.

\vspace{5px}
\noindent
For \SDO, the local convergence analysis is considerably more delicate.
Classical super-linear convergence results generally assume strict
complementarity and exploit the optimal partition. However, the strict complementarity is not always the only
hypothesis. For example, Kojima et al.~\cite{KSS98} established super-linear convergence for a predictor-corrector method under strict complementarity, non-degeneracy, and an additional tangential condition
requiring the iterates to approach the central path through a shrinking neighborhood. Potra and Sheng~\cite{PS98} subsequently showed that non-degeneracy can be removed, while retaining strict complementarity together with
the shrinking-neighborhood condition, and later weakened the latter
condition. Luo et al.~\cite{LSZ98} similarly established super-linear convergence of a symmetric primal-dual path-following method under strict complementarity together with asymptotic centering, and proposed a modification that enforces this behavior. The role of the additional centering conditions depends on the search
direction. Kojima et al.~\cite{KSS99} proposed an \AHO-based predictor-corrector method and claimed local quadratic convergence under strict complementarity. Lu and
Monteiro~\cite{LM05b} subsequently pointed out an error in that proof and showed that the original result is valid under an additional non-degeneracy assumption. They also introduced a slight
modification of the corrector step and established local quadratic convergence of the resulting \AHO-based method under strict complementarity, without imposing a separate tangential-convergence condition.

\vspace{5px}
\noindent
These additional asymptotic conditions can also be interpreted in terms of weighted central paths, and their role depends on the chosen symmetrization. For the scaling $X^{1/2}SX^{1/2}$, one classical sufficient condition is $W_k\to I$~\cite{KSS98}, while weaker formulations (e.g.,~\cite{PS98}) require only the appropriate off-diagonal block of $W_k$ to vanish asymptotically~\cite{LM07}.

The situation in the present paper is different. Our higher-order predictor is constructed from derivatives of weighted \AHO~paths. Because the weight varies between iterations, our local analysis requires uniform control of the higher-order predictors; this issue is discussed further in Remark~\ref{distinction_between_SC_regualrity}.
 
\vspace{5px}
\noindent
More recent work continues to rely on strict complementarity, sometimes together with additional assumptions on the problem, initialization, or the behavior of the iterates. For example, implementable polynomial-time methods~\cite{S25} with super-linear convergence have been obtained under strict complementarity by imposing non-degeneracy or neighborhood-narrowing
conditions, or by considering particular homogeneous or feasibility formulations. To the best of our knowledge, however, super-linear convergence of a primal-dual \IPM~for general \SDO~has not been established without strict complementarity.

\begin{remark}\label{distinction_between_SC_regualrity}
The distinction between strict complementarity and the additional
regularity conditions used in super-linear convergence analysis is important for the present paper. Strict complementarity concerns the limiting primal-dual geometry and, in classical \SDO~analyses, provides
strong regularity of the central or weighted central paths. Additional conditions such as non-degeneracy, shrinking neighborhoods, tangential convergence, or asymptotic conditions on the associated weights control how the algorithmic iterates approach these paths.

Our analysis removes the strict complementarity assumption, but replaces the classical endpoint regularity by an analytic reparametrization of weighted central paths. Since our method uses higher-order derivatives of weighted paths whose weights vary with the iterates, the local super-linear convergence result additionally requires a uniform predictor regularity condition. This condition is different from strict complementarity and from the classical tangential conditions, although they play analogous roles in controlling the local behavior of the algorithm.
\end{remark}

\begin{remark}[Relation to higher-order methods for \LCP]
Our predictor-corrector construction is inspired by the higher-order \IPM~framework of Stoer et al.~\cite{SWM98} for \LCP s. For sufficient \LCP s, Stoer and Wechs~\cite{SW98}
subsequently established polynomial iteration bounds for these higher-order predictor-corrector methods. Potra and Stoer~\cite{PS10} later developed a related class of super-linearly convergent polynomial-time higher-order \IPM s for sufficient \LCP s, obtaining arbitrarily high local convergence orders while retaining polynomial iteration complexity.

The semi-definite setting is substantially more delicate: weighted central paths depend on the matrix symmetrization, may fail to be analytic at the limiting point, and can have large ramification indices (see~\cite{BM24} or Section~\ref{ssec:EdK} of this paper). 

At the level of the corrector step, our construction is closely related
to the \AHO~centering mechanism of Kojima et al.~\cite{KSS99}. In particular, the damping rule in~\eqref{step_length} is based on the normalized product $\|\Delta X\|_F\|\Delta S\|_F$. However, our use of $\|\Delta X\|_F\|\Delta S\|_F$ is different from~\cite{KSS99}: we do not require an asymptotic bound on $\|\Delta X\|_F\|\Delta S\|_F$. Instead, for each fixed centering parameter, continuity and compactness yield the finite bound needed to prove restoration of the inner neighborhood after finitely many damped corrector steps; see
Lemma~\ref{lem:finite_corrector}.
\end{remark}

\subsection{Algebraic methods and \IPM s}
Algebraic and semi-algebraic methods have been used both to study the structure of central paths and to compute exact solutions of \SDO. This perspective is particularly
useful because, even for rational input data, an exact optimal solution need not have rational coordinates (see e.g.,~\cite[Page~17]{Kl02}). Thus, exact computation in this setting generally means representing the relevant coordinates as real
algebraic numbers, for example through univariate representations and Thom encodings (see Section~\ref{sign_cond_Thom_encod}), rather than computing finite-precision approximations.

\vspace{5px}
\noindent
The exact computational complexity of \SDO~is considerably more subtle than its convexity might suggest. In the real-number model of
computation~\cite{BCSS98}, the \SDO~feasibility problem belongs to $\mathrm{NP}\cap\mathrm{coNP}$; see, e.g.,~\cite{B17,R97}. In the classical bit
model~\cite{GJ79}, even membership of general \SDO~feasibility in $\mathrm{NP}$ is
not known~\cite{R97}. Algorithmic real algebraic geometry nevertheless provides finite exact procedures for describing solutions of polynomial systems
arising from \SDO~(see~\cite[Algorithm~14.9]{BPR06} or~\cite[Section~3.1]{BM22}). Such procedures have singly exponential complexity in the relevant dimensions and degrees.

\vspace{5px}
\noindent
Algebraic geometry has also been used to quantify geometric and algebraic properties of central paths~\cite{ABGJ18, DSV12, DMS05,HKR02, HLMZ21, HST21,KMNY91}. Complementary symbolic and numerical algebraic geometry approaches have been developed for
computing solutions of \SDO~independently of \IPM s~\cite{BPR06,HLMZ21, HNS16, CNS20}. For example, Hauenstein et al.~\cite{HLMZ21} interpret the central path as a bilinear homotopy and use analytic reparametrization
to improve numerical endgames. Algebraic techniques have further provided quantitative bounds on the degree, ramification, and
convergence rate of the central path of~\SDO 

\begin{proposition}[Theorem~1.1 in~\cite{BM22}]
Let $(X^{**},y^{**},S^{**})$ denote the limit point of the central
path. Then, for all sufficiently small $\mu>0$,
\begin{align*}
 \|X(\mu)-X^{**}\|_F
   = O\!\left(\mu^{1/\gamma}\right),
 \qquad
 \|S(\mu)-S^{**}\|_F
   = O\!\left(\mu^{1/\gamma}\right),
\end{align*}
where
\begin{align*}
 \gamma = 2^{O(m+n^2)}.
\end{align*}
\end{proposition}
\noindent
These results quantify aspects of the algebraic and geometric complexity of \SDO, but such bounds do not by themselves translate into
convergence order guarantees for \IPM s. In
particular, bounds on algebraic degree~\cite{BM22}, error bound exponents~\cite{BM24b}, ramification indices~\cite{BM24}, or the convergence rate of the central path~\cite{BM22} do not directly specify how an \IPM~should exploit this structure.

In contrast with these largely descriptive uses of algebraic geometry, the present work uses the ramification structure of weighted \AHO~paths directly in the construction and convergence analysis of a higher-order predictor-corrector \IPM.

\section{Analytic reparametrization of weighted central paths}\label{weighted_central_paths}
The existence of a weighted central path for any $W \succ 0$ is well-known. For the sake of completeness, we include the proof here.
\begin{notation}
For ease of exposition, we define and use the following linear map: 
\begin{align*}
    \mathcal{A}(X):=\begin{bmatrix} \langle A_1, X \rangle\\ \vdots\\ \langle A_m, X \rangle \end{bmatrix}.
\end{align*}
\end{notation}
\begin{proposition}[Existence and analyticity of weighted \AHO~central paths]
\label{prop:weighted_path_existence}
For every $\mu>0$ and every
$W\in\mathbb S_{++}^n$, there exists a unique strictly primal-dual
feasible solution
\begin{align*}
    (X(\mu,W),y(\mu,W),S(\mu,W))
\end{align*}
satisfying~\eqref{Weighted-Central_Path}. Further, the mapping $\xi(\mu,W)$ is well-defined and jointly analytic on $(0,\infty)\times \mathbb{S}^n_{++}$.
\end{proposition}

\begin{proof}
The proof follows from~\cite[Theorem~1(a)]{MP98}. We define the affine map
\begin{align*}
F:\mathbb{S}^n\times\mathbb{S}^n\times\mathbb{R}^m
   \longrightarrow
   \mathbb{S}^n\times\mathbb{R}^m
\end{align*}
by
\begin{align*}
F(X,S,y)
:=
\bigl(
\mathcal A^*(y)+S-C,\,
\mathcal A(X)-b
\bigr).
\end{align*}
Since $A_1,\ldots,A_m$ are linearly independent, $\mathcal A^*$ is injective, and hence $F$ is \textit{$y$-injective} according to~\cite{MP98}. Furthermore, if $F(X,S,y)=F(\tilde X,\tilde S,\tilde y)$, then
\begin{align*}
\mathcal A(X-\tilde X)=0,
\qquad
S-\tilde S=-\mathcal A^*(y-\tilde y),
\end{align*}
and therefore
\begin{align*}
\langle X-\tilde X,S-\tilde S\rangle
=
-\langle\mathcal A(X-\tilde X),y-\tilde y\rangle
=0.
\end{align*}
Thus, $F$ is \textit{$(X,S)$-equilevel-monotone}. Since $F$ is affine, $y$-injectivity also implies the \textit{$y$-boundedness} condition required
in~\cite[Theorem~1(a)]{MP98}; see also~\cite[Lemma~3]{MP98}.

\vspace{5px}
\noindent
Now, consider the \AHO~interior-point map
\begin{align*}
\Phi_{\text{\AHO}}(X,S)=\frac{XS+SX}{2}
\end{align*}
and define
\begin{align*}
H:=(\Phi_{\text{\AHO}},F), \quad \mathcal U
:=
\left\{
(X,S)\in\mathbb S_{++}^n\times\mathbb S_{++}^n
\;\middle|\;
XS+SX\in\mathbb S_{++}^n
\right\}.
\end{align*}
By Assumption~\ref{IPM_Cond}, the central path exists~\cite[Theorem~3.1]{Kl02}. Hence, for any central point, we have
\begin{align*}
F(X,S,y)=0,
\qquad
XS+SX=2\mu I_n\succ0,
\end{align*}
so that $(X,S)\in\mathcal U$. All this implies $0\in F(\mathcal U\times\mathbb R^m)$. It now follows from~\cite[Theorem~1(a)]{MP98} that
\begin{align*}
H:\mathcal U\times\mathbb R^m
\longrightarrow
\mathbb S_{++}^n
\times
F(\mathcal U\times\mathbb R^m)
\end{align*}
is a homeomorphism. Therefore, for every $\mu>0$ and $W\in\mathbb S_{++}^n$, there exists a unique
$(X(\mu,W),S(\mu,W),y(\mu,W))$ satisfying
\begin{align*}
H(X(\mu,W),S(\mu,W),y(\mu,W))
=
(\mu W,0),
\end{align*}
which proves existence and uniqueness, and continuity of weighted central paths.

\vspace{5px}
\noindent
It remains to establish joint analyticity. The map $H$ is polynomial,
hence real analytic. Furthermore, the proof of~\cite[Lemma~4]{MP98}
shows that $DH(X,S,y)$ is an isomorphism at every
$(X,S,y)\in\mathcal U\times\mathbb R^m$.
Hence, by the analytic inverse function theorem (see e.g.,~\cite[Proposition~3.24]{BPR06}), $H^{-1}$ is real
analytic locally throughout its image. Since $H$ is globally one-to-one,
these local analytic inverses agree on overlaps and therefore define a
real-analytic inverse on
\begin{align*}
\mathbb S_{++}^n
\times
F(\mathcal U\times\mathbb R^m).
\end{align*}
Finally, the map $(\mu,W)\longmapsto(\mu W,0)$ is real analytic on $(0,\infty)\times\mathbb S_{++}^n$. Thus, $\xi(\mu,W)=H^{-1}(\mu W,0)$
is jointly real analytic on $(0,\infty)\times\mathbb S_{++}^n$.
\end{proof}
\noindent
For the purpose of local convergence of \IPM s in this paper, we restrict weighted central paths to the bounded set
\begin{equation}\label{restricted_weighted_central_paths}
\begin{aligned}
\bigg\{(X,y,S) \mid \langle &A_i, X \rangle =b_i, \ \ i=1,\ldots,m, \\
&\sum_{i=1}^m y_iA_i + S -C=0,\ \ 
X S + SX - 2\mu W=0, \ \ X,S \succ 0\\
&\text{for some} \ \|W - I_n \|_F \le \beta, \ \ \mu \in (0,1]\bigg\},     
\end{aligned}
\end{equation}
where $\beta < 1$ is a fixed parameter (see~\eqref{neighborhood_definition}).

\vspace{5px}
\noindent
The semi-algebraicity of the weighted central paths together with their local boundedness implies that a weighted central path converges for every $W \succ 0$. 

\begin{proposition}\label{weighted_central_path_convergence}
The mapping
\begin{align*}
\xi:(0,\infty)\times\mathbb{S}_{++}^n
\longrightarrow
\mathbb{S}^n\times\mathbb{R}^m\times\mathbb{S}^n
\end{align*}
is jointly semi-algebraic. Moreover, for every fixed $W\succ0$, the weighted central path $\xi_W(\mu)$ converges as $\mu\downarrow0$.
\end{proposition}

\begin{proof}
By Proposition~\ref{prop:weighted_path_existence}, for every
$\mu>0$ and $W\succ0$ there exists a unique weighted central solution.
Hence, the graph of $\xi$ is precisely the semi-algebraic set
\[
\left\{
(\mu,W,X,y,S)\ \middle|\
\begin{array}{l}
\mu>0,\quad W\succ0,\quad X,S\succ0,\\
\mathcal A(X)=b,\quad
\mathcal A^*(y)+S=C,\\
XS+SX=2\mu W
\end{array}
\right\}.
\]
Therefore, $\xi$ is a semi-algebraic mapping.

\vspace{5px}
\noindent
Now, we fix $W\succ0$. We first show that $\xi_W(\mu)$ remains bounded for $0<\mu\le1$. Let $(\bar X,\bar y,\bar S)$ be any fixed strictly
primal-dual feasible point, whose existence follows from
Assumption~\ref{IPM_Cond}. For every weighted central point
$(X(\mu,W),y(\mu,W),S(\mu,W))$, primal-dual feasibility gives
\begin{align*}
\la \bar S, X(\mu,W)\ra+\la \bar X, S(\mu,W)\ra
=
\la X(\mu,W), S(\mu,W)+\la \bar X,\bar S\ra.
\end{align*}
Additionally, we have
\begin{align*}
\la X(\mu,W),S(\mu,W)\ra
=
\mu\,\trace(W).
\end{align*}
Consequently, we get
\begin{align*}
\la \bar S, X(\mu,W)\ra+ \la \bar X, S(\mu,W) \ra
=
\mu\,\trace(W)+ \la \bar X,\bar S \ra,
\end{align*}
which is uniformly bounded for $0<\mu\le1$. Since
$\bar X,\bar S\succ0$ and $X(\mu,W),S(\mu,W)\succ0$, it follows that
$X(\mu,W)$ and $S(\mu,W)$ are uniformly bounded. Further, the linear independence of $A_1,\ldots,A_m$ implies that $y(\mu,W)$ is uniformly bounded as
well. Thus, $\xi_W:(0,1]\to
\mathbb{S}^n\times\mathbb{R}^m\times\mathbb{S}^n$ is a bounded semi-algebraic function, and By~\cite[Proposition~3.18]{BPR06}, each of its coordinate functions has a finite limit as $\mu\downarrow0$. 
\end{proof}
\noindent
From Proposition~\ref{weighted_central_path_convergence}, each coordinate of $\xi_W$ is described by Puiseux series $\sum_{i=r}^{\infty} c_i \mu^{i/q}$, when $\mu$ is sufficiently small, where $q \in \mathbb{Z}_{>0}$ is a constant integer, and $c_i$ depends on $W$. Now, we prove Theorem~\ref{global_ramification}.
\begin{proof}[Proof of Theorem~\ref{global_ramification}]
We show that there are only finitely many ramification indices for weighted central paths restricted to~\eqref{restricted_weighted_central_paths}. The proof of Proposition~\ref{weighted_central_path_convergence} already establishes the uniform boundedness of weighted central paths restricted to~\eqref{restricted_weighted_central_paths}.

\vspace{5px}
\noindent
We define $Q \in \mathbb{R} [\mu,W, V_1,\ldots, V_{\bar{n}}]$ as follows
\begin{align*}
Q
:=
\|\mathcal A(X)-b\|^2
+\|\mathcal A^*(y)+S-C\|_F^2
+\|XS+SX-2\mu W\|_F^2.
\end{align*}
We choose $\varepsilon>0$ sufficiently small so that
\begin{align*}
\|(X(\mu,W),y(\mu,W),S(\mu,W))\|_F
<
\frac{1}{\varepsilon}
\end{align*}
for every $\mu\in(0,1]$ and $W\in\mathcal W_\beta$, and we introduce an
additional variable $V_{\bar n+1}$. We also define
\begin{align*}
\tilde Q
:=
Q^2+
\left(
\varepsilon^2
\sum_{i=1}^{\bar n+1}V_i^2-1
\right)^2.
\end{align*}
The projection of
\begin{align*}
\zero
\bigl(\tilde Q(\mu,W,\cdot),\mathbb R^{\bar n+1}\bigr)
\end{align*}
onto the first $\bar n$ coordinates contains the weighted central point. Applying Parameterized Bounded Algebraic Sampling
\cite[Algorithm~12.18]{BPR06} to $\tilde Q$, with
$(\mu,W)$ as parameters, produces a finite family
\begin{align*}
\mathcal U=\{u_1,\ldots,u_M\}
\end{align*}
of parameterized univariate representations
\begin{align*}
u_a=(f_a,g_a)
=
\bigl(
f_a,(g_{a,0},g_{a,1},\ldots,g_{a,\bar n+1})
\bigr),
\qquad a=1,\ldots,M.
\end{align*}
For every fixed $(\mu,W)$, the corresponding sample points meet every
semi-algebraically connected component of
\begin{align*}
\zero
\bigl(\tilde Q(\mu,W,\cdot),\mathbb R^{\bar n+1}\bigr).
\end{align*}
\noindent
Consequently, after projection onto the first $\bar n$ coordinates, the weighted central point is represented by at least one element
$u_a\in\mathcal U$. More precisely, its coordinates have the form
\begin{align*}
\hat v_i(\mu,W)
=
\frac{g_{a,i}(\mu,W,t)}
     {g_{a,0}(\mu,W,t)},
\qquad i=1,\ldots,\bar n,
\end{align*}
where $t$ is a real root of $f_a(\mu,W,T)$ and $g_{a,0}(\mu,W,t)\neq0$. Set
\begin{align*}
d:=\max_{1\le a\le M}\deg_T(f_a),
\end{align*}
and consider the finite family
\begin{align*}
\mathcal P
:=
\bigcup_{a=1}^M\operatorname{Der}_T(f_a),
\end{align*}
where
\begin{align*}
\Der_T(f_a):=\left\{f_a,\frac{\partial f_a}{\partial T},\ldots,\frac{\partial^d f_a}{\partial^d T}\right\}.
\end{align*}
\noindent
Choose a cylindrical algebraic decomposition adapted to
$\mathcal P$ (see Section~\ref{cad}). Over every cell of the induced decomposition of the
$(\mu,W)$-parameter space, the real roots of each $f_a$ are delineable and their Thom encodings remain fixed. Refining the
decomposition, if necessary, by the finitely many semi-algebraic
conditions determining which of the above sample points is the
weighted central point, we may assume that the relevant univariate
representation and Thom encoding are fixed on each cell. 
Since we are interested only in sufficiently small positive values of
$\mu$, we refine the induced decomposition at $\mu=0$. This yields a
finite semi-algebraic partition
\begin{align*}
\mathcal W_\beta
=
\bigcup_{\ell=1}^s Z_\ell
\end{align*}
with the following property. For each $\ell$, there exist an index
$a_\ell\in\{1,\ldots,M\}$ and a Thom encoding $\sigma_\ell$ such that,
for every $W\in Z_\ell$, there exists $\varepsilon_W>0$ for which
\begin{align}\label{weighted_CP_rep}
v_i(\mu,W)
=
\frac{
g_{a_\ell,i}
\bigl(\mu,W,t_{\sigma_\ell}(\mu,W)\bigr)}
{
g_{a_\ell,0}
\bigl(\mu,W,t_{\sigma_\ell}(\mu,W)\bigr)
},
\qquad i=1,\ldots,\bar n,
\end{align}
for every $0<\mu<\varepsilon_W$, where
$t_{\sigma_\ell}(\mu,W)$ is the real root of
$f_{a_\ell}(\mu,W,T)$ determined by $\sigma_\ell$. Now, fix $W\in Z_\ell$. After specializing $W$, let
\begin{align*}
f_W(\mu,T)
:=
f_{a_\ell}(\mu,W,T),
\qquad
d(W):=\deg_T(f_W).
\end{align*}
The specialization associated with the above representation is non-zero, and $d(W)\le d$. Further, we have
\begin{align*}
f_W
\bigl(\mu,t_{\sigma_\ell}(\mu,W)\bigr)=0,
\qquad
f_W(\mu,T)\in\mathbb R[\mu,T]
\subset\mathbb C[\mu,T].
\end{align*}
By Proposition~\ref{Puiseux_Newton_Thm},
$t_{\sigma_\ell}(\mu,W)$ admits a convergent Puiseux expansion with some ramification index $q(W)$ satisfying
\begin{align*}
q(W)\le d(W)\le d.
\end{align*}
\noindent
Notice that each coordinate in~\eqref{weighted_CP_rep} is a rational function of $t_{\sigma_\ell}(\mu,W)$ with coefficients in $\mathbb R(\mu)$.
Hence, all coordinates belong to the same finite Puiseux extension,
and the ramification index of every coordinate of the weighted central
path divides $q(W)$. Consequently, every coordinate ramification
index is a positive integer not exceeding $d$. Let $\bar q$ be the least common multiple of  $1,2,\ldots,d$. Then it follows that $\bar q$ is a common multiple of the ramification indices of all coordinates
of all weighted central paths with $W\in\mathcal W_\beta$. By Proposition~\ref{weighted_central_path_convergence}, for every fixed $W$ the weighted central path converges as $\mu\downarrow0$. Hence, its coordinate Puiseux expansions contain no negative powers. Therefore,
\begin{align*}
X(\mu^{\bar q},W),\quad
y(\mu^{\bar q},W),\quad
S(\mu^{\bar q},W)
\end{align*}
are analytic at $\mu=0$ for every $W\in\mathcal W_\beta$. All this implies that the set of positive integers $q$ for which $\xi_W(\mu^q)$ is analytic at $\mu=0$ for every $W\in\mathcal W_\beta$ is non-empty. We denote the smallest integer by $N$, which completes the proof.
\end{proof}

\begin{example}
Consider again the \SDO~problem in Example~\ref{derivatives_diverge}. We now introduce the one-parameter family of weights
\begin{align*}
W_\theta
:=
\kappa_\theta
\begin{pmatrix}
1 & 0 & 0\\
0 & 1+\theta^2 & \theta\\
0 & \theta & 1+\theta^2
\end{pmatrix},
\qquad
\kappa_\theta:=\frac{3}{3+2\theta^2},
\qquad
\theta\in\mathbb R .
\end{align*}
The matrix $W_\theta$ is positive definite for every $\theta\in\mathbb R$. The corresponding weighted central path is determined by
\begin{align*}
XS+SX=2\mu W_\theta,
\end{align*}
where
\begin{align*}
X=
\begin{pmatrix}
1 & X_{12} & X_{13}\\
X_{12} & 1 & X_{23}\\
X_{13} & X_{23} & 1
\end{pmatrix},
\qquad
S=
\begin{pmatrix}
-y_1 & 2 & -2\\
2 & -y_2 & -1\\
-2 & -1 & -y_3
\end{pmatrix}.
\end{align*}
The problem and $W_\theta$ are invariant under the orthogonal
transformation
\begin{align*}
P=
\begin{pmatrix}
-1&0&0\\
0&0&1\\
0&1&0
\end{pmatrix}.
\end{align*}
By uniqueness of the weighted central solution, $X_{13}=-X_{12}$ and $y_2=y_3$. Set $x:=X_{12}$, $z:=X_{23}$, and $\nu:=\kappa_\theta\mu$. Then
\begin{align*}
y_1=4x-\nu,
\qquad
y_2=y_3=2x-z-(1+\theta^2)\nu,
\end{align*}
and the remaining equations reduce to
\begin{align*}
\begin{cases}
-6x^2+xz+(\theta^2+2)\nu x+x-2z+4=0,\\[1mm]
z^2-2xz+(1+\theta^2)\nu z-\theta\nu-2x-1=0.
\end{cases}
\end{align*}
For sufficiently small $\mu$, the first equation gives
\[
z=
\frac{
6x^2-(\theta^2+2)\nu x-x-4
}{
x-2
}.
\]
Substitution into the second equation yields $F_\theta(\nu,x)=0$, where
\begin{align*}
F_\theta(\nu,x)
={}&24x^4
+\bigl(12-(4\theta^2+14)\nu\bigr)x^3\\
&+\Bigl(
(\theta^2+2)\nu^2
-(15\theta^2+\theta+17)\nu
-36
\Bigr)x^2\\
&+\Bigl(
(2\theta^4+6\theta^2+4)\nu^2
+(6\theta^2+4\theta+14)\nu
-12
\Bigr)x\\
&+12+(8\theta^2-4\theta+8)\nu .
\end{align*}
At $\nu=0$,
\begin{align*}
F_\theta(0,x)=12(x-1)(x+1)^2(2x-1),
\end{align*}
thus the weighted central path is the branch converging to $x=-1$.

The corresponding Puiseux expansion is
\begin{align}
X_{12}(\mu)
={}&-1
+\sqrt{\frac{\theta^2+\theta+1}{8}}\,
      \nu^{1/2}
+\frac{1+\theta-3\theta^2}{32}\,\nu
\nonumber\\
&+
\frac{
85\theta^4-254\theta^3+31\theta^2-294\theta-99
}{
4608\sqrt{2(\theta^2+\theta+1)}
}
\,\nu^{3/2}
+O(\nu^2),
\label{eq:weighted-elliptope-puiseux}
\end{align}
where $\nu=\kappa_\theta\mu$. The remaining coordinates are recovered from
\begin{align*}
&X_{13}(\mu)=-X_{12}(\mu),\\
&X_{23}(\mu)
=
\frac{
6X_{12}(\mu)^2
-(\theta^2+2)\nu X_{12}(\mu)
-X_{12}(\mu)-4
}{
X_{12}(\mu)-2
},
\end{align*}
and
\begin{align*}
y_1(\mu)=4X_{12}(\mu)-\nu,
\qquad
y_2(\mu)=y_3(\mu)
=
2X_{12}(\mu)-X_{23}(\mu)-(1+\theta^2)\nu.
\end{align*}
\noindent
Since
\begin{align*}
\theta^2+\theta+1
=
\left(\theta+\frac12\right)^2+\frac34>0
\qquad\forall\,\theta\in\mathbb R,
\end{align*}
the coefficient of $\nu^{1/2}$ in~\eqref{eq:weighted-elliptope-puiseux} never vanishes. Therefore,
\begin{align*}
q(\theta)=2,
\qquad \forall\,\theta\in\mathbb R.
\end{align*}
Consequently, for every interval $I\subseteq\mathbb R$, bounded or unbounded, the minimal common ramification index of the family $\{W_\theta:\theta\in I\}$ is
\begin{align*}
N_I=2.
\end{align*}
For $\theta=0$, $W_0=I_3$, and the expansion reduces to the ordinary central path expansion in Example~\ref{derivatives_diverge}.
\end{example}

\begin{remark}
For the \SDO~instance in Example~\ref{derivatives_diverge}, varying the positive definite weight cannot produce two different ramification indices on two parameter intervals. In fact, if $W=(w_{ij})\succ0$ is arbitrary, the branch converging to the non-strictly complementary optimal solution satisfies
\[
X_{12}(\mu)
=
-1+\frac14
\sqrt{w_{22}+2w_{23}+w_{33}}\,\mu^{1/2}+O(\mu).
\]
Since
\[
w_{22}+2w_{23}+w_{33}
=
(0,1,1)W(0,1,1)^T>0,
\]
the square-root term is present for every $W\succ0$. Thus the ramification index is always equal to $2$ for this particular \SDO~instance.  
\end{remark}

\section{A higher-order \IPM}\label{higher-order_IPM}
\subsection{\AHO~search directions}
Several primal-dual search directions have been developed for \SDO, including Nesterov-Todd~\cite{NT97, NT98}, Alizadeh-Haeberly-Overton (\AHO)~\cite{AHO98}, Kojima-Shindoh-Hara~\cite{KSH97}, or Monteiro-Zhang~\cite{MZ98}). We focus on the \AHO~search direction, which is obtained by applying Newton's method to~\eqref{central_path_description}. Given an initial solution $(X_0,y_0,S_0)$, this leads to the following linear system of equations
\begin{equation}\label{lin_Central_Path}
 \begin{aligned}
  \langle A_i, \Delta X \rangle &=0, \quad i=1,\dots, m,\\
  \sum_{i=1}^m \Delta y_i A_i + \Delta S&=0, \\
  \Delta X S_0 + S_0 \Delta X + X_0 \Delta S + \Delta S X_0   &=2\mu I_n - X_0S_0 - S_0X_0.
 \end{aligned}
\end{equation}
When $(X_0,y_0,S_0)$ with $X_0,S_0\succ 0$ is sufficiently close to the central path~\cite[Theorem~2.1]{MZ97}, the system~\eqref{lin_Central_Path} has a unique solution $(\Delta X, \Delta y, \Delta S)$. 

\begin{remark}[Semi-algebraic formulations]\label{other_formulations}
The centrality condition for the central path appears in several equivalent forms in the literature, for example,
\begin{align*}
XS=\mu I_n,\qquad
X^{\frac12}SX^{\frac12}=\mu I_n,\qquad
S^{\frac12}XS^{\frac12}=\mu I_n,
\end{align*}
or, more generally, through symmetrized formulations of the form
\begin{align*}
PXSP^{-1}+P^{-T}SXP^{T}=2\mu I_n,
\end{align*}
where $P$ is non-singular; see, e.g.,~\cite{MZ98}.
If the positive-definiteness constraints $X,S\succ0$ are omitted, then, for a fixed
$\mu>0$, these equations generally define different zero sets. In contrast,
when restricted to $X,S\succ0$, the corresponding centrality systems are
equivalent and determine the same unique central solution for each $\mu>0$.
Consequently, they give rise to the same semi-algebraic central path function.
\end{remark}

\vspace{5px}
\noindent
Analogously, \AHO-type search directions can be obtained by the application of the Newton's method to the parameterized family of weighted central paths~\eqref{Weighted-Central_Path}. The following proposition shows that the existence of such \AHO-type directions still follows from~\cite[Theorem~2.1]{MZ97}, when $(X_0,y_0,S_0)$ with $X_0,S_0\succ 0$ is sufficiently close to the central path.
\begin{notation}
We define the linear operator
\begin{align*}
    \mathcal J_{X,S}:
    \mathbb S^n\times\mathbb R^m\times\mathbb S^n
    \longrightarrow
    \mathbb R^m\times\mathbb S^n\times\mathbb S^n
\end{align*}
by
\begin{align*}
    \mathcal J_{X,S}(\Delta X,\Delta y,\Delta S)
    :=
    \left(
        \mathcal A(\Delta X),\,
        \mathcal A^*(\Delta y)+\Delta S,\,
        \Delta XS + S\Delta X+ X\Delta S+\Delta SX
    \right),
\end{align*}
where
\begin{align*}
    \mathcal A^*(\Delta y)
    :=
    \sum_{i=1}^m \Delta y_i A_i
\end{align*}
is the adjoint of the linear map $\mathcal{A}$.
\end{notation}
\begin{proposition}[Well-definedness of weighted \AHO~systems]
\label{lem:weighted_AHO_well_defined}
Let $\mu>0$, and let $(X,y,S)$ be a strictly primal-dual feasible point with
$X,S\succ0$. Suppose that
\begin{align}\label{neighborhood_condition}
    \left\|
        \frac{XS+SX}{2\mu}-I_n
    \right\|_F
    <\frac12.
\end{align}
\noindent
Then $\mathcal J_{X,S}$ is non-singular. Consequently, for every
$H\in\mathbb S^n$, the linear system
\begin{align*}
\begin{aligned}
    \langle A_i,\Delta X\rangle &=0,
        && i=1,\ldots,m,\\
    \sum_{i=1}^m \Delta y_iA_i+\Delta S &=0,\\
    \Delta XS + S\Delta X+ X\Delta S+\Delta SX &=H
\end{aligned}
\end{align*}
has a unique solution
$(\Delta X,\Delta y,\Delta S)$. In particular, for every $W\succ0$, the weighted \AHO~Newton system
\begin{equation*}
\begin{aligned}
    \langle A_i,\Delta X\rangle &=0,
        \quad i=1,\ldots,m,\\
    \sum_{i=1}^m \Delta y_iA_i+\Delta S &=0,\\
    \Delta XS + S\Delta X+ X\Delta S+\Delta SX
        &=2\mu W-XS-SX
\end{aligned}
\end{equation*}
has a unique solution whenever the neighborhood condition~\eqref{neighborhood_condition} holds.
\end{proposition}
\begin{proof}
Using the identity
\begin{align*}
    \left\|
        \frac{XS+SX}{2\mu}-I_n
    \right\|_F^2
    =
    \left\|
        \frac{S^{1/2}XS^{1/2}}{\mu}-I_n
    \right\|_F^2
    +
    \frac{1}{4\mu^2}\|XS-SX\|_F^2,
\end{align*}
we obtain
\begin{align*}
    \left\|
        \frac{S^{1/2}XS^{1/2}}{\mu}-I_n
    \right\|_F
    \le
    \left\|
        \frac{XS+SX}{2\mu}-I_n
    \right\|_F
    <\frac12.
\end{align*}
Hence, the sufficient condition of
\cite[Theorem~2.1]{MZ97}
for existence and uniqueness of the \AHO~direction is satisfied. All this means that the coefficient matrix representing $\mathcal J_{X,S}$ is non-singular. Thus, for every $H\in\mathbb S^n$, the system
\begin{align*}
    \mathcal J_{X,S}(\Delta X,\Delta y,\Delta S)
    =
    (0,0,H)
\end{align*}
has a unique solution. Taking $H=2\mu W-XS-SX$ gives the weighted \AHO~Newton direction.
\end{proof}
\subsection{A higher-order \IPM~based on \AHO~search directions}
Our construction is inspired by the higher-order predictor-corrector \IPM~for \LCP~in~\cite{SWM98} and by the \AHO~centering directions in~\cite{KSS99}.

\vspace{5px}
\noindent
Fix a positive integer reparametrization exponent $\rho$ and a Taylor
order $p$. We consider the family of reparametrized weighted \AHO~central paths
\begin{equation}\label{reparam-Weighted-Central_Path}
\begin{aligned}
 \langle A_i, X \rangle &= b_i,
 \quad i=1,\ldots,m,\\
 \sum_{i=1}^m y_i A_i + S &= C,\\
 XS+SX &= 2t^\rho W, \quad X,S \succ 0,
\end{aligned}
\end{equation}
where $W \in \mathcal{W}_\beta$. For a fixed
$W \succ 0$, we denote the solution of~\eqref{reparam-Weighted-Central_Path}
by
\begin{align*}
(X(t),y(t),S(t)):=(X(t^\rho,W),y(t^\rho,W),S(t^\rho,W)).
\end{align*}
\noindent
Throughout the construction below, we assume that the chosen reparametrization exponent $\rho$ provides the derivatives required
through order $p$ at the current iterate. Admissibility of $\rho$ is not required for the construction itself; it enters later in the convergence analysis (see Section~\ref{sec:local_convergence}).

\vspace{5px}
\noindent
\paragraph{\textbf{Predictor step}}
The predictor step of our higher-order \IPM~involves computing the $p$th-order Taylor approximation of~\eqref{reparam-Weighted-Central_Path} at $t_0$ as follows
\begin{equation}\label{Taylor_Approx_Central_Path}
\begin{aligned}
 \big(\hat{X}_{t_0}(t), \hat{y}_{t_0}(t), \hat{S}_{t_0}(t)\big) &= \big(X(t_0), y(t_0),S(t_0)\big) \\
 &\qquad+  \sum_{j=1}^p \frac{1}{j !} (t - t_0)^j \big(X^{(j)}(t_0), y^{(j)}(t_0), S^{(j)}(t_0)\big),
\end{aligned}
\end{equation}
where the higher-order derivatives are obtained by solving 
\begin{equation}\label{HOSys}
    \begin{aligned}
        \langle A_i, X^{(j)} \rangle &=0, \ \ \  i=1,\dots,m,\\[-1\jot]
        \sum_{i=1}^m y_i^{(j)} A_i + S^{(j)} &=0,\\[-1\jot]
        X^{(j)} S(t)+X(t) S^{(j)} + S^{(j)} X(t) + S(t) X^{(j)}&=\\
        &\begin{cases} 2\rho (\rho - 1)\cdots (\rho-j+1) t^{\rho - j} W + \xi(j) \ \ &j\le \rho,\\
        \xi(j), \ \ &j>\rho,
        \end{cases}
    \end{aligned}
\end{equation}
in which
\begin{align*}
\xi(j):=-\sum_{k=1}^{j-1} {j\choose k}\big( X^{(k)}S^{(j-k)} + S^{(k)}X^{(j-k)}\big).
\end{align*}
We maintain feasibility and centrality by keeping the iterates in the following neighborhood of the central path: 
\begin{align}\label{neighborhood_definition}
    \mathcal{N}_t(\beta):=\bigg\{(X,y,S) \mid  \vartheta_t(X,S) \le \beta \bigg\}, \quad \vartheta_t(X,S):=\left\|\frac{XS + SX}{2t^{\rho}} - I_n \right\|_F, \qquad \beta < 1.
\end{align}
\begin{lemma}\label{positive_definite_argument}
If $\big(\hat{X}_{t_0}(t), \hat{y}_{t_0}(t), \hat{S}_{t_0}(t)\big) \in \mathcal{N}_t(\beta)$ for all $t \in [\bar{t},t_0]$, then $\big(\hat{X}_{t_0}(t), \hat{y}_{t_0}(t), \hat{S}_{t_0}(t)\big)$ is primal-dual feasible, and we have $\hat{X}_{t_0}(t), \hat{S}_{t_0}(t) \succ 0$ for all $t \in [\bar{t},t_0]$. 
\end{lemma}

\begin{proof}
Since $\vartheta_t\bigl(\hat X_{t_0}(t),\hat S_{t_0}(t)\bigr)
\le \beta<1$, we have
\begin{align*}
\hat X_{t_0}(t)\hat S_{t_0}(t)
+\hat S_{t_0}(t)\hat X_{t_0}(t)\succ0
\end{align*}
for every $t\in[\bar t,t_0]$. We first show that $\hat X_{t_0}(t)$ and
$\hat S_{t_0}(t)$ are non-singular throughout this interval. Indeed, if $\hat X_{t_0}(t)v=0$ for some $v\neq0$, then
\begin{align*}
v^T\!\left(
\hat X_{t_0}(t)\hat S_{t_0}(t)
+\hat S_{t_0}(t)\hat X_{t_0}(t)
\right)v=0,
\end{align*}
contradicting the positive definiteness above. Hence
$\hat X_{t_0}(t)$ is non-singular, and the same argument applies to $\hat S_{t_0}(t)$. Since $\hat X_{t_0}(t)$ and $\hat S_{t_0}(t)$ are continuous symmetric matrix-valued functions of $t$, and
\begin{align*}
\hat X_{t_0}(t_0)=X(t_0)\succ0,
\quad
\hat S_{t_0}(t_0)=S(t_0)\succ0,
\end{align*}
their inertia cannot change on $[\bar t,t_0]$ without one of them
becoming singular. Therefore,
\begin{align*}
\hat X_{t_0}(t)\succ0,
\quad
\hat S_{t_0}(t)\succ0,
\qquad
t\in[\bar t,t_0].
\end{align*}
\noindent
Finally, primal-dual feasibility follows directly from~\eqref{Taylor_Approx_Central_Path}-\eqref{HOSys}.
\end{proof}
\noindent
To ensure adequate reduction in the duality gap, the iterates from the predictor step stay in the wider neighborhood $\mathcal{N}(2\beta)$, where the next $t$ is obtained from
\begin{align}\label{mu_reduction}
    t_+:=\inf_t\Bigg \{\bar{t} > 0 \mid \vartheta_t(\hat{X}_{t_0}(t), \hat{S}_{t_0}(t)) \le 2 \beta, \ \ \forall \ \bar{t} \le t \le t_0\Bigg\}.
\end{align}
\begin{remark}
By Proposition~\ref{lem:weighted_AHO_well_defined}, when $\beta<\frac14$, the \AHO~direction is well-defined throughout the larger predictor neighborhood $\mathcal{N}(2\beta)$.
\end{remark}

\begin{remark}[Numerical computation of $t_+$]
\label{rem:line-search}
The convergence analysis of our \IPM~uses the exact boundary parameter $t_+$ defined by~\eqref{mu_reduction}. However, in our numerical implementation, we approximate $t_+$ by a geometric line search. Starting from $t_0$, we successively reduce the path parameter by a fixed factor $\alpha_{\rm ls}\in(0,1)$ and retain the last accepted
value before the next trial point leaves $\mathcal N(2\beta)$. We denote this numerical approximation by $t_+^{\rm ls}$. In our experiments, we take $\alpha_{\rm ls}=0.9999$ (see Section~\ref{sec:numerical_demonstrations}).
\end{remark}

\vspace{5px}
\paragraph{\textbf{Corrector step}}
The corrector step of our higher-order \IPM~regains centrality (which is classically an important factor for the super-linear convergence of \IPM s). This step, starting from 
\begin{align*}
(X_0^c,y_0^c,S_0^c):=(\hat{X}_{t_0}(t_+),\hat{y}_{t_0}(t_+),\hat{S}_{t_0}(t_+)),
\end{align*}
performs multiple damped Newton corrections as follows
\begin{align}\label{damped_corrector}
X^c_{r+1}(\alpha):=X^c_{r} + \alpha \Delta X_r, \quad y^c_{r+1}(\alpha):=y^c_r + \alpha \Delta y_r, \quad S^c_{r+1}(\alpha):=S^c_r + \alpha \Delta S_r,
\end{align}
where $(\Delta X_r,\Delta y_r,\Delta S_r)$ is the solution of
\begin{equation}\label{higher_centring_direction}
 \begin{aligned}
  \langle A_i, \Delta X_r \rangle &=0, \ \ \  i=1,\ldots,m,\\
  \sum_{i=1}^m (\Delta y_r)_i A_i +\Delta S_r&=0,\\
  \Delta X_r S_r^c + S_r^c \Delta X_r + X_r^c\Delta S_r + \Delta S_r X_r^c    &=2t_+^{\rho} I -X_r^c S_r^c -S_r^c X_r^c,
 \end{aligned}
 r = 0,1,2,\ldots
 \end{equation}
and $X^c_{r+1} := X^c_{r+1}(\alpha_r)$, $y^c_{r+1} := y^c_{r+1}(\alpha_r)$, $S^c_{r+1} := S^c_{r+1}(\alpha_r)$, where 
\begin{align}\label{step_length}
    \alpha_r = \begin{cases} 1 \ \ &\delta_r \le \gamma_r/2\\ \cfrac{\gamma_r}{2\delta_r} \ \ &\delta_r > \gamma_r/2  \end{cases}, \quad \gamma_r:=\vartheta_{t_+}(X^c_r,S^c_{r}), \quad  \delta_r:=\frac{\|\Delta X_r\|_F \|\Delta S_r\|_F}{t_+^{\rho}}.
\end{align}
The corrector steps are continued until $\vartheta_{t_+}(X^c_r,S^c_r) \le \beta$. We call the new corrected point $(X^+,y^+,S^+):=(X^c_r,y^c_r,S^c_r)$. In what follows (Lemmas~\ref{corrector_back_to_neighborhood}-\ref{lem:finite_corrector}), we show that the iterate $(X^+,y^+,S^+)$ is strictly feasible, and can be obtained using finitely many corrector steps.   

\vspace{5px}
\noindent
The following lemma shows that each corrector step strictly reduces the centrality while keeping the strict feasibility of the iterate.
\begin{lemma}
\label{corrector_back_to_neighborhood}
Recall $t_+>0$ from~\eqref{mu_reduction}. Suppose that, at the $r$th
corrector step, $(X_r^c,y_r^c,S_r^c)$ is strictly primal-dual feasible and
\begin{align*}
    \beta
    <
    \gamma_r=\vartheta_{t_+}(X^c_r,S^c_{r})
    \le 2\beta.
\end{align*}
Let $(\Delta X_r,\Delta y_r,\Delta S_r)$ be the solution of
\eqref{higher_centring_direction}, and let $\delta_r$ and $\alpha_r$ be
defined as in~\eqref{step_length}. Then
\begin{align*}
    (X_{r+1}^c,y_{r+1}^c,S_{r+1}^c)
    :=
    \bigl(
        X_{r+1}^c(\alpha_r),
        y_{r+1}^c(\alpha_r),
        S_{r+1}^c(\alpha_r)
    \bigr)
\end{align*}
is strictly primal-dual feasible and
\begin{align*}
    \gamma_{r+1}
    :=
    \vartheta_{t_+}(X_{r+1}^c,S_{r+1}^c)
    \le
    \begin{cases}
        \delta_r \le \dfrac{\gamma_r}{2},
        & \text{if } \delta_r\le \dfrac{\gamma_r}{2},\\[2mm]
        \displaystyle
        \gamma_r\left(1-\frac{\gamma_r}{4\delta_r}\right)
        <\gamma_r,
        & \text{if } \delta_r>\dfrac{\gamma_r}{2}.
    \end{cases}
\end{align*}
Further, we have
\begin{align*}
    X_{r+1}^c(\bar\alpha)\succ0,
    \qquad
    S_{r+1}^c(\bar\alpha)\succ0
    \qquad
    \text{for every } \bar\alpha\in[0,\alpha_r].
\end{align*}
In particular,
\begin{align*}
    \gamma_{r+1}<\gamma_r\le2\beta,
\end{align*}
and thus every accepted corrector iterate remains in $\mathcal N_{t_+}(2\beta)$.
\end{lemma}

\begin{proof}
Set $\tau:=t_+^\rho$ and define
\begin{align*}
E_r:=
    \frac{X_r^cS_r^c+S_r^cX_r^c}{2\tau}
    -I_n,
\end{align*}
so that
\begin{align*}
    \|E_r\|_F=\gamma_r.
\end{align*}
For $\bar\alpha\in[0,1]$, define
\begin{align*}
E_{r+1}(\bar\alpha):=
    \frac{
        X_{r+1}^c(\bar\alpha)S_{r+1}^c(\bar\alpha)
        +
        S_{r+1}^c(\bar\alpha)X_{r+1}^c(\bar\alpha)
    }{2\tau}
    -I_n.
\end{align*}
Using \eqref{damped_corrector} and the third equation of
\eqref{higher_centring_direction}, we obtain
\begin{align*}
    E_{r+1}(\bar\alpha)
    =
    (1-\bar\alpha)E_r
    +
    \frac{\bar\alpha^2}{2\tau}
    \left(
        \Delta X_r\Delta S_r
        +
        \Delta S_r\Delta X_r
    \right).
\end{align*}
Therefore,
\begin{align}
    \|E_{r+1}(\bar\alpha)\|_F
    &\le
    (1-\bar\alpha)\gamma_r
    +
    \frac{\bar\alpha^2}{2\tau}
    \left\|
        \Delta X_r\Delta S_r
        +
        \Delta S_r\Delta X_r
    \right\|_F
    \notag\\
    &\le
    (1-\bar\alpha)\gamma_r
    +
    \bar\alpha^2\delta_r.
\label{eq:corrector_residual_bound}
\end{align}
\noindent
If $\delta_r\le\gamma_r/2$, then $\alpha_r=1$, and hence
\begin{align*}
\gamma_{r+1} \le
    \delta_r
    \le
    \frac{\gamma_r}{2}.
\end{align*}
If $\delta_r>\gamma_r/2$, then
\begin{align*}
    \alpha_r
    =
    \frac{\gamma_r}{2\delta_r}
    <1,
\end{align*}
and consequently
\begin{align*}
    \gamma_{r+1}
    &\le
    \left(
        1-\frac{\gamma_r}{2\delta_r}
    \right)\gamma_r
    +
    \frac{\gamma_r^2}{4\delta_r^2}\delta_r\\
    &=
    \gamma_r
    \left(
        1-\frac{\gamma_r}{4\delta_r}
    \right)
    <\gamma_r.
\end{align*}
\noindent
Now, it remains to prove that positive definiteness is preserved along the
accepted corrector segment. We define
\begin{align*}
    q_r(\bar\alpha)
    :=
    (1-\bar\alpha)\gamma_r
    +
    \bar\alpha^2\delta_r.
\end{align*}
By the definition of $\alpha_r$,
\begin{align*}
    2\bar\alpha\delta_r
    \le
    \gamma_r
    \qquad
    \text{for every }
    \bar\alpha\in[0,\alpha_r].
\end{align*}
Thus, we have
\begin{align*}
    q_r'(\bar\alpha)
    =
    -\gamma_r+2\bar\alpha\delta_r
    \le0,
\end{align*}
and \eqref{eq:corrector_residual_bound} yields
\begin{align*}
    \|E_{r+1}(\bar\alpha)\|_F
    \le
    q_r(\bar\alpha)
    \le
    q_r(0)
    =
    \gamma_r
    \le
    2\beta
    <1.
\end{align*}
Therefore, we get
\begin{align*}
    I_n+E_{r+1}(\bar\alpha)\succ0,
\end{align*}
and hence
\begin{align*}
    X_{r+1}^c(\bar\alpha)S_{r+1}^c(\bar\alpha)
    +
    S_{r+1}^c(\bar\alpha)X_{r+1}^c(\bar\alpha)
    \succ0.
\end{align*}
\noindent
Since $X_r^c,S_r^c\succ0$, analogous to the continuity argument of Lemma~\ref{positive_definite_argument},
neither $X_{r+1}^c(\bar\alpha)$ nor
$S_{r+1}^c(\bar\alpha)$ can become singular for
$\bar\alpha\in[0,\alpha_r]$. Hence
\begin{align*}
    X_{r+1}^c(\bar\alpha)\succ0,
    \qquad
    S_{r+1}^c(\bar\alpha)\succ0
\end{align*}
throughout the accepted segment. Finally, the first two equations of
\eqref{higher_centring_direction} imply that the primal and dual affine constraints are preserved. Thus,
$(X_{r+1}^c,y_{r+1}^c,S_{r+1}^c)$ is strictly primal-dual feasible.
\end{proof}
\noindent
Next, we show that the damped corrector steps preserve the complementarity gap.
\begin{lemma}
\label{corrector_gap}
Let $(X(t),y(t),S(t))$ be the reparametrized weighted central path
through the current corrected iterate at $t_0$, with weight
\begin{align*}
    W
    :=
    \frac{
        X(t_0)S(t_0)+S(t_0)X(t_0)
    }{2t_0^\rho}.
\end{align*}
Assume that
\begin{align*}
    \frac{\langle X(t_0),S(t_0)\rangle}{n}
    =
    t_0^\rho
\end{align*}
and that $p\ge\rho$. Let $(\hat X_{t_0}(t),\hat y_{t_0}(t),\hat S_{t_0}(t))$ be the $p$th-order Taylor predictor defined in~\eqref{Taylor_Approx_Central_Path}. Then
\begin{equation}
\label{eq:predictor_exact_gap}
    \frac{
        \langle
            \hat X_{t_0}(t),
            \hat S_{t_0}(t)
        \rangle
    }{n}
    =
    t^\rho.
\end{equation}
In particular, for $(X_0^c,y_0^c,S_0^c)=\bigl(
        \hat X_{t_0}(t_+),
        \hat y_{t_0}(t_+),
        \hat S_{t_0}(t_+)\bigr)$, we have
\begin{align*}
    \frac{\langle X_0^c,S_0^c\rangle}{n}
    =
    t_+^\rho.
\end{align*}
\noindent
Further, let $(X_r^c,y_r^c,S_r^c)$ be any corrector iterate generated
by~\eqref{damped_corrector}-\eqref{step_length}. Then we have 
\begin{align*}
    \frac{\langle X_r^c,S_r^c\rangle}{n}
    =
    t_+^\rho
    \qquad
    \text{for every corrector iterate } r\ge0.
\end{align*}
\end{lemma}

\begin{proof}
By the definition of $W$,
\begin{align*}
    \trace(W)
    =
    \frac{
        \trace
        \bigl(
            X(t_0)S(t_0)+S(t_0)X(t_0)
        \bigr)
    }{2t_0^\rho}
    =
    \frac{\langle X(t_0),S(t_0)\rangle}{t_0^\rho}
    =
    n.
\end{align*}
Hence, along the weighted central path $\xi_W$, the centrality condition 
\begin{align*}
X(t)S(t)+S(t)X(t)=2t^\rho W
\end{align*}
implies
\begin{align*}
    \langle X(t),S(t)\rangle
    =
    t^\rho\operatorname{Tr}(W)
    =
    nt^\rho.
\end{align*}
\noindent
The Taylor predictor is primal-dual feasible, since the zeroth-order term is feasible and all higher-order derivatives satisfy the homogeneous primal and dual feasibility equations. Therefore,
\begin{align*}
    \langle
        \hat X_{t_0}(t),
        \hat S_{t_0}(t)
    \rangle
    =
    \langle
        C,\hat X_{t_0}(t)
   \rangle
    -
    b^T\hat y_{t_0}(t).
\end{align*}
Since the right-hand side is a linear functional of the predictor
variables, it is precisely the $p$th-order Taylor polynomial at $t_0$ of
\begin{align*}
    \langle C,X(t)\rangle-b^T y(t)
    =
    \langle X(t),S(t)\rangle
    =
    nt^\rho.
\end{align*}
Because $p\ge\rho$, the $p$th-order Taylor polynomial of $nt^\rho$ is the polynomial itself. Thus,
\begin{align*}
    \langle
        \hat X_{t_0}(t),
        \hat S_{t_0}(t)
    \rangle
    =
    nt^\rho,
\end{align*}
which proves \eqref{eq:predictor_exact_gap}.

\vspace{5px}
\noindent
We now consider one damped corrector step. Set
\begin{align*}
    \tau:=t_+^\rho,
    \qquad
    \nu_r
    :=
    \frac{\langle X_r^c,S_r^c\rangle}{n}.
\end{align*}
Taking the trace of the third equation in~\eqref{higher_centring_direction} gives
\begin{align*}
    \left\langle X_r^c,\Delta S_r\right\rangle
    +
    \left\langle S_r^c,\Delta X_r\right\rangle
    =
    n(\tau-\nu_r).
\end{align*}
Furthermore, the first two equations in~\eqref{higher_centring_direction} imply
\begin{align*}
    \mathcal A(\Delta X_r)=0,
    \qquad
    \Delta S_r=-\mathcal A^*(\Delta y_r),
\end{align*}
and hence
\begin{align*}
    \langle\Delta X_r,\Delta S_r\rangle
    =
    -\left\langle
        \Delta X_r,\mathcal A^*(\Delta y_r)
    \right\rangle
    =
    -\left\langle
        \mathcal A(\Delta X_r),\Delta y_r
    \right\rangle
    =
    0.
\end{align*}
Using
\[
    X_{r+1}^c
    =
    X_r^c+\alpha_r\Delta X_r,
    \qquad
    S_{r+1}^c
    =
    S_r^c+\alpha_r\Delta S_r,
\]
we obtain
\begin{align*}
    \frac{\langle X_{r+1}^c,S_{r+1}^c\rangle}{n}
    &=
    \nu_r
    +
    \frac{\alpha_r}{n}
    \left(
        \langle X_r^c,\Delta S_r\rangle
        +
        \langle S_r^c,\Delta X_r\rangle
    \right)
    +
    \frac{\alpha_r^2}{n}
    \langle\Delta X_r,\Delta S_r\rangle\\
    &=
    \nu_r+\alpha_r(\tau-\nu_r)\\
    &=
    (1-\alpha_r)\nu_r+\alpha_r\tau.
\end{align*}
Therefore, if $\nu_r=\tau$, then
\begin{align*}
    \frac{\langle X_{r+1}^c,S_{r+1}^c\rangle}{n}
    =
    \tau
    =
    t_+^\rho.
\end{align*}
Since the predictor point satisfies
\begin{align*}
    \frac{\langle X_0^c,S_0^c\rangle}{n}
    =
    t_+^\rho,
\end{align*}
the conclusion follows by induction on $r$.
\end{proof}
\noindent
Now, we show that only finitely many damped corrector steps are needed to ensure the desired centrality tolerance.
\begin{lemma}[Finite restoration of centrality]
\label{lem:finite_corrector}
Recall $t_+>0$ from~\eqref{mu_reduction}, assume that $p\ge\rho$, and let
\begin{align*}
    (X_0^c,y_0^c,S_0^c)
    =
    \bigl(
        \hat X_{t_0}(t_+),
        \hat y_{t_0}(t_+),
        \hat S_{t_0}(t_+)
    \bigr)
\end{align*}
be the predictor point, where $\vartheta_{t_+}(X_0^c,S_0^c)\le 2\beta$, and $\beta < \frac14$. Then the damped corrector iteration~\eqref{damped_corrector}-\eqref{step_length} reaches
$\mathcal N_{t_+}(\beta)$ after finitely many inner iterations, i.e., there exists $r_*<\infty$ such that
\begin{align*}
    \gamma_{r_*}=\vartheta_{t_+}(X_{r^*}^c,S_{r^*}^c)\le\beta.
\end{align*}
\end{lemma}

\begin{proof}
Set $\tau:=t_+^\rho$. From Lemma~\ref{corrector_gap} we have
\begin{align*}
    \frac{\langle X_r^c,S_r^c\rangle}{n}
    =
    \tau
    \qquad
    \text{for every corrector iterate }r.
\end{align*}
Further, Lemma~\ref{corrector_back_to_neighborhood} implies that every accepted corrector iterate is strictly primal-dual feasible, remains in
$\mathcal N_{t_+}(2\beta)$, and satisfies $\gamma_{r+1}<\gamma_r$ whenever $\gamma_r>\beta$. Now, it remains to show that the strict decrease cannot continue indefinitely
while $\gamma_r>\beta$. To this end, consider
\begin{align*}
\mathcal K_\tau
:=
\left\{
\begin{array}{l|l}
(X,y,S)
&
\begin{array}{l}
\langle A_i,X\rangle=b_i,\quad i=1,\ldots,m,\\
\displaystyle \sum_{i=1}^m y_iA_i+S=C,\\
X,S\succ0,\quad \langle X,S\rangle=n\tau,\\
\beta\le\vartheta_{t_+}(X,S)\le2\beta
\end{array}
\end{array}
\right\}.
\end{align*}
We first show that $\mathcal K_\tau$ is compact. Let $(\bar X,\bar y,\bar S)$ be any fixed strictly primal-dual feasible
point (see Assumption~\ref{IPM_Cond}). For every primal-dual feasible $(X,y,S) \in \kappa_{\tau}$, we have
\begin{align*}
    \langle\bar S,X\rangle
    +
    \langle\bar X,S\rangle
    =
    \langle X,S\rangle
    +
    \langle\bar X,\bar S\rangle = n\tau +  \langle\bar X,\bar S\rangle,
\end{align*}
which implies
\begin{align*}
    \lambda_{\min}(\bar S)\operatorname{Tr}(X)
    +
    \lambda_{\min}(\bar X)\operatorname{Tr}(S)
    \le
    n\tau+\langle\bar X,\bar S\rangle.
\end{align*}
Since $\bar X,\bar S\succ0$ and $X,S\succeq0$, this yields uniform Frobenius norm bounds on
$X$ and $S$. Furthermore, since $A_1,\ldots,A_m$ are linearly independent, the boundedness of $S$ implies the boundedness of $y$. 

\vspace{5px}
\noindent
To prove closedness of $\mathcal{K}_{\tau}$, let a sequence in $\mathcal K_\tau$ converge to
$(X,y,S)$. The affine feasibility equations and
$\langle X,S\rangle=n\tau$ pass to the limit, and $X,S\succeq0$.
Furthermore, $\beta \le \vartheta_{t_+}(X,S)\le2\beta<\frac12$, and hence $XS+SX\succ0$. This also implies $X,S \succ 0$, because otherwise there would exist $v\ne0$ with $Xv=0$, which
would give $v^T(XS+SX)v=0$. Hence $(X,y,S)\in\mathcal K_\tau$, and thus
$\mathcal K_\tau$ is compact.

\vspace{5px}
\noindent
By Proposition~\ref{lem:weighted_AHO_well_defined}, for each $(X,y,S)\in\mathcal K_\tau$, the \AHO~centering system~\eqref{higher_centring_direction} is
non-singular, because $\mathcal K_\tau\subseteq\mathcal N_{t_+}(2\beta)$.
Therefore, by the Implicit Function Theorem, its solution $(\Delta X,\Delta y,\Delta S)$ depends continuously on $(X,y,S)$, that also implies 
\begin{align*}
    \delta(X,y,S)
    :=
    \frac{\|\Delta X\|_F\|\Delta S\|_F}{\tau}
\end{align*}
is continuous on $\mathcal K_\tau$. Therefore, there exists a finite $D_\tau\ge0$ such that
\begin{align*}
    \delta(X,y,S)\le D_\tau
    \qquad
    \forall (X,y,S)\in\mathcal K_\tau.
\end{align*}
\noindent
If $D_\tau=0$, then Lemma~\ref{corrector_back_to_neighborhood} gives
$\gamma_{r+1}=0$ whenever $\gamma_r>\beta$, and the corrector phase terminates. Now, suppose that $D_\tau>0$. As long as
$\gamma_r>\beta$, Lemma~\ref{corrector_back_to_neighborhood} yields either
\begin{align*}
    \gamma_{r+1}\le\frac{\gamma_r}{2},
\end{align*}
or
\begin{align*}
    \gamma_{r+1}
    \le
    \gamma_r
    \left(
        1-\frac{\gamma_r}{4\delta_r}
    \right)
    \le
    \gamma_r
    \left(
        1-\frac{\beta}{4D_\tau}
    \right),
\end{align*}
which can be written as $\gamma_{r+1}\le q_\tau\gamma_r$, with
\begin{align*}
    q_\tau
    :=
    \max\left\{
        \frac12,\,
        1-\frac{\beta}{4D_\tau}
    \right\}
    <1.
\end{align*}
Since $\gamma_0\le2\beta$, it follows that $\gamma_r\le2\beta q_\tau^r$, and thus there exists a finite $r_*$ such that $2q_\tau^{r_*}\le1$, and then $\gamma_{r_*}\le\beta$. This completes the proof.
\end{proof}
\noindent
As a result of Lemmas~\ref{positive_definite_argument}-\ref{lem:finite_corrector}, the corrected point
\begin{align*}   
(X^+,y^+,S^+)
    :=
    (X_{r_*}^c,y_{r_*}^c,S_{r_*}^c)
\end{align*}
is strictly primal-dual feasible, and satisfies
\begin{align*}
    (X^+,y^+,S^+)\in\mathcal N_{t_+}(\beta), \quad     \frac{\langle X^+,S^+\rangle}{n}
    =
    t_+^\rho.
\end{align*}
\noindent
Further, the new iterate $(X^+,y^+,S^+)$ also leads to a positive definite weight matrix 
\begin{align*}
W_+:=\big(X^+ S^+ + S^+ X^+\big)/2t_+^{\rho},
\end{align*}
which defines a new weighted central path passing through the current solution.

\vspace{5px}
\noindent
\paragraph{\textbf{Numerical algorithm}}
Algorithm~\ref{alg:higher-order-IPM-imp} is a numerical realization, in which the higher-order directions and the corrector systems are computed using numerical linear algebra. For the theoretical analysis, $t_+$ denotes the exact boundary parameter defined by~\eqref{mu_reduction}; in the numerical implementation of Section~\ref{sec:numerical_demonstrations}, it is approximated by the geometric line search described in Remark~\ref{rem:line-search}.

\begin{algorithm}[]
\caption{Higher-order \IPM~for \SDO}
\label{alg:higher-order-IPM-imp}
\begin{algorithmic}

\Statex \textbf{Input:} \parbox[t]{0.87\linewidth}{
The degree $\rho$ of reparametrization, the order $p$ of approximation,
the neighborhood parameter $0 < \beta < \frac{1}{4}$, the target tolerance $\zeta>0$, an initial solution $(X_0,y_0,S_0)$ with $t_0^{\rho} = \la X_0, S_0\ra/n$ such that $\vartheta_{t_0}(X_0,S_0) \le \beta$.}

\Statex \textbf{Initialize:}
Set $k = 0$ 

\While{$\frac{\langle X_k,S_k\rangle}{n}>\zeta$}

    \State
    $t_k\gets
    \left(
    \frac{\langle X_k,S_k\rangle}{n}
    \right)^{1/\rho}$

    \State
    $W_k\gets
    \frac{X_kS_k+S_kX_k}{2t_k^\rho}$

    \For{$j=1,\ldots,p$}
        \State
        $(X^{(j)},y^{(j)},S^{(j)})
        \gets$
        Solve \eqref{HOSys} with $t=t_k$
    \EndFor

    \State
    $t_+\gets$
    Determine the boundary parameter from
    \eqref{mu_reduction}

    \State
    \textbf{Predictor:} $\big(
    \hat X_{t_k}(t_+),
    \hat y_{t_k}(t_+),
    \hat S_{t_k}(t_+)
    \big)
    \gets$
    Compute \eqref{Taylor_Approx_Central_Path}
    at $t=t_+$
     \State
    $(X_0^c,y_0^c,S_0^c)\gets$  $\big(
    \hat X_{t_k}(t_+),
    \hat y_{t_k}(t_+),
    \hat S_{t_k}(t_+)
    \big)$
    \State
Set $r = 0$ 
\While{$\vartheta_{t_+}(X_r^c,S_r^c)>\beta$}

    \State
    \textbf{Corrector:} $(X_{r+1}^c,y_{r+1}^c,S_{r+1}^c)\gets$
   Solve \eqref{higher_centring_direction}
    starting from
    $(X_r^c,y_r^c,S_r^c)$ at $t=t_+$ and then compute~\eqref{damped_corrector} using~\eqref{step_length}
      \State
    $r\gets$  $r+1$
\EndWhile
  \State
    $(X_{k+1},y_{k+1},S_{k+1})\gets$  $(X_r^c,y_r^c,S_r^c)$
  \State $k \gets$ $k+1$  
\EndWhile

\Statex \textbf{Return:} $(X_k,y_k,S_k)$

\end{algorithmic}
\end{algorithm}

\vspace{5px}
\noindent
In the following sections, we analyze the local convergence and iteration complexity of Algorithm~\ref{alg:higher-order-IPM-imp}. 

\section{Proof of the main results}\label{Proof_of_main_results}
In this section, we prove Theorems~\ref{Algorithm_Convergence}-\ref{iteration_complexity_result}, i.e., the local convergence and the iteration complexity of Algorithm~\ref{alg:higher-order-IPM-imp}.

\subsection{Local convergence}\label{sec:local_convergence}
The higher-order predictor is constructed from weighted central paths whose weight varies with the iterates. Accordingly, point-wise differentiability of each weighted path is not sufficient for the convergence analysis; a uniform control of the predictor residual along the relevant family of weights is required. We formalize this property next.

\begin{definition}[Uniform predictor regularity]\label{predictor_regularity}
Let $\rho\in\mathbb Z_{>0}$ and $p\geq 1$, and let
$\mathcal W\subseteq\mathbb S^n_{++}$.
We say that the family of reparametrized weighted central paths is
\emph{uniformly $(p,\rho)$-predictor regular on $\mathcal W$} if there
exist constants $\tilde t>0$ and $c>0$ such that, for every
$W\in\mathcal W$, every $0<t_0\leq\tilde t$, and every
$0<t\leq t_0$, the $p$th-order Taylor predictor satisfies
\begin{align*}
\left\|
\hat X(t)\hat S(t)
+\hat S(t)\hat X(t)
-2t^\rho W
\right\|_F
\leq
c|t-t_0|^{p+1}.
\end{align*}
\end{definition}
\noindent
We now prove the super-linear convergence by deriving bounds on a feasible solution of~\eqref{mu_reduction}.

\begin{proof}[Proof of Theorem~\ref{Algorithm_Convergence}]
Let
\begin{align*}
h:=t-t_0,
\qquad
X^{(0)}:=X(t_0),
\qquad
S^{(0)}:=S(t_0).
\end{align*}
By the Taylor expansion~\eqref{Taylor_Approx_Central_Path}, we have
\begin{align*}
\hat X(t)
=
\sum_{j=0}^p \frac{h^j}{j!}X^{(j)}(t_0),
\qquad
\hat S(t)
=
\sum_{j=0}^p \frac{h^j}{j!}S^{(j)}(t_0).
\end{align*}
Therefore,
\begin{align}
\hat X(t)\hat S(t)
+\hat S(t)\hat X(t)
={}&
\sum_{r=0}^{2p} h^r
\sum_{\ell=\max\{0,r-p\}}^{\min\{p,r\}}
\frac{
X^{(\ell)}S^{(r-\ell)}
+
S^{(\ell)}X^{(r-\ell)}
}{
\ell!(r-\ell)!
},
\label{Taylor_product}
\end{align}
where all the derivatives on the right-hand side are evaluated at
$t_0$. For every $0\leq r\leq p$, the coefficient of $h^r$ in
\eqref{Taylor_product} is
\begin{align*}
\frac{1}{r!}
\sum_{\ell=0}^r
\binom{r}{\ell}
\left(
X^{(\ell)}S^{(r-\ell)}
+
S^{(\ell)}X^{(r-\ell)}
\right).
\end{align*}
By~\eqref{HOSys}, this is precisely
\begin{align*}
\frac{1}{r!}
\frac{d^r}{dt^r}
\left(2t^\rho W\right)\bigg|_{t=t_0}.
\end{align*}
Hence, when $p\geq\rho$,
\begin{align}
\hat X(t)\hat S(t)
+\hat S(t)\hat X(t)
=
2t^\rho W+R_p(t),
\label{Taylor_product_remainder}
\end{align}
where
\begin{align*}
R_p(t)
:=
\sum_{r=p+1}^{2p} h^r
\sum_{\ell=r-p}^{p}
\frac{
X^{(\ell)}S^{(r-\ell)}
+
S^{(\ell)}X^{(r-\ell)}
}{
\ell!(r-\ell)!
}.
\end{align*}
Here we used the fact that, since $p\geq\rho$, the $p$th-order
Taylor expansion of the polynomial $2t^\rho W$ about $t_0$ is exact.

\vspace{5px}
\noindent
Now, we restrict the reparametrized weighted central paths to
$(0,\tilde t)$ for some $\tilde t>0$. By the assumption, there exist a
constant $c>0$, independent of $t$ and $t_0$, such that
\begin{align*}
\|R_p(t)\|
\leq
c|t-t_0|^{p+1},
\qquad
0<t\leq t_0\leq\tilde t.
\end{align*}
It follows from~\eqref{Taylor_product_remainder} that
\begin{align*}
\left\|
\frac{
\hat X(t)\hat S(t)
+\hat S(t)\hat X(t)
}{2t^\rho}
-W
\right\|
\leq
\frac{c|t-t_0|^{p+1}}{2t^\rho}.
\end{align*}

Further,
\begin{align}\label{triangle_ineq}
   \bigg\|\frac{\hat{X}(t)\hat{S}(t) + \hat{S}(t)\hat{X}(t)}{2t^{\rho}} - I_n\bigg\| &\le \bigg\|\frac{\hat{X}(t)\hat{S}(t) + \hat{S}(t)\hat{X}(t)}{2t^{\rho}} - W\bigg\| + \|W - I_n\| \nonumber\\
   &\le \bigg\|\frac{\hat{X}(t)\hat{S}(t) + \hat{S}(t)\hat{X}(t)}{2t^{\rho}} - W\bigg\| + \beta,
\end{align}
where the last inequality follows from $(X(t_0),y(t_0),S(t_0)) \in \mathcal{N}(\beta)$. Since the right-hand side of~\eqref{triangle_ineq} is a monotonically decreasing continuous function of $t$ over $(0,t_0]$ taking values in $[0,\infty)$, by the Intermediate Value Theorem there exists a $\bar{t} < t_0$ such that 
\begin{align*}
&\frac{(t_0-\bar{t})^{p+1}}{2\bar{t}^{\rho}}c = \beta,\\
&\frac{(t_0-t)^{p+1}}{2t^{\rho}}c \le \beta, \quad \forall t \in [\bar{t}, t_0]. \nonumber
\end{align*}
Then it follows from~\eqref{mu_reduction} and~\eqref{triangle_ineq} that $\bar{t}$ is a feasible solution for~\eqref{mu_reduction}, and thus $t_+ \le \bar{t}$. Moreover, from
\begin{align*}
\frac{c(t_0-\bar t)^{p+1}}{2\bar t^\rho}=\beta
\end{align*}
we obtain
\begin{align*}
\frac{2\beta\bar t^\rho}{ct_0^{p+1}}
=
\left(\frac{t_0-\bar t}{t_0}\right)^{p+1}
\leq1,
\end{align*}
and therefore
\begin{align*}
t_+
\leq\bar t
\leq
\left(\frac{c}{2\beta}\right)^{1/\rho}
t_0^{(p+1)/\rho}.
\end{align*}
Since $p\geq\rho$, we have $(p+1)/\rho>1$. In particular, for
\[
t_0
\leq
\left(\frac{2\beta}{c}\right)^{1/(p+1-\rho)},
\]
the predictor step yields a strict reduction in $t$, and
\begin{align*}
t_+
=
O\left(t_0^{(p+1)/\rho}\right).
\end{align*}
By definition, $\mu_k=t_0^\rho$ and by Lemma~\ref{corrector_gap}, $\mu_{k+1}=t_+^\rho$. Hence
\begin{align*}
\mu_{k+1}
=
O\left(
\left(\mu_k\right)^{(p+1)/\rho}
\right),
\end{align*}
which proves $Q$-super-linear convergence of the complementarity gaps
with order
\begin{align*}
\gamma=\frac{p+1}{\rho}.
\end{align*}
\end{proof}
\noindent
The predictor regularity condition in Theorem~\ref{Algorithm_Convergence} is formulated directly in terms of the Taylor residuals of the weighted central paths. We next record several sufficient conditions under which
this assumption is automatically satisfied. We begin with a standard joint-smoothness condition, and then specialize it to the strictly complementary setting and to semi-algebraic families of weighted paths.
\begin{corollary}[Joint smoothness]\label{joint_smoothness}
Suppose that $p \ge \rho$, and for some $\tilde t>0$, the mapping
\begin{align*}
(t,W)\longmapsto \xi_W(t^\rho)
\end{align*}
extends to a $\mathcal{C}^p$-smooth mapping on
$[0,\tilde t]\times\mathcal W_{\beta}$, where $\mathcal W_{\beta}$ is compact. Then the family is uniformly $(p,\rho)$-predictor regular on $\mathcal W_{\beta}$.
\end{corollary}

\begin{proof}
Since $\mathcal W_\beta$ is compact and the reparametrized weighted
central path is $\mathcal{C}^p$-smooth on
$[0,\tilde t]\times\mathcal W_\beta$, there exists a constant $M>0$ such that
\begin{align*}
\max_{0\leq j\leq p}
\left\{
\|X^{(j)}(t,W)\|,
\|S^{(j)}(t,W)\|
\right\}
\leq M
\end{align*}
for all $(t,W)\in[0,\tilde t]\times\mathcal W_\beta$. Let $h:=t-t_0$. By definition,
\begin{align*}
\hat X(t)
=
\sum_{j=0}^p
\frac{h^j}{j!}X^{(j)}(t_0,W),
\qquad
\hat S(t)
=
\sum_{j=0}^p
\frac{h^j}{j!}S^{(j)}(t_0,W).
\end{align*}
Therefore,
\begin{align*}
&\hat X(t)\hat S(t)
+\hat S(t)\hat X(t)\\
&\quad =
\sum_{r=0}^{2p} h^r
\sum_{\substack{0\leq j\leq p\\0\leq r-j\leq p}}
\frac{
X^{(j)}(t_0,W)S^{(r-j)}(t_0,W)
+
S^{(j)}(t_0,W)X^{(r-j)}(t_0,W)
}{
j!(r-j)!
}.
\end{align*}
For $0\leq r\leq p$, the inner sum is precisely
\begin{align*}
\frac{1}{r!}
\frac{d^r}{dt^r}
\left(
X(t,W)S(t,W)+S(t,W)X(t,W)
\right)\bigg|_{t=t_0}.
\end{align*}
Since the weighted central path satisfies
\begin{align*}
X(t,W)S(t,W)+S(t,W)X(t,W)=2t^\rho W,
\end{align*}
we obtain
\begin{align*}
\sum_{r=0}^{p}
\frac{h^r}{r!}
\frac{d^r}{dt^r}
\left(2t^\rho W\right)\bigg|_{t=t_0}
=
2t^\rho W,
\end{align*}
where the last equality follows from $p\geq\rho$. Consequently,
\begin{align*}
&\hat X(t)\hat S(t)
+\hat S(t)\hat X(t)
-2t^\rho W\\
&\quad =
\sum_{r=p+1}^{2p} h^r
\sum_{j=r-p}^{p}
\frac{
X^{(j)}(t_0,W)S^{(r-j)}(t_0,W)
+
S^{(j)}(t_0,W)X^{(r-j)}(t_0,W)
}{
j!(r-j)!
}.
\end{align*}
Without loss of generality, we may take $\tilde t\leq1$.
Since $0<t\leq t_0\leq\tilde t$, we have $|h|\leq1$, and hence
$|h|^r\leq |h|^{p+1}$ for every $r\geq p+1$. Using the
sub-multiplicativity of the matrix norm and the uniform bound above,
\begin{align*}
&
\left\|
\hat X(t)\hat S(t)
+\hat S(t)\hat X(t)
-2t^\rho W
\right\|\\
&\quad\leq
2M^2 |h|^{p+1}
\sum_{r=p+1}^{2p}
\sum_{j=r-p}^{p}
\frac{1}{j!(r-j)!}.
\end{align*}
Thus, defining
\begin{align*}
c:=
2M^2
\sum_{r=p+1}^{2p}
\sum_{j=r-p}^{p}
\frac{1}{j!(r-j)!},
\end{align*}
we obtain
\begin{align*}
\left\|
\hat X(t)\hat S(t)
+\hat S(t)\hat X(t)
-2t^\rho W
\right\|
\leq
c|t-t_0|^{p+1},
\end{align*}
where $c$ is independent of $t_0$ and $W$.
\end{proof}
\noindent
The preceding criterion immediately recovers the classical strictly complementary regime.
\begin{corollary}[Strict complementarity]
Suppose that the strict complementarity condition holds. Then, for
every $p\geq1$, the \AHO~weighted central paths are uniformly
$(p,1)$-predictor regular on $\mathcal W_\beta$. Consequently,
\begin{align*}
\mu_{k+1}=O(\mu_k^{p+1}).
\end{align*}
In particular, choosing $p=1$ yields quadratic convergence of the
normalized complementarity gaps.
\end{corollary}

\begin{proof}
Under strict complementarity, the \AHO~weighted central paths depend analytically jointly on the path parameter and the weight matrix (see~\cite{PS04}). Hence, after possibly restricting the path parameter to a sufficiently small interval, the mapping
\begin{align*}
(\mu,W)\longmapsto \xi_W(\mu)
\end{align*}
is jointly $\mathcal{C}^p$-smooth on the compact set
$[0,\tilde\mu]\times\mathcal W_\beta$ for every $p\geq1$.
Corollary~\ref{joint_smoothness} therefore implies uniform $(p,1)$-predictor regularity.
The result now follows from Theorem~\ref{Algorithm_Convergence} with $\rho=1$.
\end{proof}
\noindent
We next give a variant adapted to the semi-algebraic structure of the weighted path family established earlier in the paper.
\begin{corollary}[Cell-wise algebraic regularity]
Let $\mathcal W\subset\mathbb S^n_{++}$ be compact and suppose that
$\mathcal W$ admits a finite semi-algebraic partition
\begin{align*}
\mathcal W=C_1\cup\cdots\cup C_r
\end{align*}
such that, for every $i$, the mapping
\begin{align*}
(t,W)\longmapsto\xi_W(t^\rho),
\qquad
(t,W)\in(0,\tilde t]\times C_i,
\end{align*}
admits a $\mathcal{C}^p$ extension to an open neighborhood of
\begin{align*}
[0,\tilde t]\times\overline{C_i}.
\end{align*}
Then the family is uniformly $(p,\rho)$-predictor regular on $\mathcal W$.
\end{corollary}

\begin{proof}
For each $i$, the set
$[0,\tilde t]\times\overline{C_i}$ is compact. Hence, by the assumed $\mathcal{C}^p$ extension, all derivatives up to order $p$ of the reparametrized
weighted central path are uniformly bounded on this set. Therefore, by the proof of Corollary~\ref{joint_smoothness}, there exists $c_i>0$ such that
\begin{align*}
\left\|
\hat X(t)\hat S(t)
+\hat S(t)\hat X(t)
-2t^\rho W
\right\|
\leq
c_i|t-t_0|^{p+1}
\end{align*}
for every $W\in C_i$ and $0<t\leq t_0\leq\tilde t$.
Since the partition is finite, taking $c:=\max_{1\leq i\leq r}c_i$ proves uniform $(p,\rho)$-predictor regularity on $\mathcal W$.
\end{proof}

\begin{remark}[Non-admissible reparametrizations]
\label{rem:nonadmissible_rho}
Analyticity is stronger than what is needed to construct a predictor
of a fixed finite order. To see this, consider a Puiseux term $c\mu^\alpha$ appearing in a coordinate of a weighted central path. Under the reparametrization $\mu=t^\rho$, this term becomes $ct^{\rho\alpha}$. If $\rho$ is admissible, all such exponents are non-negative integers, and derivatives of arbitrary order are available at $t=0$. If $\rho$ is not admissible, some exponent $\rho\alpha$ remains non-integral, and the reparametrized path generally possesses only finitely many bounded derivatives at $t=0$. A $p$th-order predictor may still be constructed as long as the derivatives through order $p$ exist and satisfy the uniform bounds required in the preceding analysis.

\vspace{5px}
\noindent
There is, however, an additional restriction for super-linear convergence. The proof of
Theorem~\ref{Algorithm_Convergence} requires $p\geq\rho$, and yields the convergence order $(p+1)/\rho$. Thus, a non-admissible reparametrization can provide the same
super-linear guarantee only if it retains at least $\rho$ bounded
derivatives. This need not occur. For example, if a relevant coordinate contains a non-zero Puiseux term
$c\mu^{1/N}$, then after reparametrization it becomes $ct^{\rho/N}$. If $\rho$ is not a multiple of $N$, this term provides only
$\lfloor\rho/N\rfloor$ bounded derivatives. Since, for $N>1$, $\left\lfloor\frac{\rho}{N}\right\rfloor<\rho$, the requirements for constructing a $p$th-order predictor and for obtaining the super-linear convergence guarantee cannot simultaneously be satisfied in this case.

\vspace{5px}
\noindent
Therefore, an admissible reparametrization removes the finite-differentiability
obstruction altogether. Among all admissible choices, $\rho=N$ is
the smallest and, for a fixed Taylor order $p$, gives the largest
guaranteed convergence order.
\end{remark}
\noindent
We now prove super-linear convergence of simultaneously vanishing blocks of $X_k$ and $S_k$ with respect to the optimal partition of~\SDO. To this end, we need to define the notion of the optimal partition (see e.g.,~\cite[Lemma~2.3]{Kl02} or~\cite{MT19}).
\begin{definition}\label{optimal_partition}
Let $(X^*,y^*,S^*) \in \ri(\mathcal{P}^* \times \mathcal{D}^*)$, where $\ri(\cdot)$ denotes the relative interior of a convex set, and let $\mathcal{B}:=\mathcal{R}(X^*)$, $\mathcal{N}:=\mathcal{R}(S^*)$, and $\mathcal{T}\!:=\big(\mathcal{R}\big(X^*\big) + \mathcal{R}\big(S^*\big)\big)^{\perp}$, where $\mathcal{R}(\cdot)$ is the column space of a matrix. Then the 3-tuple $\big(\mathcal{B}, \mathcal{T},\mathcal{N}\big)$ is called the \textit{optimal partition} of $(\mathrm{P})-(\mathrm{D})$.  

\vspace{5px}
\noindent
The optimal partition $\big(\mathcal{B}, \mathcal{T},\mathcal{N}\big)$
is independent of the choice of $(X^*,y^*,S^*)$~\cite[Lemma~2.3(i)]{Kl02}. By the complementarity condition, $X^*$ and $S^*$ commute, and thus they have a common eigenvector basis $Q$, whose columns can be partitioned into $Q=[Q_{\mathcal{B}},Q_{\mathcal{T}},Q_{\mathcal{N}}]$ according to the partitions $\mathcal{B}$, $\mathcal{T}$, and $\mathcal{N}$.
\end{definition}

\begin{proof}[Proof of Theorem~\ref{convergence_of_vanishing_blocks}]
An orthogonal transformation of $(X^*,y^*,S^*) \in \ri(\mathcal{P}^* \times \mathcal{D}^*)$ with respect to $Q$ reveals the optimal partition as
\begin{align*}
Q^T X^* Q=\begin{pmatrix} U_{X^*}  \  & 0  \  & 0 \\ 0  \  & 0  \  & 0 \\ 0  \  & 0  \  & 0 \end{pmatrix}, \qquad Q^T S^* Q=\begin{pmatrix} 0  & 0  & 0 \\ 0  &  0  &  0 \\ 0  &  0  &  U_{S^*} \end{pmatrix}, 
\end{align*}
where $U_{X^*},U_{S^*} \succ 0$. Now, it follows from~\cite[Theorem~4.3]{MT19} that  when $k$ is large enough, there exist positive constants $\kappa,\kappa'$ independent of $t_{k+1}$ (but dependent on $n$ and $\beta$) and a positive exponent $\gamma$ such that
\begin{align*}
\kappa' t_{k+1}^{(1-\gamma)\rho}\le \|Q^T_{\mathcal{T}} X_{k+1}Q_{\mathcal{T}}\| \le \kappa t_{k+1}^{\gamma \rho},
\end{align*}
where
\begin{align*}
   \frac{1}{2^{n-1}} \le \gamma \le \frac12.
\end{align*}
Further, using~\cite[Theorem~4.3]{MT19} again and uniform boundedness, there exist positive constants $\eta,\eta'$ independent of $t_{k+1}$ (dependent on $n$ and a condition number defined in~\cite[Page~232]{MT19}) such that 
\begin{align*}
\eta' t_{k+1}^{ \rho} \le \|Q^T_{\mathcal{N}} X_{k+1}Q_{\mathcal{N}}\| \le \eta t_{k+1}^{ \rho}.
\end{align*}
\noindent
Now, for sufficiently large $k$, it follows from Theorem~\ref{Algorithm_Convergence} that 
\begin{align*}
\|Q^T_{\mathcal{T}} X_{k+1}Q_{\mathcal{T}}\| \le \kappa t_{k+1}^{\gamma \rho} &\le \kappa \Bigg(\frac{c}{2\beta}\Bigg)^{\gamma}t_k^{(p+1)\gamma} \\
&\le \kappa \Bigg(\frac{c}{2\beta}\Bigg)^{\gamma} \Big(\frac{1}{\kappa'}\Big)^{\frac{(p+1)\gamma}{\rho(1-\gamma)}} \|Q^T_{\mathcal{T}} X_{k}Q_{\mathcal{T}}\|^{\frac{(p+1)\gamma}{\rho(1-\gamma)}},\\
&\le \kappa \Bigg(\frac{c}{2\beta}\Bigg)^{\gamma} \Big(\frac{1}{\kappa'}\Big)^{\frac{(p+1)\gamma}{\rho(1-\gamma)}} \|Q^T_{\mathcal{T}} X_{k}Q_{\mathcal{T}}\|^{\frac{p+1}{\rho(2^{n-1}-1)}},
\end{align*}
where the last inequality follows from $\|Q^T_{\mathcal{T}} X_{k}Q_{\mathcal{T}}\| \to 0$ and $\gamma/(1-\gamma) \ge 1/(2^{n-1}-1)$.
Analogously, we can derive
\begin{align*}
\|Q^T_{\mathcal{N}} X_{k+1}Q_{\mathcal{N}}\| \le \eta t_{k+1}^{\rho} &\le \eta \Bigg(\frac{c}{2\beta}\Bigg)t_k^{(p+1)} \\
&\le \eta \Bigg(\frac{c}{2\beta}\Bigg) \Big(\frac{1}{\eta'}\Big)^{\frac{p+1}{\rho}} \|Q^T_{\mathcal{N}} X_{k}Q_{\mathcal{N}}\|^{\frac{p+1}{\rho}}.
\end{align*}
The proof is similar for $Q^T_{\mathcal{B}} S_kQ_{\mathcal{B}}$ and $Q^T_{\mathcal{T}} S_kQ_{\mathcal{T}}$. 
\end{proof}

\subsection{Iteration complexity}\label{sec:pitercomp} 
We first derive a quantitative reduction bound for the normalized complementarity gap on any region where the uniform
$(p,\rho)$-predictor regularity estimate holds. We then combine this bound with the super-linear recurrence of Theorem~\ref{Algorithm_Convergence} to obtain a doubly logarithmic local outer iteration bound.

\begin{proposition}[Predictor reduction]
\label{prop:predictor-reduction}
Suppose that the family of reparametrized weighted central paths is
uniformly $(p,\rho)$-predictor regular on $\mathcal{W}_\beta$ with
constants $\tilde t>0$ and $c>0$. Let $p\geq\rho$, and consider an iterate with $0<t_k\leq\tilde t$. Let $\bar t_k\in(0,t_k)$ be the unique solution of
\begin{align*}
\frac{c(t_k-\bar t_k)^{p+1}}{2\bar t_k^\rho}=\beta,
\end{align*}
as in the proof of Theorem~\ref{Algorithm_Convergence}. If $\alpha_k:=\bar t_k/t_k$, then
\begin{align}\label{strong_bound}
\alpha_k
\leq
\frac{1}{
1+
\left(
\frac{2\beta}
     {c\,t_k^{p+1-\rho}}
\right)^{1/\rho}
}.
\end{align}
Further, we have
\begin{align*}
\mu_{k+1}
\leq
\left[
1+
\left(
\frac{2\beta}
     {c\,t_k^{p+1-\rho}}
\right)^{1/\rho}
\right]^{-\rho}
\mu_k.
\end{align*}
\end{proposition}

\begin{proof}
Apply the argument in the proof of Theorem~\ref{Algorithm_Convergence} with $t_0=t_k$ and with the weight associated with the current iterate. In particular, inequality~\eqref{triangle_ineq} shows that, if $\bar t_k\in(0,t_k)$ satisfies
\begin{align*}
\frac{c(t_k-\bar t_k)^{p+1}}{2\bar t_k^\rho}=\beta,
\end{align*}
then every $t\in[\bar t_k,t_k]$ is feasible for the predictor problem~\eqref{mu_reduction}. Since the left-hand side is continuous and strictly
decreasing from $+\infty$ to $0$ as $\bar t_k$ increases from $0$ to $t_k$, such $\bar t_k$ is unique. Hence $t_+\le \bar t_k$. Now, set
\begin{align*}
\alpha_k:=\frac{\bar t_k}{t_k}\in(0,1).
\end{align*}
The defining equation for $\bar t_k$ gives
\begin{align*}
\frac{(1-\alpha_k)^{p+1}}{\alpha_k^\rho}
=
\frac{2\beta}{c\,t_k^{p+1-\rho}}.
\end{align*}
Since $p\ge\rho$,
\begin{align*}
\left(\frac{1-\alpha_k}{\alpha_k}\right)^\rho
\ge
\frac{2\beta}{c\,t_k^{p+1-\rho}},
\end{align*}
and therefore
\begin{align*}
\alpha_k
\le
\frac{1}{
1+
\left(
\frac{2\beta}{c\,t_k^{p+1-\rho}}
\right)^{1/\rho}
}.
\end{align*}
\noindent
Finally, Lemma~\ref{corrector_gap} gives $\mu_{k+1}=t_+^\rho$, while by definition
$\mu_k=t_k^\rho$. Thus,
\begin{align*}
\mu_{k+1}
=t_+^\rho
\le \bar t_k^\rho
=\alpha_k^\rho\mu_k,
\end{align*}
and substituting the preceding bound on $\alpha_k$ yields the claimed estimate.
\end{proof}

\begin{remark}[Lower-order predictors]
\label{rem:p-less-rho}
Theorem~\ref{Algorithm_Convergence} requires $p\geq\rho$ and hence does not provide a super-linear convergence guarantee when $p<\rho$. We leave the behavior of the method in this regime for future investigation.
\end{remark}

\begin{remark}[Global reduction versus polynomial complexity]
\label{rem:global-reduction}
Proposition~\ref{prop:predictor-reduction} shows that, whenever the uniform
$(p,\rho)$-predictor regularity estimate is valid, Algorithm~\ref{alg:higher-order-IPM-imp} makes a quantitative reduction in the normalized complementarity gap. In particular, the reduction becomes increasingly aggressive as $t_k\downarrow0$.

This estimate, however, does not by itself yield the conventional polynomial iteration complexity of path-following \IPM s. Indeed, the bound on the predictor reduction depends on the constant $c$ in the bound~\eqref{strong_bound}. Consequently, the estimate does not provide a uniform reduction of the form
\begin{align*}
\mu_{k+1}
\leq
\left(1-\frac{\delta}{\sqrt n}\right)\mu_k
\end{align*}
with a constant $\delta>0$ independent of the problem data. Thus, Proposition~\ref{prop:predictor-reduction} should be interpreted as a quantitative
progress estimate for the higher-order predictor rather than as a polynomial outer iteration complexity bound.
\end{remark}
\noindent
We now turn to the local regime. When $p\geq\rho$, Theorem~\ref{Algorithm_Convergence} gives the stronger recurrence
\begin{align*}
\mu_{k+1}=O\!\left(\mu_k^{(p+1)/\rho}\right),
\end{align*}
from which we obtain a doubly logarithmic local outer iteration bound.

\begin{proof}[Proof of Theorem~\ref{iteration_complexity_result}]
Recall $d=p+1-\rho \ge 1$. By Proposition~\ref{prop:predictor-reduction},
\begin{align*}
\alpha_k
\le
\frac{1}{
1+\left(
\frac{2\beta}{c t_k^d}
\right)^{1/\rho}
}.
\end{align*}
Since $t_k\leq t_0$ and $d\geq1$,
\begin{align*}
\frac{2\beta}{c\,t_k^d}
\geq
\frac{2\beta}{c\,t_0^d},
\end{align*}
and therefore $\alpha_k\leq q_0$. Now, using Proposition~\ref{prop:predictor-reduction} we get
\begin{align*}
\mu_{k+1}
\leq
\alpha_k^\rho\mu_k
\leq
q_0^\rho\mu_k.
\end{align*}
Inductively, we get $\mu_k\leq q_0^{\rho k}\mu_0$. Hence, $\mu_k\leq\eta$ whenever
\begin{align*}
q_0^{\rho k}\mu_0\leq\eta,
\end{align*}
which gives
\begin{align*}
k
\geq
\frac{\log(\mu_0/\eta)}
     {\rho\log(1/q_0)}.
\end{align*}

\vspace{5px}
\noindent
For the local phase, Theorem~\ref{Algorithm_Convergence} gives
\begin{align*}
\mu_{k+1}
\leq
C\mu_k^q,
\qquad
q:=\frac{p+1}{\rho}>1,
\end{align*}
for all sufficiently small $\mu_k$. Choose $\mu_{\mathrm{loc}}>0$ sufficiently small so that $C\mu_{\mathrm{loc}}^{q-1}<1$. Setting $v_k:=C^{1/(q-1)}\mu_k$ gives
\begin{align*}
v_{k+1}\leq v_k^q.
\end{align*}
Thus, the normalized complementarity gap decreases doubly exponentially with the local outer iteration counter, and the number of additional outer iterations required to reach $\mu_k\leq\zeta$ is
\begin{align*}
O\!\left(
\log\log\frac{1}{\zeta}
\right).
\end{align*}
\end{proof}
\noindent
We can also comment on the costs of predictor and corrector steps per iteration.
\begin{notation}
The symmetric Kronecker product of any two square matrices $K_1$ and $K_2$ is defined as a mapping 
\begin{align*}
(K_1 \otimes_s K_2)\svectorize(H):= \frac12 \svectorize\big(K_2 H K_1^T + K_1HK_2^T\big),
\end{align*}
where $H$ is any symmetric matrix (see e.g.,~\cite{Kl02} for more details).
\end{notation}

\begin{remark}[Cost of higher-order directions]\label{system_solving}
At a fixed outer iterate, all derivative systems in~\eqref{HOSys} have the
same coefficient matrix and differ only in their right-hand sides. More specifically, to obtain higher-order derivatives, we iteratively solve the system \eqref{HOSys} at each iteration of \IPM, where the coefficient matrix is fixed in each iteration. The full Newton system is given by
\begin{equation}\label{sys:FNS}
    \begin{aligned}
        \begin{bmatrix}
         0 & \mathbf{A}^T & \mathcal{I} \\
        \mathbf{A} & 0   & 0 \\
        S \otimes_s I & 0   & X \otimes_s I
        \end{bmatrix}
        \begin{bmatrix}
            {\rm \svectorize}(X^{(j)}) \\
            y^{(j)} \\
            {\rm \svectorize}(S^{(j)})
        \end{bmatrix}
        =
        \begin{bmatrix}
            0 \\
            0 \\
            \svectorize(f(W,\xi(j)))
        \end{bmatrix} \qquad \forall j = 1,\dots, p,
    \end{aligned}
\end{equation}
where $\mathbf{A} = [{\rm \svectorize}(A_1), \dots,  {\rm \svectorize}(A_m)]^T$, and
\begin{equation*}
   f(W,\xi(j)) = \begin{cases} 2\rho (\rho - 1)\cdots (\rho-j+1) t^{\rho - j} W + \xi(j) \ \ &j\le \rho,\\
        \xi(j) \ \ &j>\rho.
        \end{cases}
\end{equation*}
As an alternative, we can build the normal equations system: 
\begin{align}\label{sys:NES}
        \left(\mathbf{A} \mathbf{E}^{-1} \mathbf{F} \mathbf{A}^T \right) y^{(j)} =- \mathbf{A} \mathbf{E}^{-1} \svectorize(f(W,\xi(j))), \qquad \forall j = 1,\dots, p,
\end{align}
where $\mathbf{E} = S \otimes_s I$ and $\mathbf{F} = X \otimes_s I$ (notice that both mappings $\mathbf{E}$ and $\mathbf{F}$ are invertible at the interior solution $(X,y,S)$). Given that $\mathbf{A}\mathbf{E}^{-1} \mathbf{F} \mathbf{A}^T$ is positive definite, we can perform a one-time factorization to solve~\eqref{sys:NES}. We then perform forward-backward substitution for computing higher-order derivatives.

Thus, the additional cost of increasing the predictor order is primarily the cost of these
additional right-hand sides and back-solves, rather than $p$ independent factorizations. Each corrector step similarly requires the
solution of one \AHO~Newton system. We do not, however, obtain here a uniform bound on the number of inner corrector steps per outer
iteration.
\end{remark}

\section{Numerical experiments}\label{sec:numerical_demonstrations}
In this section, we elaborate on the implementation details of Algorithm~\ref{alg:higher-order-IPM-imp}. After fixing $(\rho,p)$ at some pre-determined pair of integers, we initiate Algorithm~\ref{alg:higher-order-IPM-imp} with a feasible solution $(X_0, y_0, S_0) \in \mathcal{N}(\beta)$ with the predictor-neighborhood parameter $\beta = 0.5$. The value $\beta=0.5$ used in the numerical experiments is larger than the sufficient neighborhood parameter required by the theoretical well-definedness argument in Proposition~\ref{lem:weighted_AHO_well_defined}. This is an implementation choice intended to permit larger predictor steps.

\vspace{5px}
\paragraph{\textbf{Predictor step}} The predictor step consists of computing the higher-order derivatives of a reparametrized weighted central path. To that end, we solve~\eqref{sys:FNS} directly. Once the higher-order information is computed, we construct the $p$th-order Taylor approximation of the weighted central
path passing through the current point
$(X_k,y_k,S_k)$.

The theoretical algorithm uses the boundary parameter $t_+$ defined
by~\eqref{mu_reduction}. Rather than solving this
one-dimensional semi-algebraic problem exactly, in our implementation we approximate $t_+$ by a geometric line search. Starting from $t=t_k$, we successively reduce $t$ according to
\begin{align*}
t\leftarrow\alpha_{\rm ls}t,
\qquad
\alpha_{\rm ls}=0.9999,
\end{align*}
and retain the last accepted value before the next trial point leaves $\mathcal N(2\beta)$. Since we use $\beta=0.5$ in the numerical experiments, we additionally
check the semi-definite constraints during the forward-tracking line search. In particular, a trial point is accepted only if it remains in $\mathcal{N}(2\beta)$ and its primal and dual matrices remain
positive definite. Thus, the line search is prevented from overshooting the positive semi-definite cone.

\vspace{5px}
\paragraph{\textbf{Corrector step}} At the point obtained from the predictor step, we perform, possibly multiple, damped centering steps toward the central path according to~\eqref{damped_corrector}-\eqref{step_length}. For this purpose, we set the centrality tolerance, which determines the desired proximity to the central path, to $10^{-1}$. Upon completion of the centering step, the resulting interior solution is used as the base point for constructing the Taylor approximation in the subsequent iteration. Although the
implementation allows for damped corrector steps, all centering steps taken in the experiments reported below were full steps
($\alpha=1$); no damping was required.

\vspace{5px}
\noindent
The above procedure is repeated until the complementarity gap falls below the prescribed tolerance $\zeta = 10^{-10}$. 

\vspace{5px}
\noindent
To evaluate the performance of our higher-order \IPM~with respect to the convergence order, we compared it to MOSEK\footnote{\url{https://www.mosek.com/}} and SDPA\footnote{\url{https://sdpa.sourceforge.net/}} through CVXPY\footnote{\url{https://www.cvxpy.org/}}, and SeDuMi\footnote{\url{https://sedumi.ie.lehigh.edu/}} and SDPT3\footnote{\url{https://www.math.cmu.edu/~reha/sdpt3.html}}. While MOSEK serves as our primary benchmark, we occasionally include additional solvers for further comparison. MOSEK employs a homogeneous self-dual (HSD) embedding together with a Mehrotra-type predictor-corrector scheme. SeDuMi likewise uses a predictor-corrector approach. SDPT3 provides both an HSD-embedded path-following method (SDPT3-HSD) and an infeasible path-following method (SDPT3), while SDPA is based on a Mehrotra-type predictor-corrector primal-dual \IPM. Unless otherwise stated, the gap values reported for SDPT3 in the figures correspond to the original problem formulation, irrespective of whether the infeasible or HSD variant is used. Our primary objective is to compare the asymptotic complementarity gap reduction and outer iteration counts of these
established solvers with those of our higher-order \IPM, rather than wall-clock performance. A higher-order iteration generally requires
more linear algebra than a conventional first-order predictor-corrector
iteration. Our algorithm is implemented in Python, and the source code is available at \url{https://github.com/psm-optimizes/HOIPM4SDO.git}.

\subsection{Initialization}  
Solvers based on HSD embedding are initialized with a point $(X_0, y_0, S_0)$ such that $\frac{\langle X_0,S_0 \rangle}{n}=1$. To facilitate a fair comparison of complementarity gap reduction, we initialize our \IPM~with a point having the same initial gap. Moreover, we set the optimality tolerance to $10^{-10}$ for both our higher-order \IPM~and all competing solvers.

\vspace{5px}
\noindent
The figures in this section illustrate the convergence history of our algorithm, either in isolation or in comparison with other solvers. In each figure, the horizontal axis (\emph{Iteration}) represents the iteration count, while the vertical axis (\emph{Normalized Complementarity Gap (log)}) displays the normalized complementarity gap $\frac{\langle X_k, S_k \rangle}{n}$ at iteration $k$ on a logarithmic scale.

\vspace{5px}
\noindent
For the prescribed maximum order $p$ of the Taylor approximation, we choose $p\in\{\rho,2\rho,3\rho\}$, unless stated otherwise.
In the numerical implementation, the computation of higher-order derivatives may terminate before reaching $p$ if the scaled Taylor coefficient at a given order falls below the threshold $10^{-16}$. Thus, $p$ denotes the prescribed maximum Taylor order, while the effective order used at a particular iteration may be smaller.

\subsection{Modification of the Newton systems}\label{modified_system}
Due to the structure of the right-hand side (RHS) in the linear systems used to compute higher-order derivatives, which contains the factor $\rho(\rho-1)\cdots(\rho-j+1)$, the right-hand side might have a significantly large value, even for moderately large reparametrization exponents. To improve the numerical accuracy of the
linear solves, we employ two scaling strategies. First, we scale the system with $t_0^j$ for computing the $j$th-order derivative. As our \IPM~starts at a complementarity gap of $1$ and proceeds downward, this scaling helps reduce the magnitude of the RHS. Second, we further scale the system using the $1/j!$ coefficient from the Taylor expansion. For convenience, let
\begin{align*}
(\bar X^{(j)},\bar y^{(j)},\bar S^{(j)})
:=
\frac{t_0^j}{j!}
(X^{(j)},y^{(j)},S^{(j)}).
\end{align*}
Specifically, instead of~\eqref{sys:FNS}, we solve
\begin{align*}
\begin{bmatrix}
0 & \mathbf A^T & I\\
\mathbf A & 0 & 0\\
\mathbf E & 0 & \mathbf F
\end{bmatrix}
\begin{bmatrix}
\operatorname{svec}(\bar X^{(j)})\\
\bar y^{(j)}\\
\operatorname{svec}(\bar S^{(j)})
\end{bmatrix}
=
\frac{t_0^j}{j!}
\begin{bmatrix}
0\\
0\\
\svectorize(f(W,\xi^{(j)}))
\end{bmatrix},
\qquad j=1,\ldots,p.
\end{align*}
\noindent
Analogously, instead of the normal equations system~\eqref{sys:NES}, we solve
\begin{align*}
\left(\mathbf A \mathbf E^{-1} \mathbf F \mathbf A^T\right)\bar y^{(j)}
=
-\frac{t_0^j}{j!}
\mathbf A \mathbf E^{-1}\svectorize(f(W,\xi^{(j)})),
\qquad j=1,\ldots,p.
\end{align*}
\noindent
After the scaled systems are solved, the resulting quantities are
multiplied by $j!/t_0^j$ to recover the true derivatives. Together,
these modifications significantly improve the numerical precision of
the computed derivatives by reducing the magnitude of the
right-hand sides of the linear systems.

\subsection{Experiments}
We evaluate our higher-order \IPM~on five sets of benchmark problems. In Section~\ref{ssec:3-ell}, we consider the 3-elliptope problem, whose central path has ramification index $2$. In Section~\ref{ssec:EdK}, we solve three instances from \cite[Example 3.3, p.~58]{Kl02}, for which a lower bound on the global ramification index is known. We then turn to problem classes for which no prior information on the ramification index is available. In Section~\ref{ssec:GenP}, we study a set of generated problems that fail to satisfy the strict complementarity condition. We next consider problems with a positive duality gap. Since the homogeneous self-dual embedding of such problems fails strict complementarity \cite{KRT98}, Section~\ref{ssec:PDG} reports results for solving their self-dual embedding counterparts. Finally, in Section~\ref{ssec:SSC}, we evaluate our higher-order \IPM~on problems that satisfy strict complementarity and therefore do not require reparametrization.

\vspace{5px}
\noindent
For our experiments, the two main input parameters are the candidate reparametrization exponent $\rho$ and the Taylor order $p$. When $\rho$ is admissible and the uniform $(p,\rho)$-predictor regularity assumption of Theorem~\ref{Algorithm_Convergence} holds, choosing $p\geq \rho$ yields a $Q$-order of at least $(p+1)/\rho>1$. Accordingly, we take
\begin{align*}
p\in\{\rho,2\rho,3\rho\}.
\end{align*}
When an admissible reparametrization exponent is not known, we test the candidate values $\rho=2,3,4$, motivated by the heuristic choices
suggested in~\cite [Remark~2]{HLMZ21}. For several of the tested choices of $(\rho,p)$, we observe behavior consistent with super-linear convergence. This speed-up substantially reduces the number of outer iterations required to reach the prescribed tolerance relative to the other
solvers considered in our study.

\subsubsection{Experiments on the 3-elliptope problem}\label{ssec:3-ell}
\begin{figure}
    \centering
    \includegraphics[width=.75\linewidth]{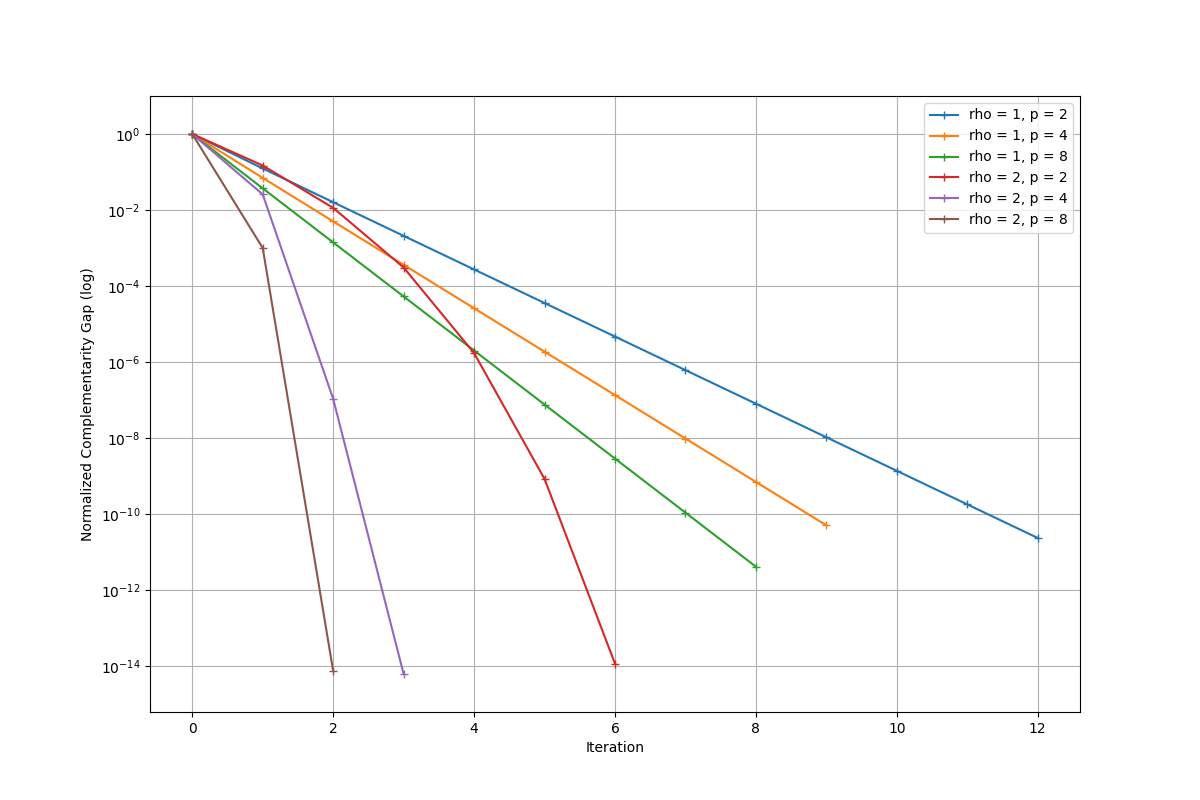}
    \caption{Impact of higher-order derivatives for solving the 3-elliptope problem with and without a reparametrization}
    \label{fig:3ell-rho1vs2Multip}
\end{figure}

\begin{figure}
    \centering
    \includegraphics[width=.75\linewidth]{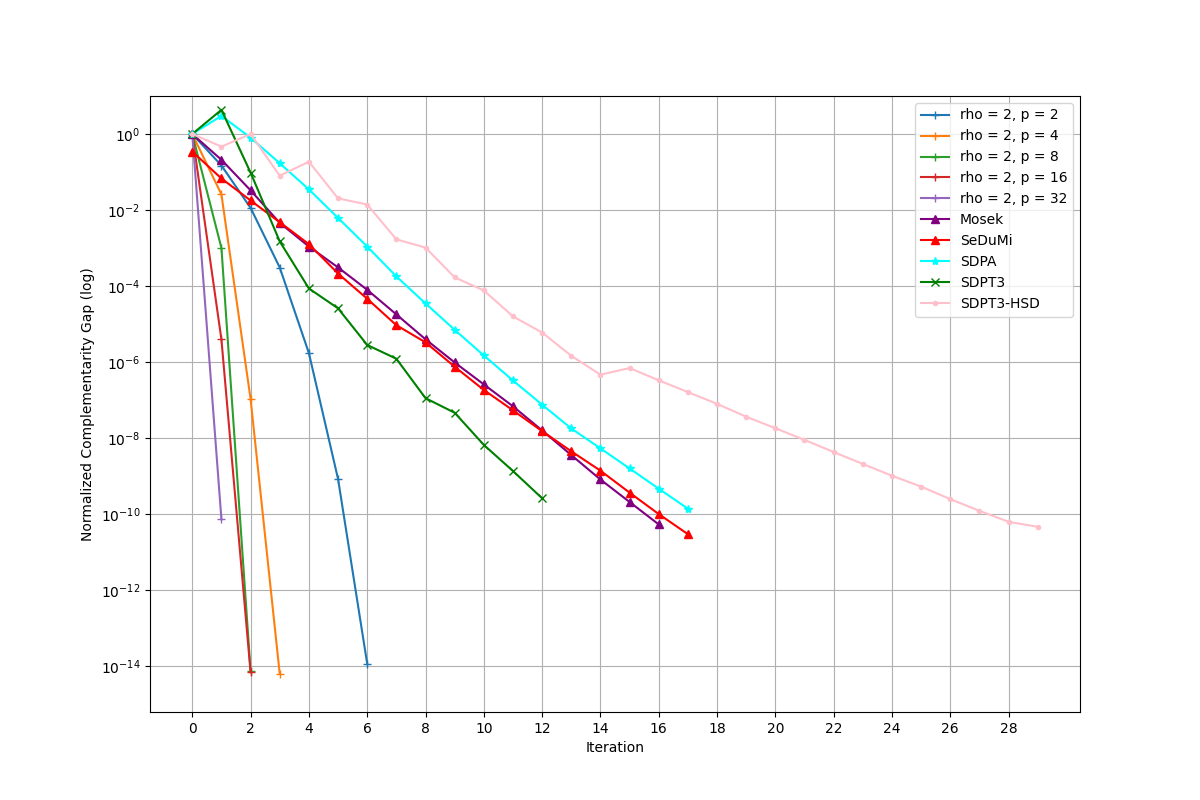}
    \caption{Increasing the order of the Taylor approximation helps reduce the number of iterations}
    \label{fig:3ell-rho2MultipvsSolvers}
\end{figure}

We recall Example~\ref{derivatives_diverge}, whose ordinary central path has ramification
index $2$. Motivated by this value, we compare $\rho=1$ (no reparametrization) with the candidate reparametrization exponent
$\rho=2$. As shown in Figure~\ref{fig:3ell-rho1vs2Multip}, reparametrization is essential for improving the convergence rate: when $\rho = 1$, increasing the order of the Taylor approximation may reduce the number of iterations, but it does not improve the asymptotic convergence rate. In contrast, when $\rho = 2$, higher-order Taylor approximations lead to both a faster convergence rate and a smaller number of iterations. 
Figure~\ref{fig:3ell-rho2MultipvsSolvers} further demonstrates a comparison between our higher-order \IPM~with $\rho = 2$ versus several solvers. We observe that our higher-order \IPM~requires fewer outer iterations and exhibits a steeper reduction in the normalized complementarity gap than the other solvers considered on the 3-elliptope problem. Particularly, with $p=8$, the prescribed tolerance is reached in
two outer iterations. As discussed in Remark~\ref{system_solving}, in addition to the full Newton system, one can exploit the normal equations system, and our algorithm requires one factorization per iteration, followed by $p$ back-solves to compute the higher-order derivatives. This strategy, which uses one factorization followed by multiple back-solves, has also been advocated as a means of reducing the number of \IPM~iterations~\cite[Page~31]{AGMX96}. Although our preference is to keep the order of the Taylor approximation small, allowing more solves per iteration can further decrease the total number of iterations. In fact, as illustrated in Figure~\ref{fig:3ell-rho2MultipvsSolvers}, increasing $p$ sufficiently enables the prescribed tolerance to be reached in a single outer iteration.

\begin{figure}
    \centering
    \includegraphics[width=.75\linewidth]{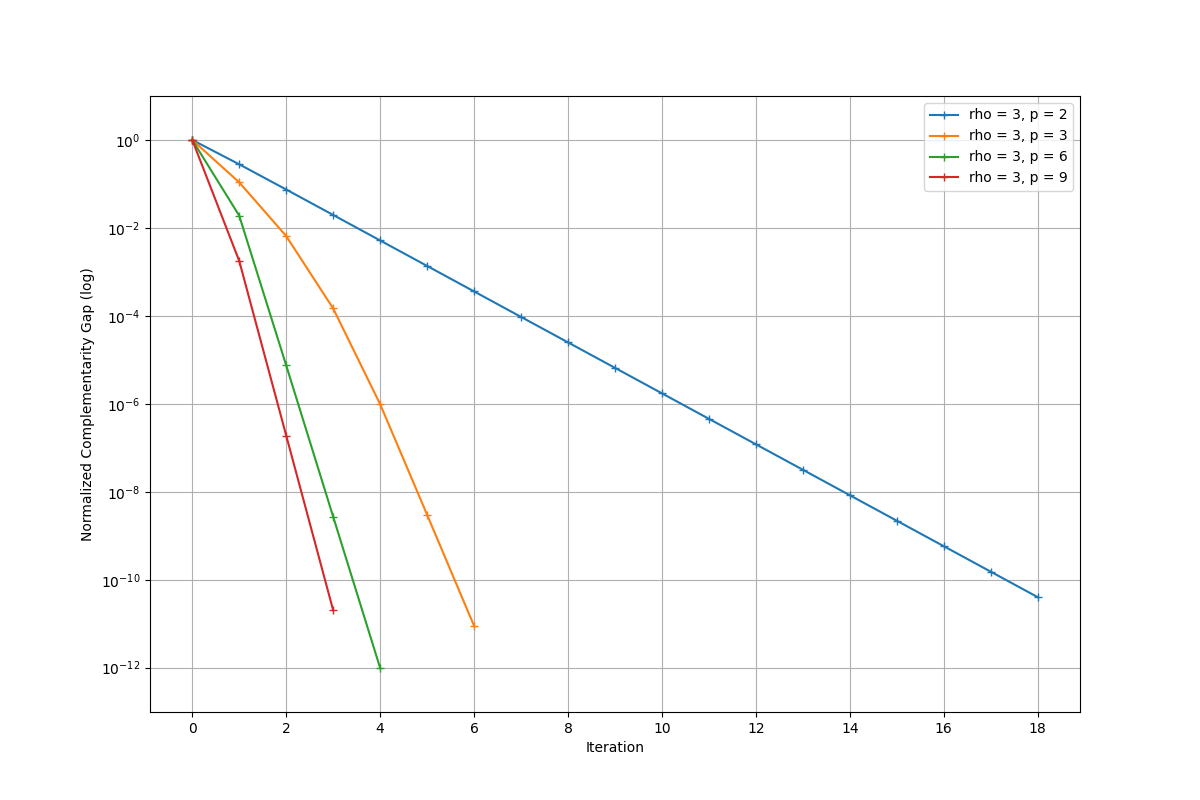 }
    \caption{Using $\rho =3$ for the 3-elliptope problem}
    \label{fig:3ell-rho3}
\end{figure}

\vspace{5px}
\noindent
Finally, recall the Puiseux expansion~\eqref{coordinate_Puiseux_expansion} of the central path near $\mu = 0$. We note that with a reparametrization $\mu \mapsto \mu^3$, we get
\begin{align*}
X_{12}(\mu)=-1+\frac{\sqrt{8}}{8} \mu^{3/2}+\frac{1}{32} \mu^3 -\frac{11\sqrt{8}}{2048} \mu^{9/2}-\frac{3}{512} \mu^{6}-\frac{121\sqrt{8}}{1048576} \mu^{15/2}+ \ldots, 
\end{align*}
whose second derivative does not exist at $\mu=0$. Consequently, choosing $\rho=3$ restores first-order differentiability of the central path, but does not recover analyticity at $\mu=0$. Figure~\ref{fig:3ell-rho3} shows that we get speed-ups even with $\rho =3$. Still, as we approach optimality, the diagnostic corresponding to $\omega=1.05$ eventually ceases to decrease, so the reported
iterates do not provide evidence of a $Q$-order exceeding $1.05$. We examine this behavior in more detail in Section~\ref{sec:convseq}.

\subsubsection{Experiments on a scalable instance}\label{ssec:EdK}
\begin{figure}
    \centering
    \includegraphics[width=.75\linewidth]{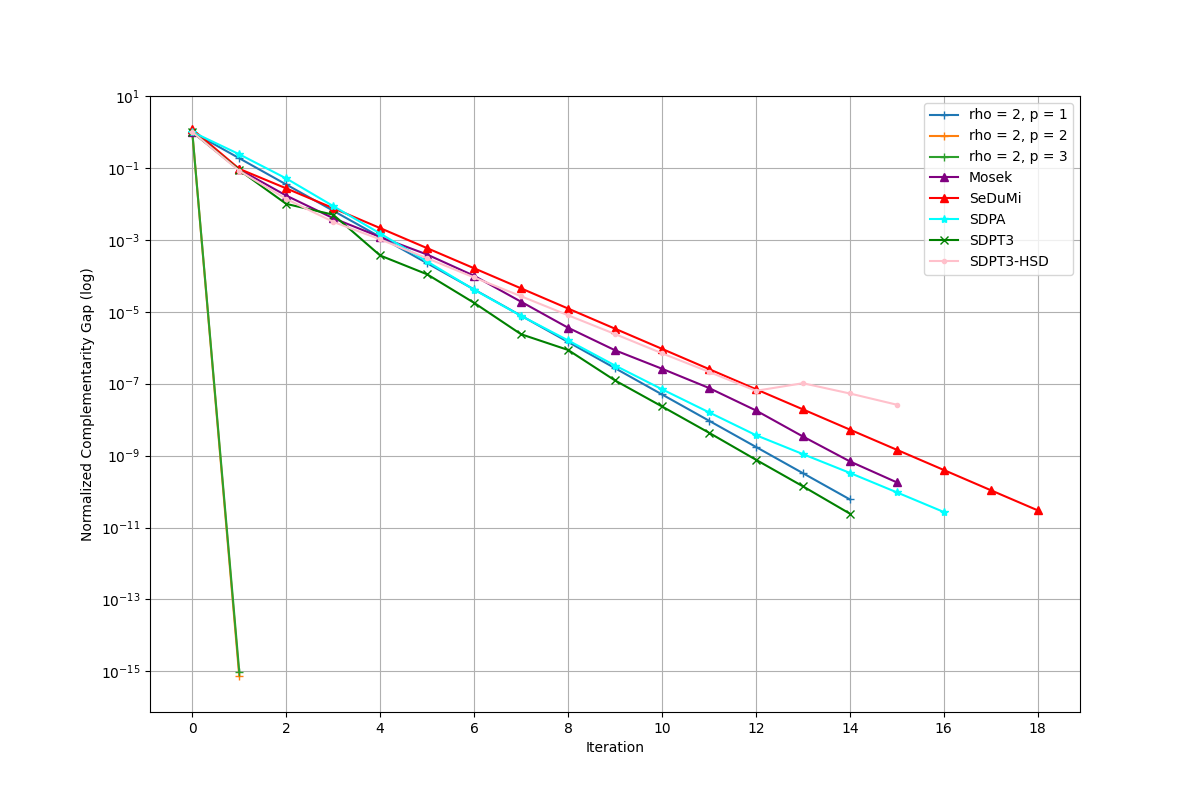}
    \caption{Impact of reparametrization for the scalable problem of size $n=3$}
    \label{fig:EdK3-rho2}
\end{figure}

In this section, we consider the following \SDO\ problem (from \cite{Kl02}) in dual form: 
\begin{align*}
\max\Bigg\{-y_n \mid S=
    \begingroup 
    \setlength{\tabcolsep}{.75pt} 
    \renewcommand{\arraystretch}{.8}
    \begin{pmatrix} 1 & y_1 & y_2 &\ldots & y_{n-1}\\y_1 & y_2 & 0 & \ldots & 0\\y_2 & 0 & y_3 & \ddots & \vdots\\ \vdots & \vdots & \ddots & \ddots & 0\\y_{n-1} & 0 & \ldots & 0 & y_n \end{pmatrix}
    \endgroup \succeq 0\Bigg\},
    \end{align*}
which has a unique optimal solution $y^{*}_i=0$ for $i=1,\ldots,n$. It is easy to check that $y_2(\mu)=c\mu^{2^{-(n-2)}}$ for some constant $c$ independent of $\mu$, which yields the lower bound $2^{n-2}$ on the ramification index of the central path. Since the central path corresponds to $W=I$, this also gives the lower bound $N\geq 2^{n-2}$ on the global ramification index of the weighted family. We begin with the case $n=3$, for which the above bound gives $N\geq 2$. Motivated by this lower bound, we first test the candidate reparametrization exponent $\rho=2$. Figure~\ref{fig:EdK3-rho2} illustrates the effect of the reparametrization with $\rho=2$ for several orders of Taylor approximation. Remarkably, orders $p=2$ and $p=3$ reach the prescribed tolerance in a single outer iteration. In contrast, when $p=1$, our \IPM~exhibits a convergence rate comparable to that of the other solvers, namely, linear convergence. 

\begin{figure}
    \centering
    \includegraphics[width=.75\linewidth]{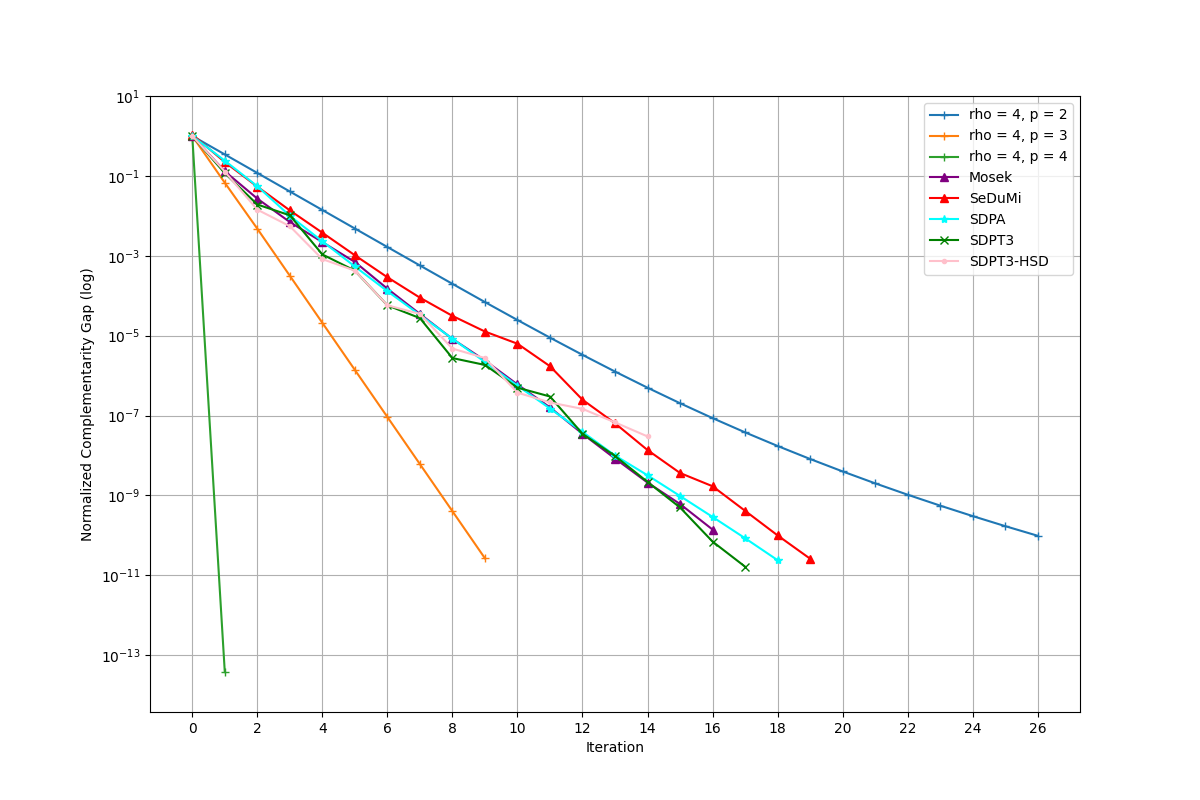}
    \caption{Impact of reparametrization for the scalable problem of size $n=4$ }
    \label{fig:EdK4}
\end{figure}

\begin{figure}
    \centering
    \includegraphics[width=.75\linewidth]{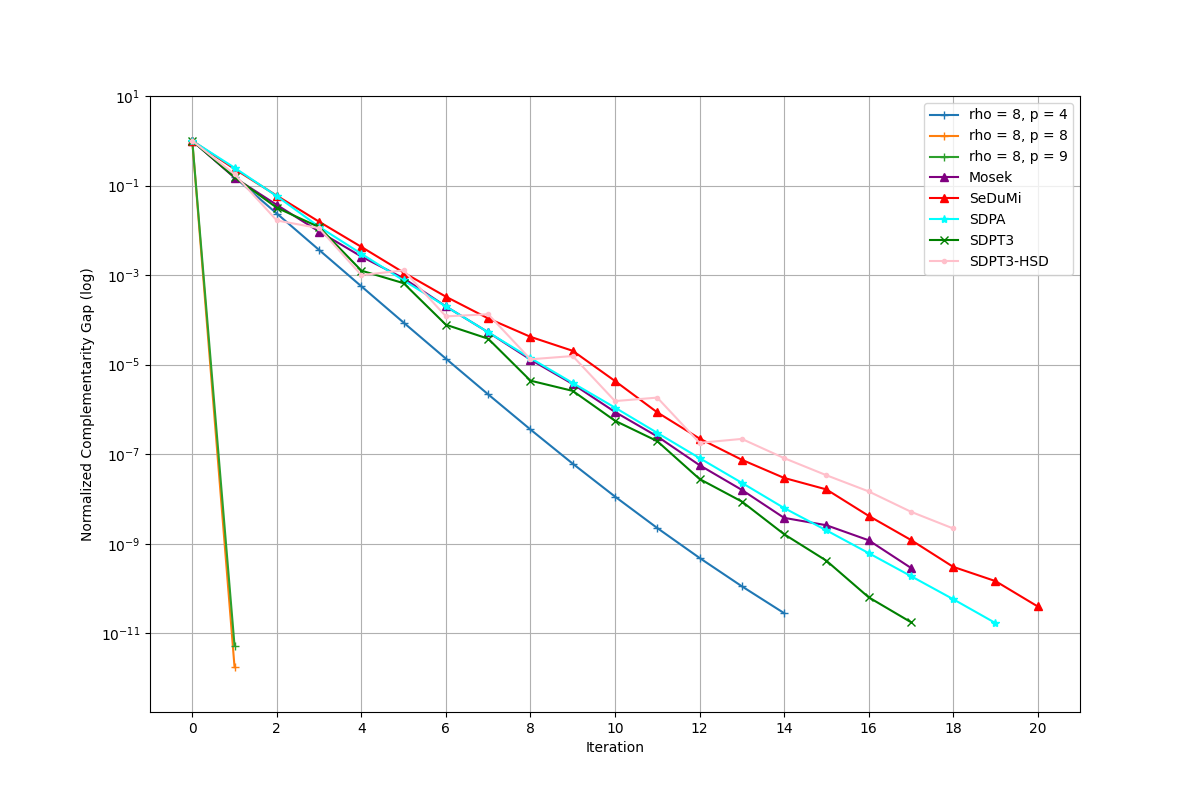}
    \caption{Impact of reparametrization for the scalable problem of size $n=5$}
    \label{fig:EdK5}
\end{figure}

\vspace{5px}
\noindent
We next consider the cases $n=4,5$, for which the above lower bound gives $N\geq4$ and $N\geq8$, respectively. Motivated by these bounds, we test $\rho=4$ for $n=4$ and $\rho=8$ for $n=5$. For $n=4$, the fourth-order Taylor approximation reaches the prescribed tolerance in a single outer iteration (see Figure~\ref{fig:EdK4}), while $p=2,3$ exhibit nearly sub-linear and linear convergence rates,
respectively. For $n=5$, Figure~\ref{fig:EdK5} shows that the prescribed tolerance is reached in a single outer iteration for $p=8,9$. We further observed, in additional experiments not shown here, that
the tested choices $\rho<7$ did not exhibit super-linear convergence.

\subsubsection{Random problems failing strict complementarity}\label{ssec:GenP}
In this section, we extend our experiments to randomly generated instances that violate the strict complementarity condition. The problems are generated using \cite[Algorithm~11]{MFAT24}, which constructs \SDO{} instances using a known strictly feasible solution and an optimal solution. We restrict our attention to problems whose partition $\mathcal{T}$ is non-empty (see Definition~\ref{optimal_partition}) and additionally enforce the primal and dual non-degeneracy conditions of~\cite{AHO97}. We generate a set of \SDO~problems that satisfy both primal and dual non-degeneracy conditions, implying that primal and dual optimal solutions are unique (see~\cite[Theorems~7 and 10]{AHO97}), while failing the strict complementarity condition, i.e., $X^* + S^* \nsucc 0$. The optimal partition configurations for all test instances are summarized in Table~\ref{tab:ND}, where $|\cdot|$ means the dimension of the corresponding subspace of the optimal partition.

\vspace{5px}
\noindent
Since the global ramification index is unavailable for these instances, we test several candidate reparametrization exponents. From this point onward, we restrict our comparisons to MOSEK, owing to its consistently strong performance in the preceding experiments. We also tested SDPA and observed that its performance is at best similar to that of MOSEK.

\begin{figure}[htbp]
\centering

\begin{minipage}[t]{0.48\textwidth}
  \centering
  \includegraphics[width=\linewidth]{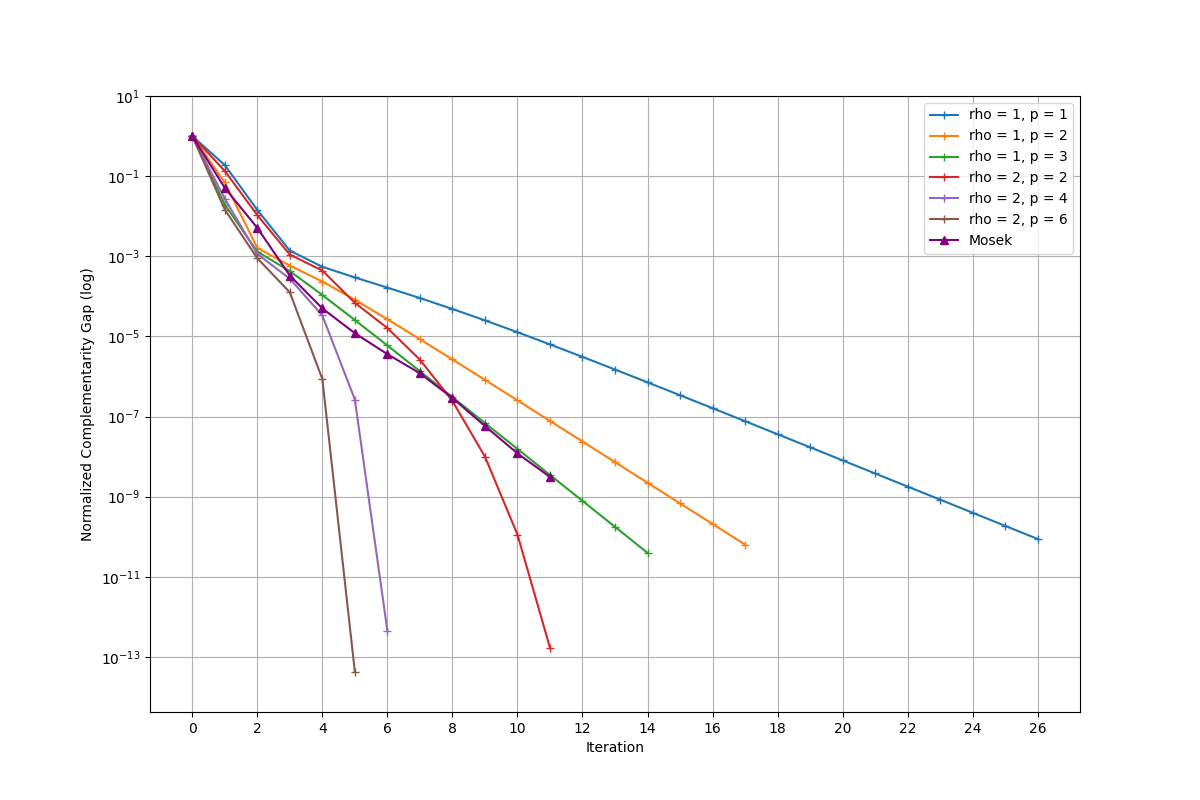}
  \par\vspace{0.3em}\small ND-1:‌ $\rho =$ 1 vs 2
\end{minipage}\hfill
\begin{minipage}[t]{0.48\textwidth}
  \centering
  \includegraphics[width=\linewidth]{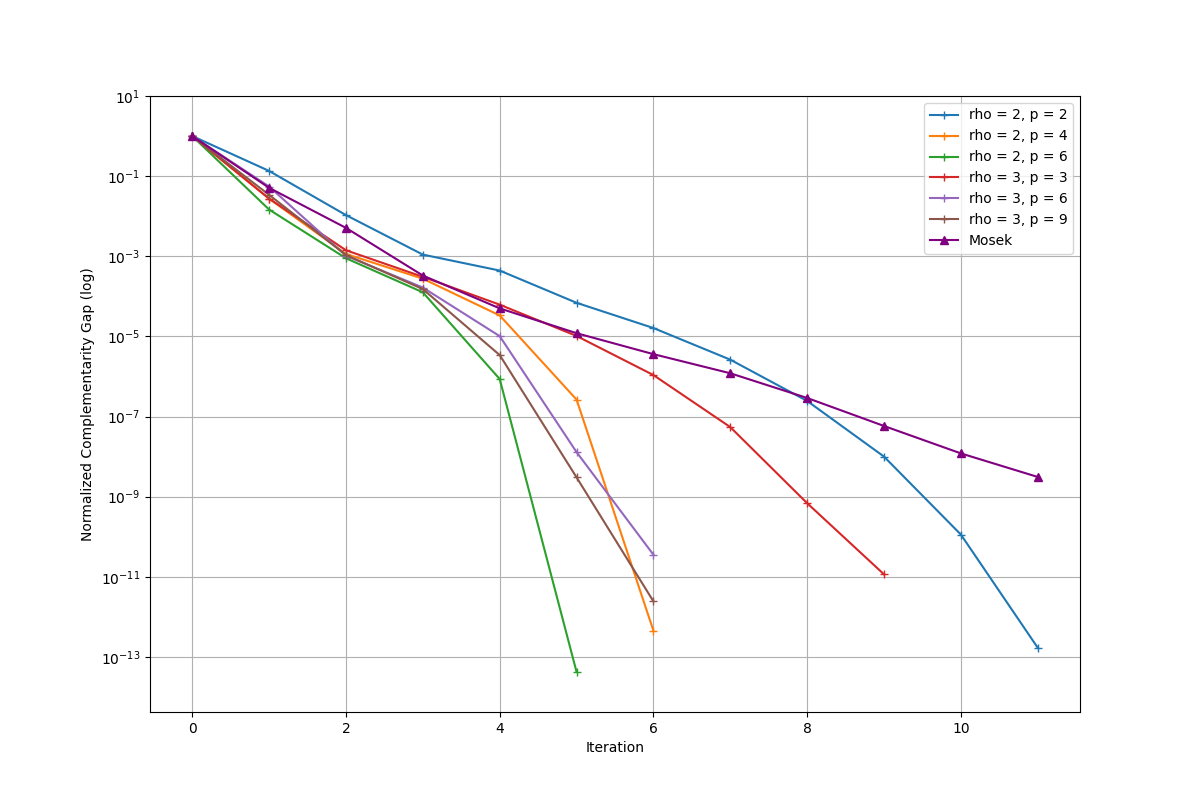}
  \par\vspace{0.3em}\small ND-1:‌ $\rho =$ 2 vs 3
\end{minipage}

\vspace{0.8em}

\begin{minipage}[t]{0.48\textwidth}
  \centering
  \includegraphics[width=\linewidth]{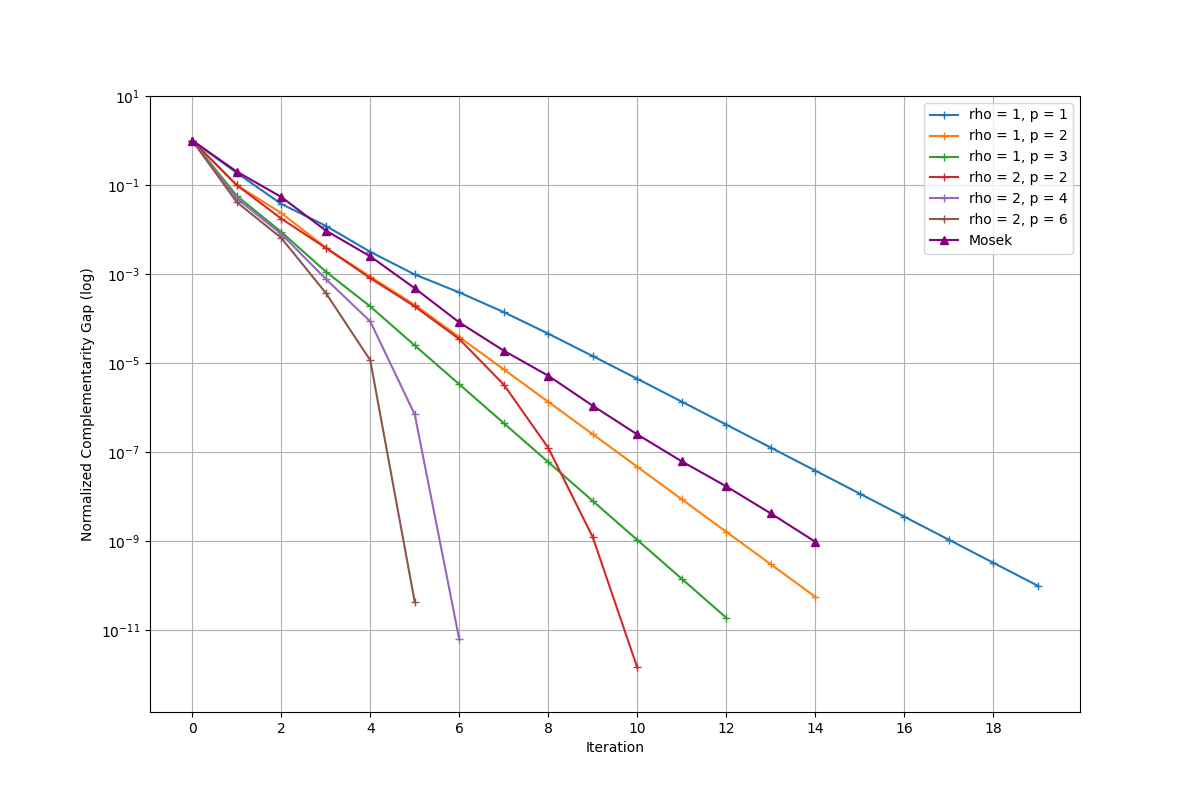}
  \par\vspace{0.3em}\small ND-2:‌ $\rho =$ 1 vs 2
\end{minipage}\hfill
\begin{minipage}[t]{0.48\textwidth}
  \centering
  \includegraphics[width=\linewidth]{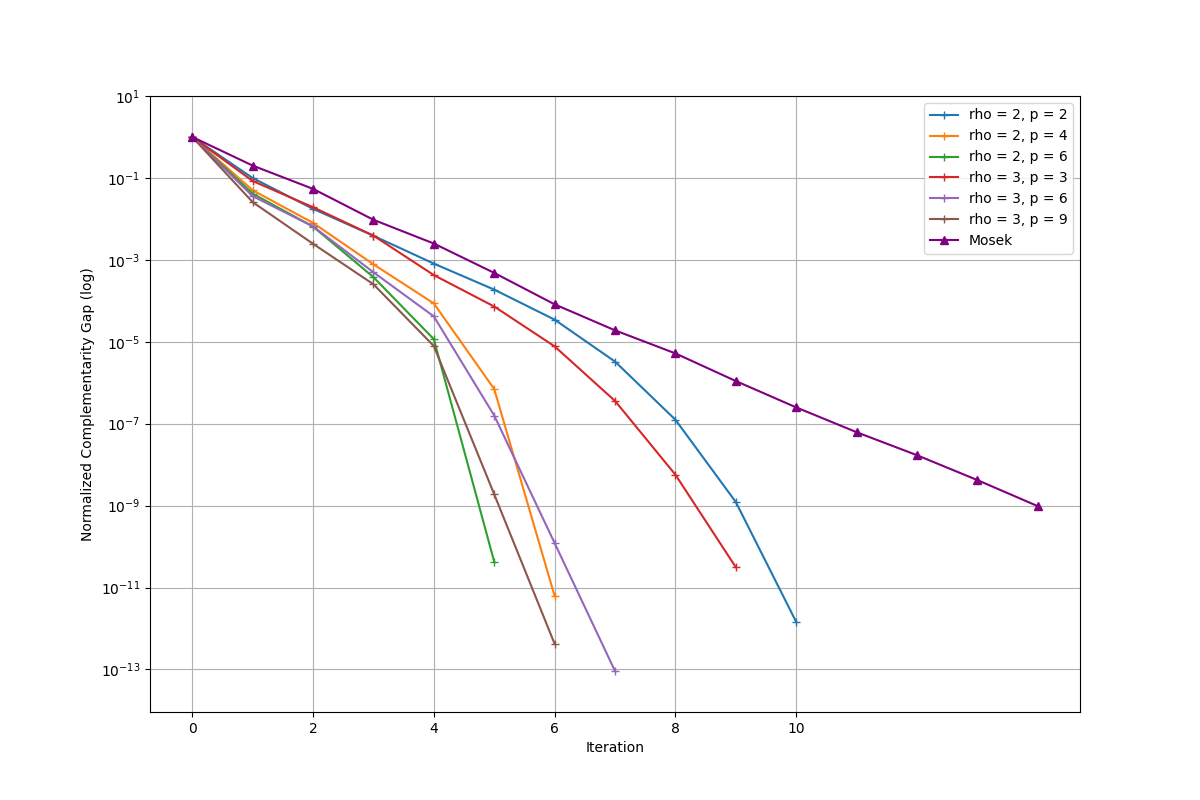}
  \par\vspace{0.3em}\small ND-2:‌ $\rho =$ 2 vs 3
\end{minipage}

\vspace{0.8em}

\begin{minipage}[t]{0.48\textwidth}
  \centering
  \includegraphics[width=\linewidth]{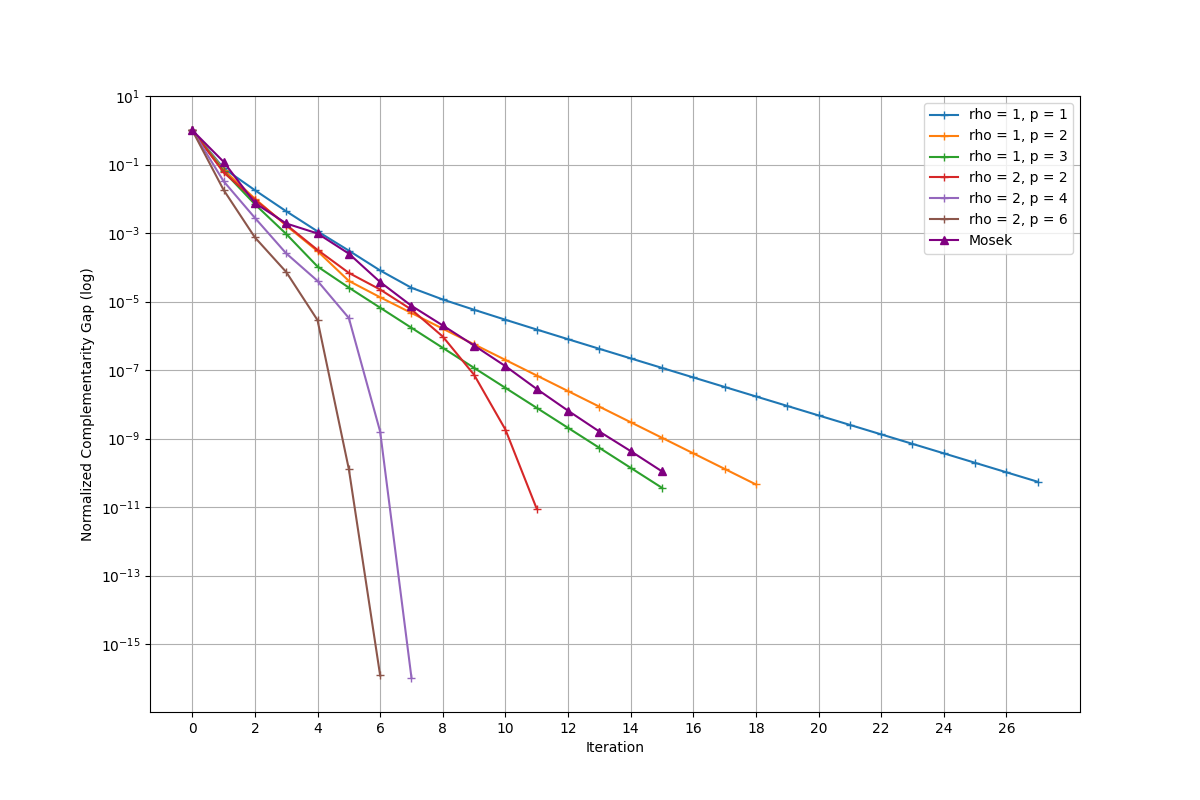}
  \par\vspace{0.3em}\small ND-3:‌ $\rho =$ 1 vs 2
\end{minipage}\hfill
\begin{minipage}[t]{0.48\textwidth}
  \centering
  \includegraphics[width=\linewidth]{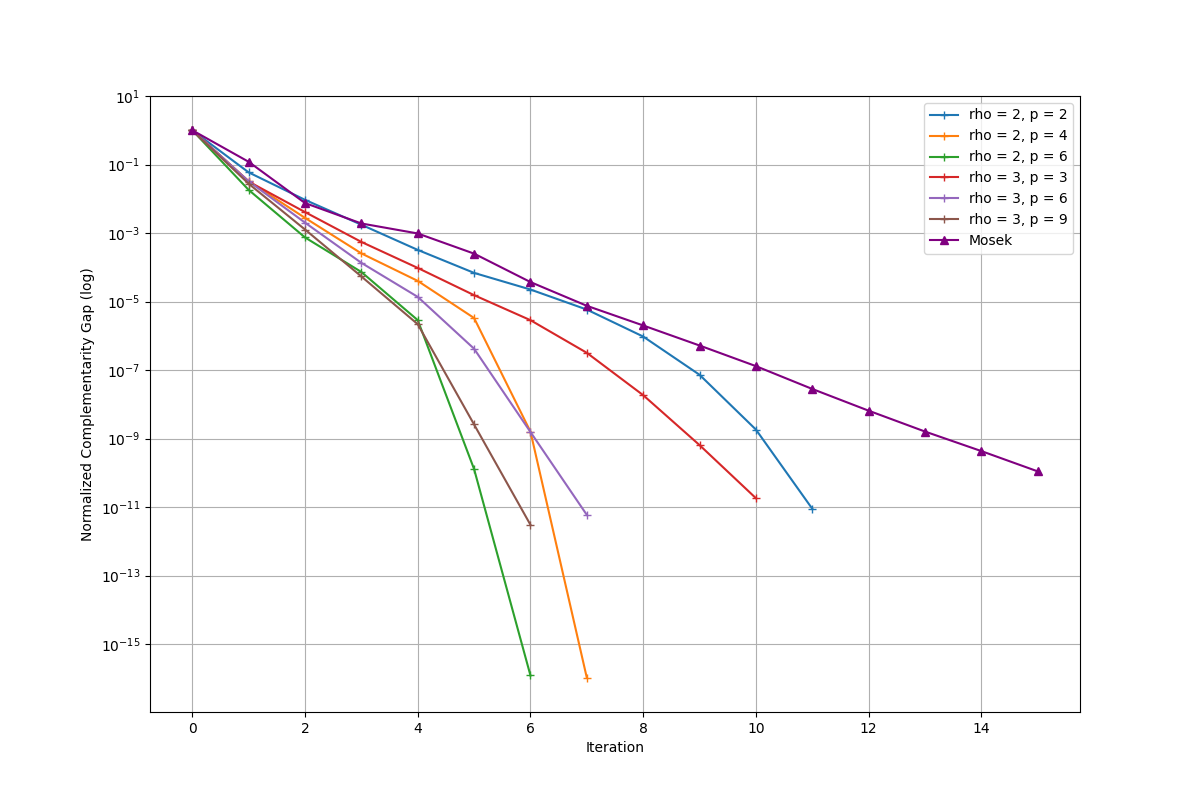}
  \par\vspace{0.3em}\small ND-3:‌ $\rho =$ 2 vs 3
\end{minipage}

\caption{Convergence of random instances with primal-dual non-degeneracy (ND1-ND3)}
\label{fig:NDConv1}
\end{figure}

\begin{figure}[htbp]
\centering

\begin{minipage}[t]{0.48\textwidth}
  \centering
  \includegraphics[width=\linewidth]{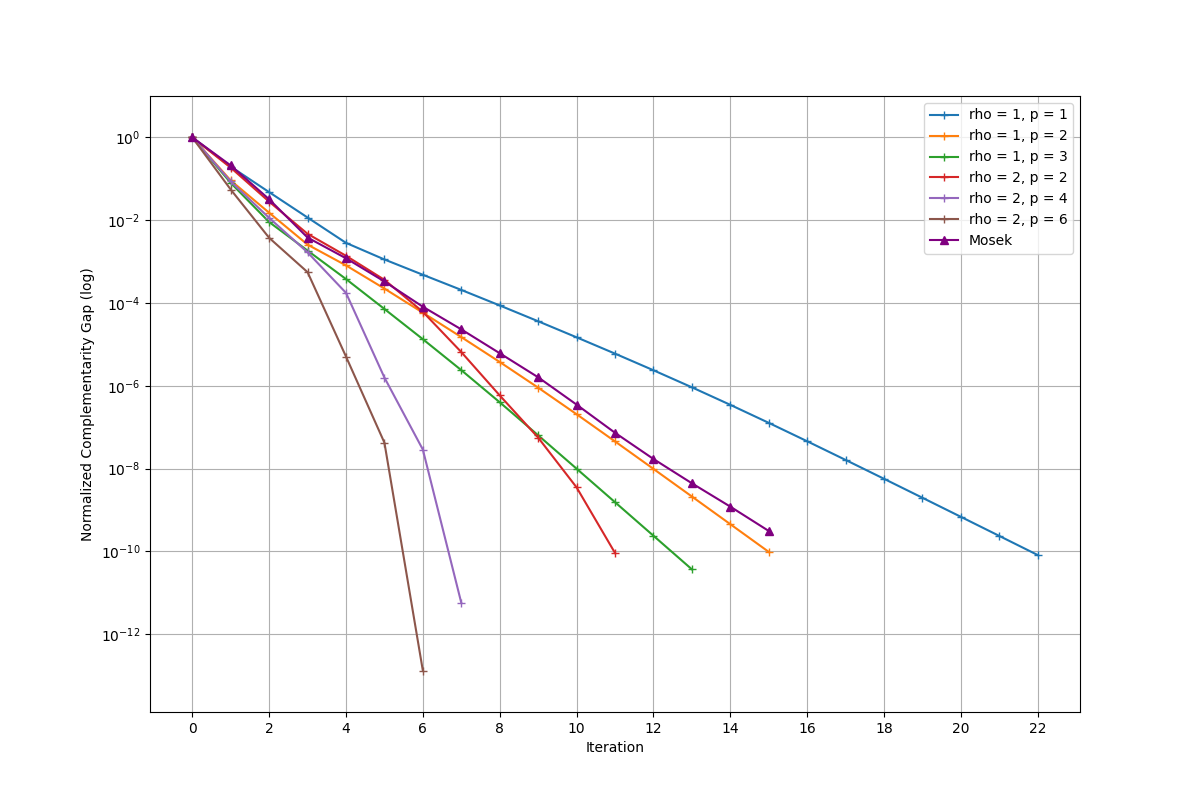}
  \par\vspace{0.3em}\small ND-4:‌ $\rho =$ 1 vs 2
\end{minipage}\hfill
\begin{minipage}[t]{0.48\textwidth}
  \centering
  \includegraphics[width=\linewidth]{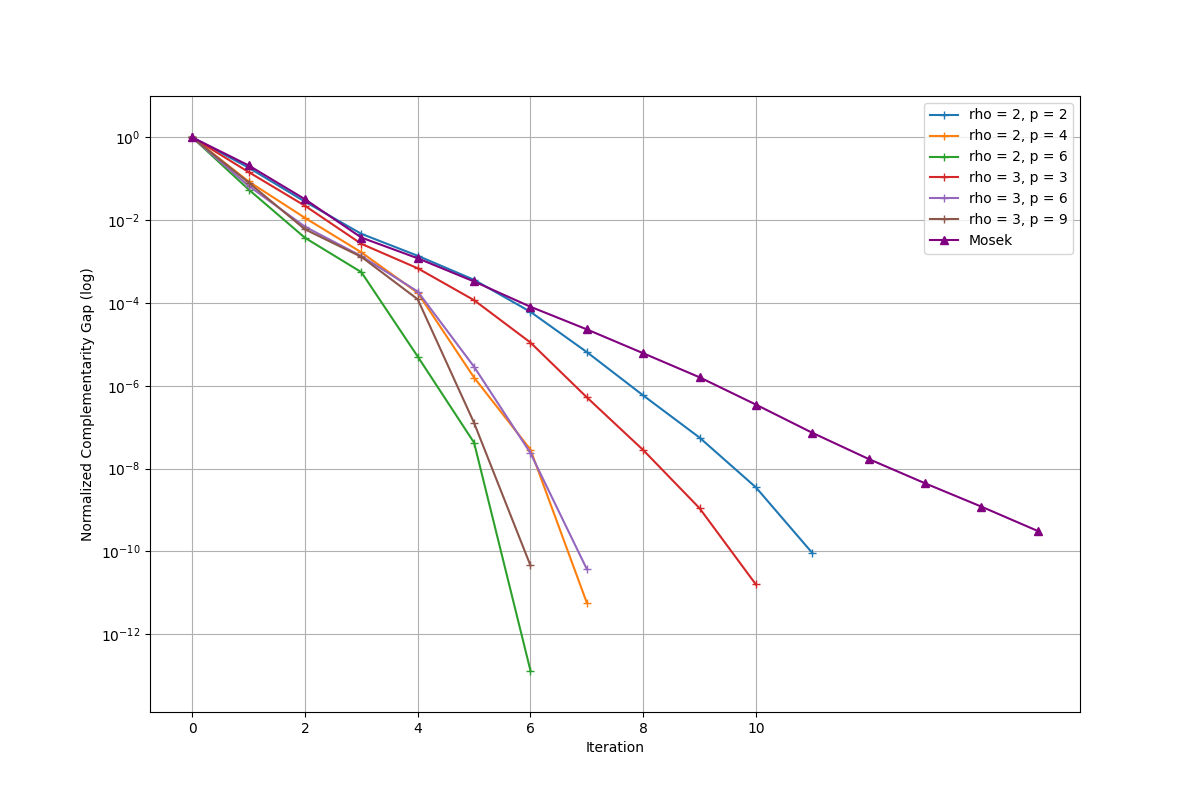}
  \par\vspace{0.3em}\small ND-4:‌ $\rho =$ 2 vs 3
\end{minipage}

\vspace{0.8em}

\begin{minipage}[t]{0.48\textwidth}
  \centering
  \includegraphics[width=\linewidth]{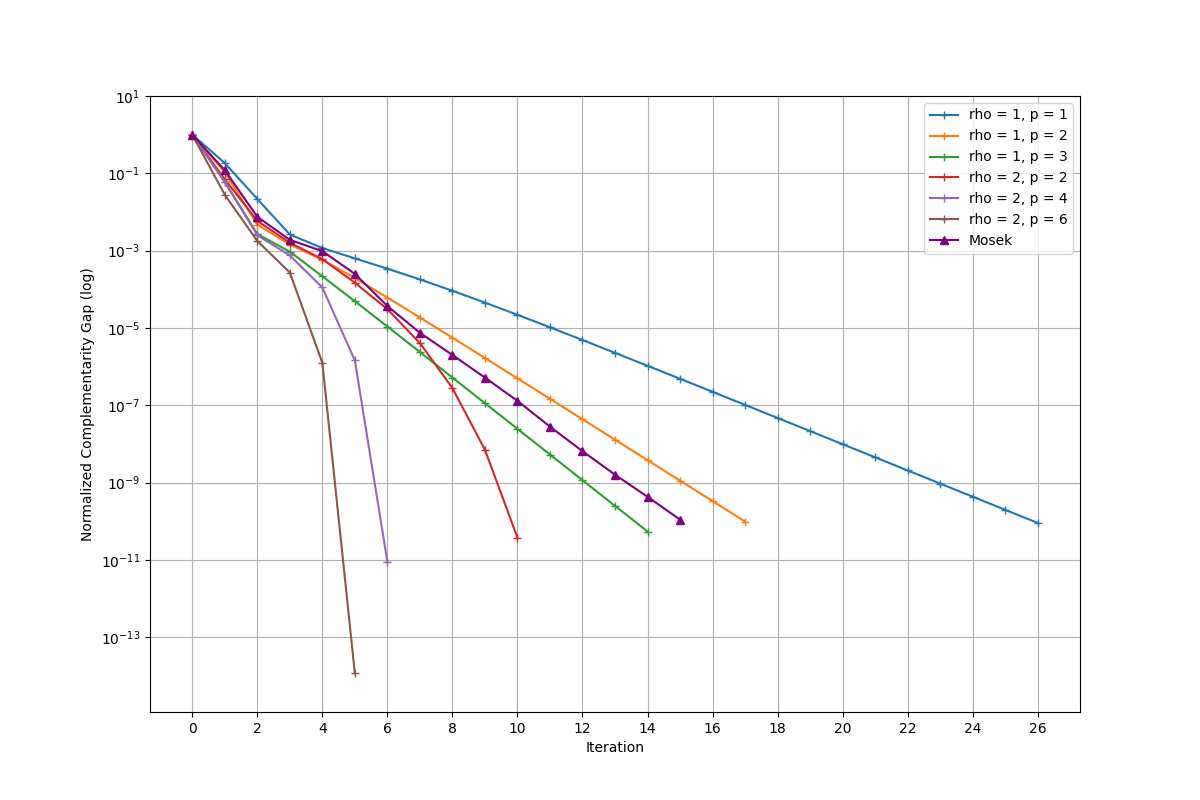}
  \par\vspace{0.3em}\small ND-5:‌ $\rho =$ 1 vs 2
\end{minipage}\hfill
\begin{minipage}[t]{0.48\textwidth}
  \centering
  \includegraphics[width=\linewidth]{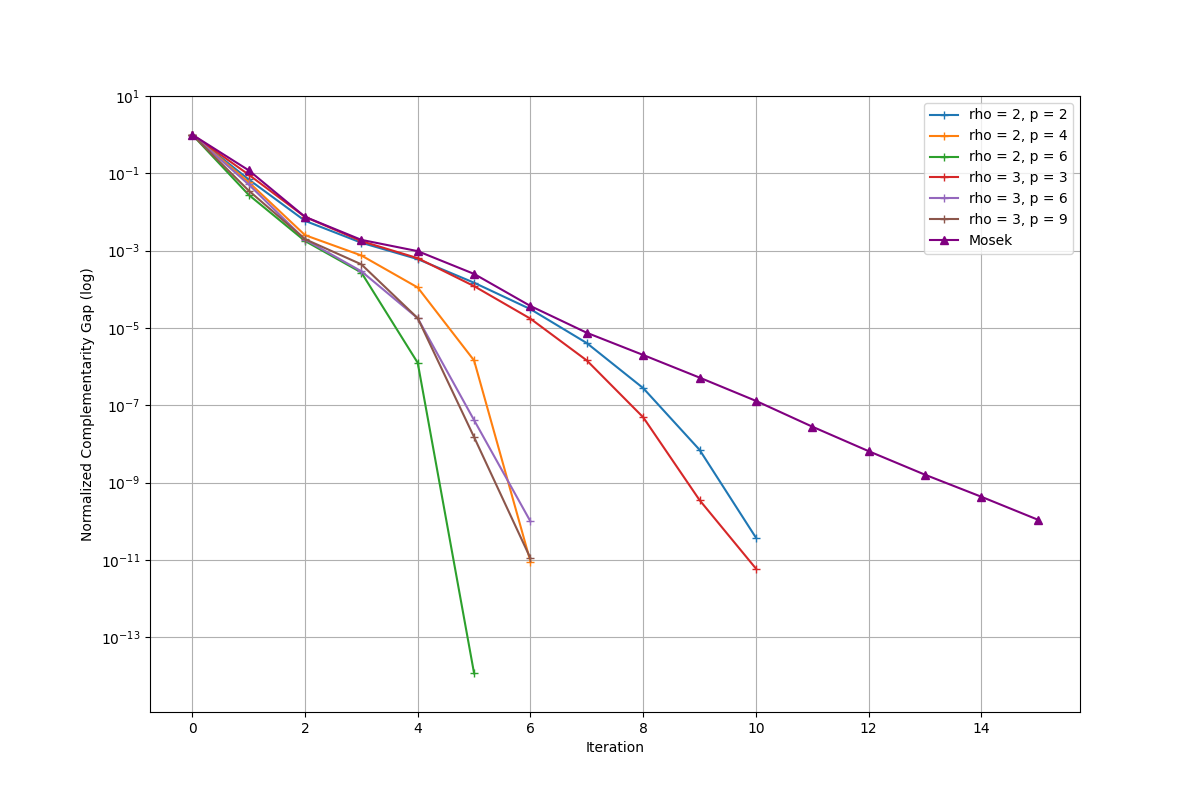}
  \par\vspace{0.3em}\small ND-5:‌ $\rho =$ 2 vs 3
\end{minipage}

\vspace{0.8em}

\begin{minipage}[t]{0.48\textwidth}
  \centering
  \includegraphics[width=\linewidth]{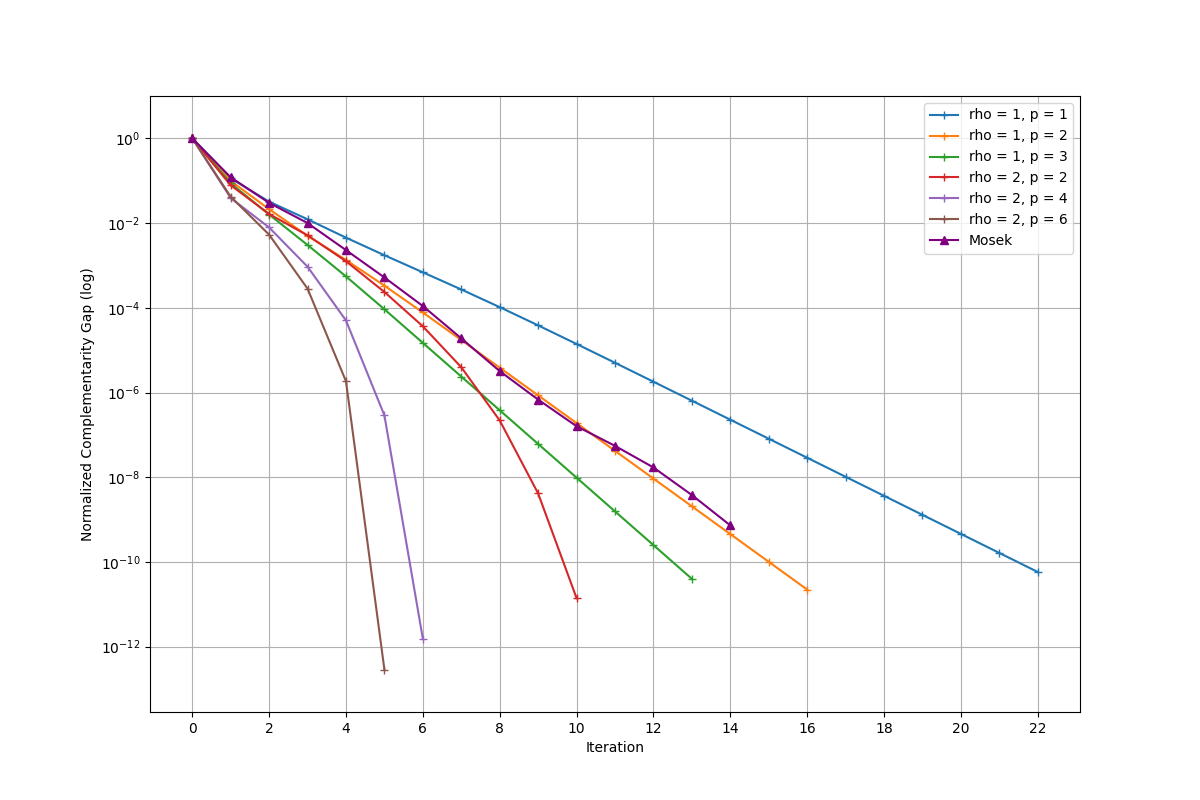}
  \par\vspace{0.3em}\small ND-6:‌ $\rho =$ 1 vs 2
\end{minipage}\hfill
\begin{minipage}[t]{0.48\textwidth}
  \centering
  \includegraphics[width=\linewidth]{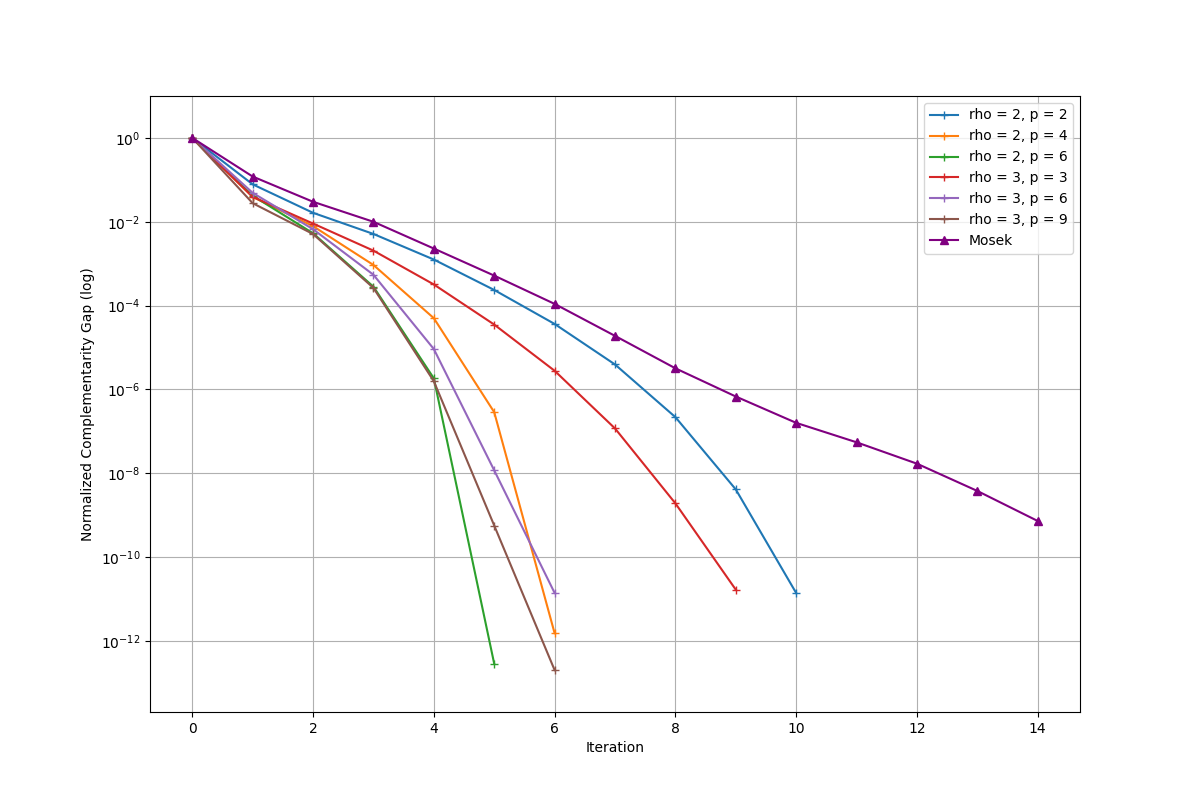}
  \par\vspace{0.3em}\small ND-6:‌ $\rho =$ 2 vs 3
\end{minipage}

\caption{Convergence of random instances with primal-dual non-degeneracy (ND4-ND6)}
\label{fig:NDConv2}
\end{figure}

\vspace{5px}
\noindent
The convergence histories of our \IPM~on six non-degenerate instances are shown in Figures~\ref{fig:NDConv1} and~\ref{fig:NDConv2}. For each instance, we first compare the performance of the method without reparametrization against the reparametrized variant with $\rho = 2$ for several preselected values of $p$. Consistent with the previous examples, increasing the order of the Taylor approximation in the absence of reparametrization does not improve the asymptotic convergence rate; it only reduces the total number of iterations required. In contrast, with $\rho=2$ we observe behavior consistent with super-linear convergence. As the iterates approach optimality, the rate of decrease of the complementarity gap accelerates. We also compared the performance of reparametrizations with $\rho = 2$ and $\rho = 3$. While the case $\rho = 2$ exhibits the expected behavior that higher-order approximations lead to steeper convergence, the same trend is not observed for $\rho = 3$. In particular, increasing the approximation order does not yield comparable improvements in the convergence speed.

\vspace{5px}
\noindent
It is worth emphasizing that our \IPM~preserves its super-linear convergence behavior up to the target complementarity gap of $\langle X, S \rangle / n = 10^{-10}$. Moreover, in several instances, our method attains a smaller final normalized complementarity gap than MOSEK under comparable tolerance settings.

\begin{table}[htbp]
\centering
\caption{Random instances with primal and dual non-degeneracy}
\renewcommand{\arraystretch}{1.25}
\begin{tabular}{cccccccc}
\hline
Instance & $m$ & $n$ & $|\mathcal{B}|$ & $|\mathcal{N}|$ & $|\mathcal{T}|$ & $\|X^*\|$ & $\|S^*\|$  \\ \hline
ND-1        &  3  &  4  &  1  &  2  &  1  & 6.45E-01 & 3.08E+00 \\ \hline
ND-2        &  4  &  6  &  1  &  4  &  1  & 6.93E-01 & 3.33E+00 \\ \hline
ND-3        &  5  &  6  &  1  &  4  &  1  & 2.64E-01 & 3.51E+00 \\ \hline
ND-4        &  6  &  7  &  1  &  4  &  2  & 8.02E-01 & 7.27E+00 \\ \hline
ND-5        &  6  &  8  &  1  &  5  &  2  & 4.05E-01 & 6.75E+00 \\ \hline
ND-6        &  7  & 10  &  1  &  7  &  2  & 2.92E-01 & 6.00E+00 \\ \hline
\end{tabular}
\label{tab:ND}
\end{table}

\subsubsection{Problems with positive duality gaps}\label{ssec:PDG}
\begin{figure}[htbp]
\centering
\begin{minipage}[t]{0.48\textwidth}
  \centering
  \includegraphics[width=\linewidth]{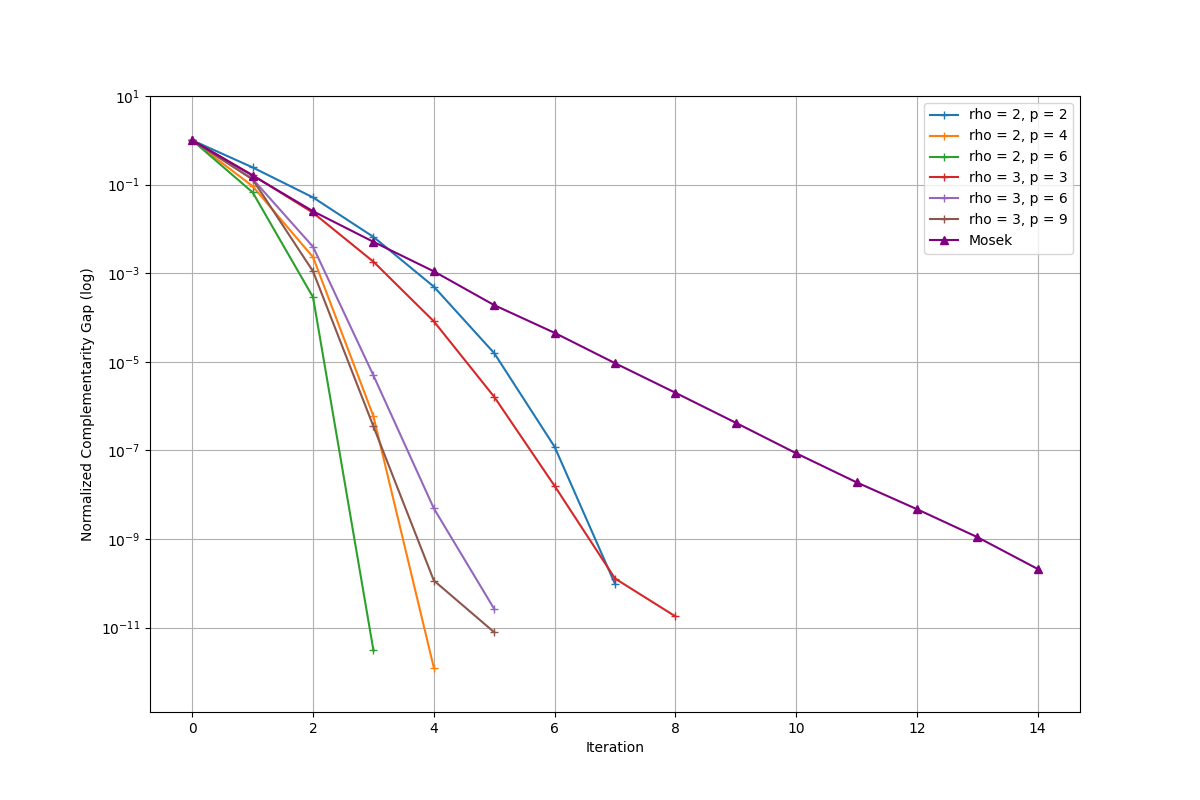}
  \par\vspace{0.3em}\small PDG: Instance 1
\end{minipage}\hfill
\begin{minipage}[t]{0.48\textwidth}
  \centering
  \includegraphics[width=\linewidth]{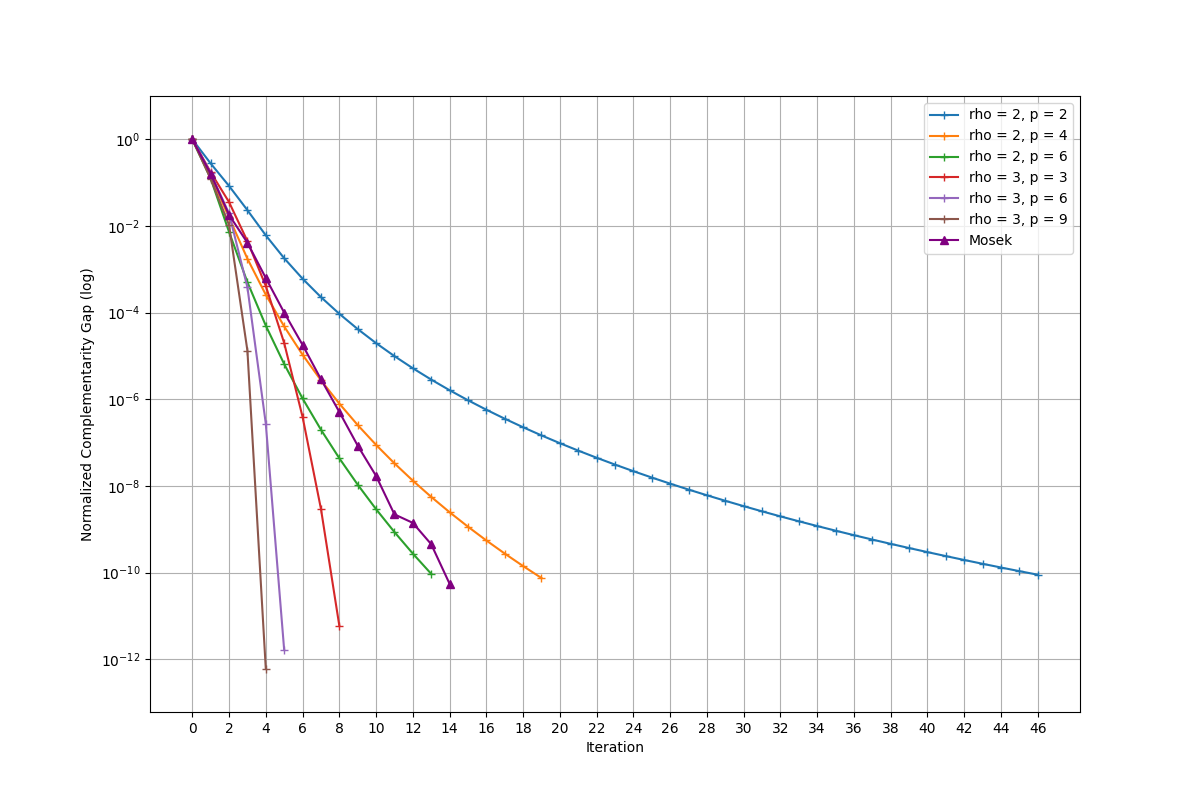}
  \par\vspace{0.3em}\small PDG: Instance 2
\end{minipage}

\vspace{0.8em}

\begin{minipage}[t]{0.48\textwidth}
  \centering
  \includegraphics[width=\linewidth]{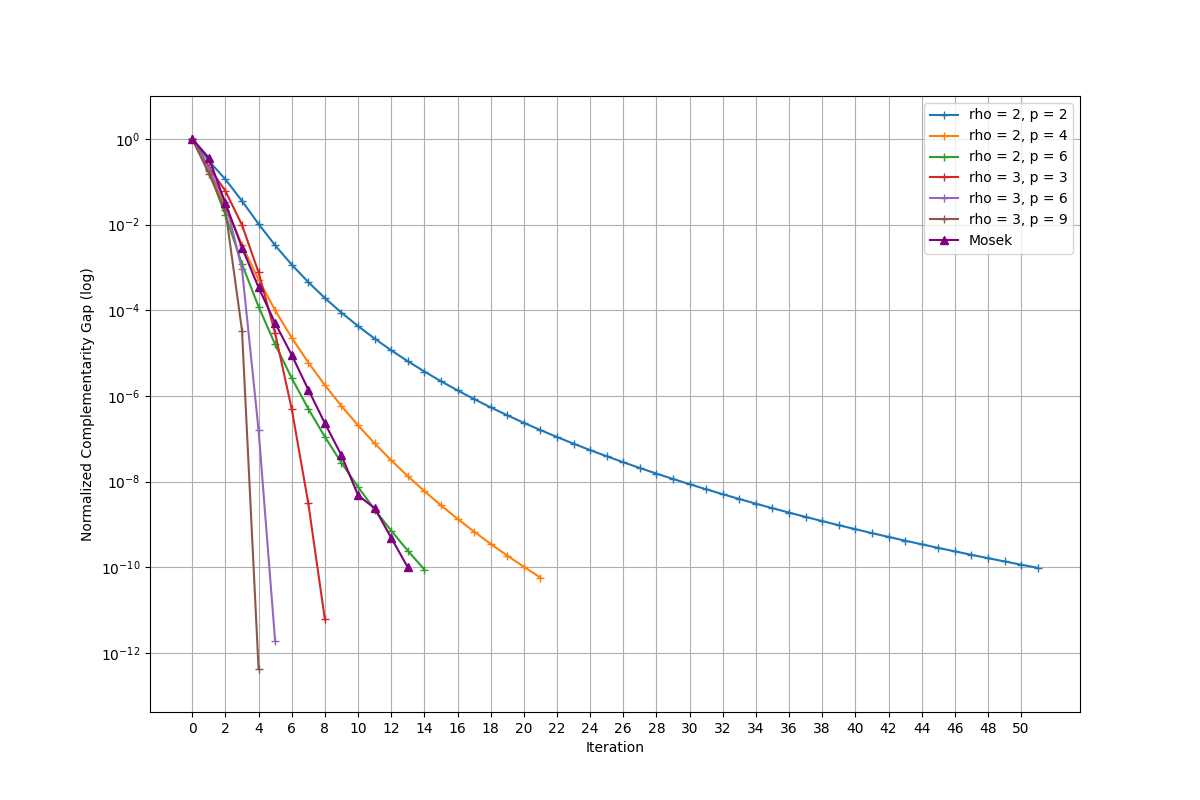}
  \par\vspace{0.3em}\small PDG: Instance 3
\end{minipage}\hfill
\begin{minipage}[t]{0.48\textwidth}
  \centering
  \includegraphics[width=\linewidth]{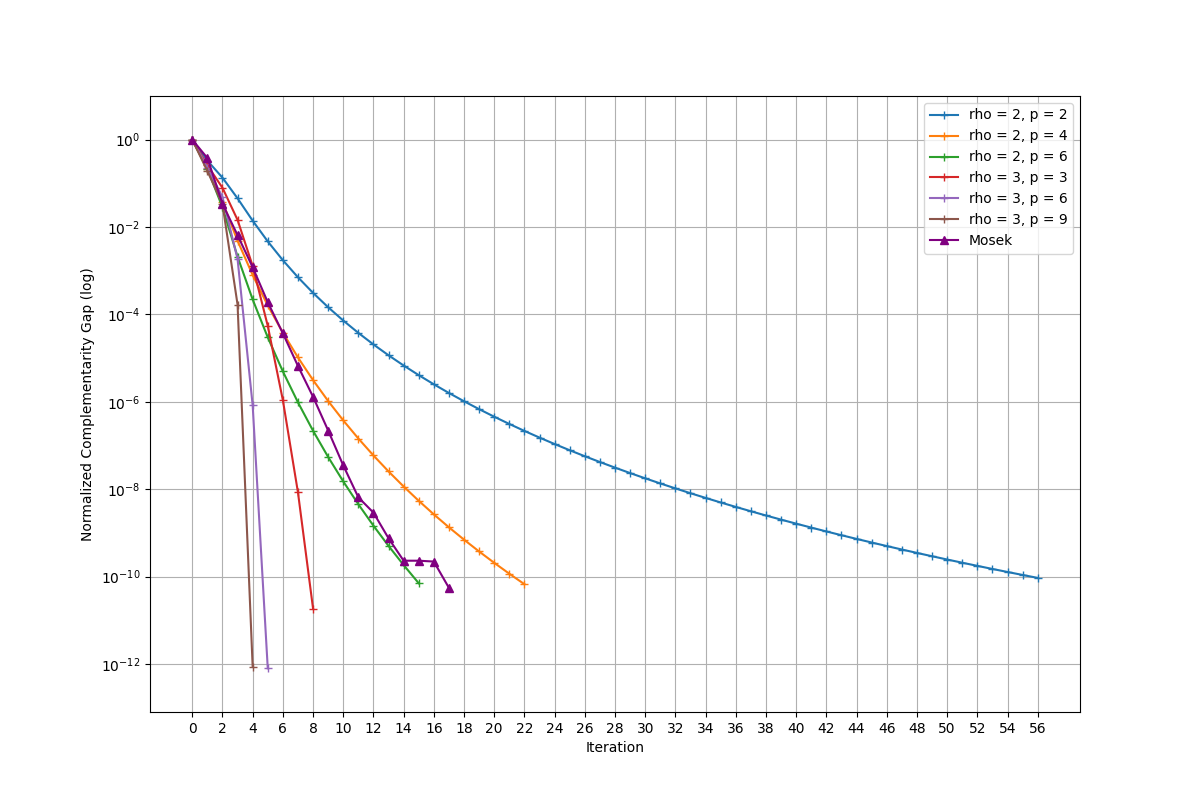}
  \par\vspace{0.3em}\small PDG: Instance 4
\end{minipage}

\vspace{0.8em}

\begin{minipage}[t]{0.48\textwidth}
  \centering
  \includegraphics[width=\linewidth]{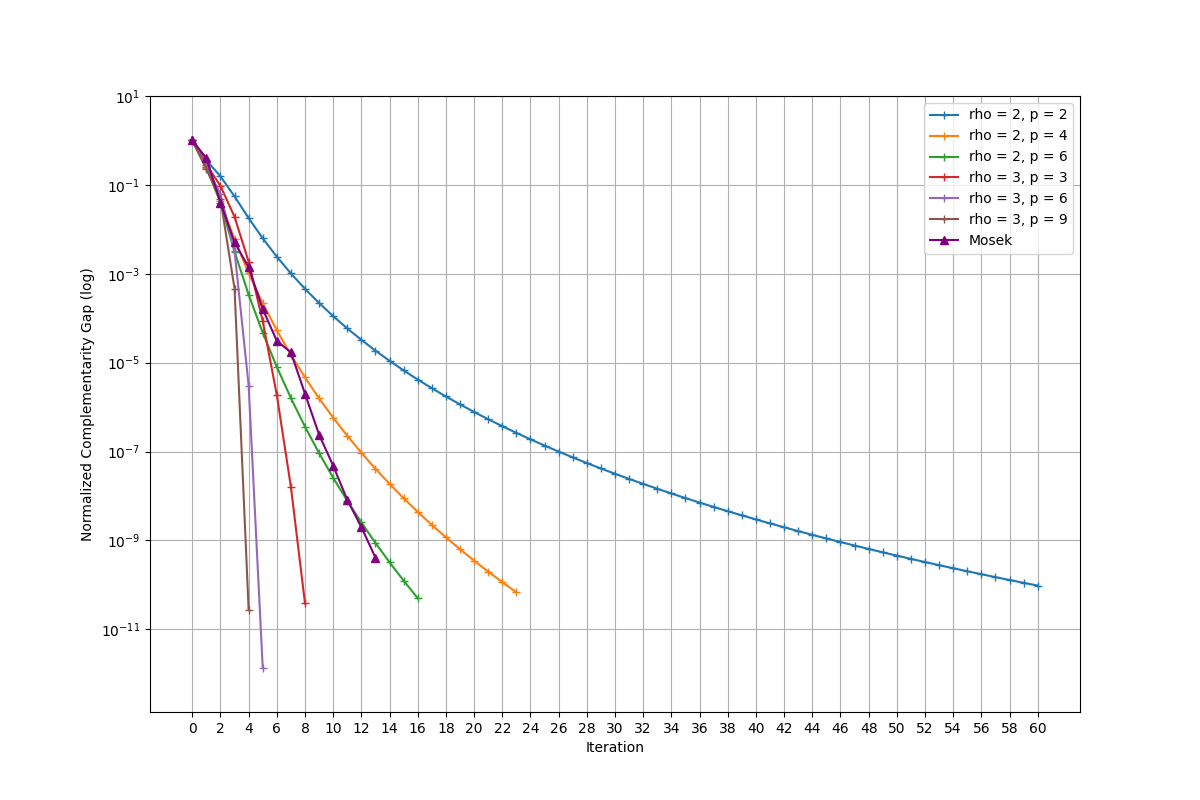}
  \par\vspace{0.3em}\small PDG: Instance 5
\end{minipage}\hfill
\begin{minipage}[t]{0.48\textwidth}
  \centering
  \includegraphics[width=\linewidth]{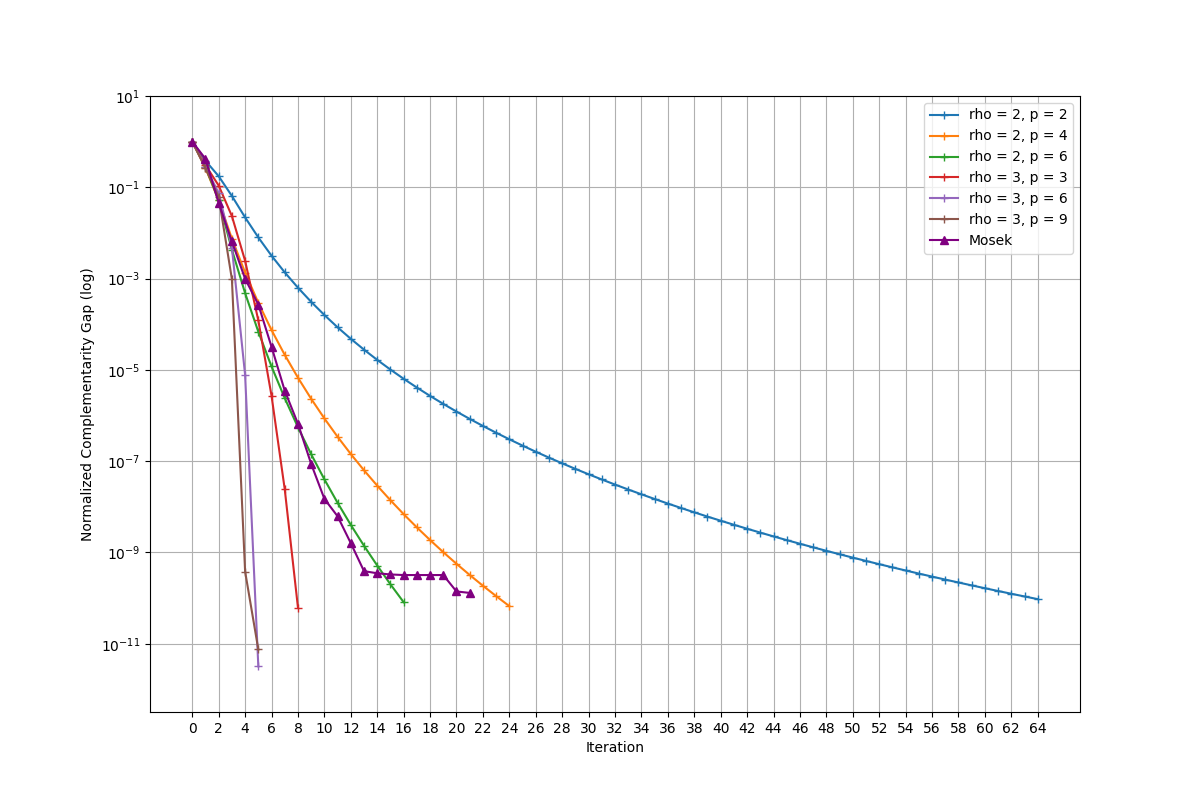}
  \par\vspace{0.3em}\small PDG: Instance 6
\end{minipage}

\caption{Convergence on PDG Problems}
\label{fig:PDGConv}
\end{figure}

In this section, we consider problems with a positive duality gap (PDG). It is known that, when formulated in homogeneous self-dual (HSD) embedding form, such problems fail to satisfy the strict complementarity condition~\cite{TW12}.  We consider a collection of PDG instances from~\cite{P18}, which are described by the author as clean examples with a finite duality gap.  Each instance is first transformed into its HSD embedding, and the resulting embedded problem is then solved using our higher-order \IPM. Table~\ref{tab:PDGdim} summarizes the dimensions of the original PDG problems, $(m,n)$, together with the dimensions of their corresponding embedded formulations, $(m_{\mathrm{HSD}}, n_{\mathrm{HSD}})$. The convergence histories of our \IPM~on several embedded PDG instances are presented in Figure~\ref{fig:PDGConv}.

\begin{table}[H]
\caption{Problems with positive duality gap}
\renewcommand{\arraystretch}{1.25}
\begin{tabular}{cccccc}
\hline
Instance & $m$ & $n$ & $m_{HSD}$ & $n_{HSD}$ \\ \hline
1        &  2  &  3  &    12     &    18     \\ \hline
2        &  3  &  4  &    18     &    24     \\ \hline
3        &  4  &  5  &    25     &    30     \\ \hline
4        &  5  &  6  &    33     &    36     \\ \hline
5        &  6  &  7  &    42     &    42     \\ \hline
6        &  7  &  8  &    52     &    48     \\ \hline
\end{tabular}
\label{tab:PDGdim}
\end{table}

\vspace{5px}
\noindent
Among the tested reparametrization exponents, $\rho=3$ is the smallest value for which we observe sustained super-linear behavior on all but Instance~1. Under the choice $\rho=3$, our method requires fewer outer iterations than MOSEK and exhibits a steeper asymptotic reduction in the normalized complementarity gap on these instances. We further observe that our algorithm converges to a primal-dual optimal solution that fails the strict complementarity condition and differs from the solution returned by MOSEK. This discrepancy is not unexpected, as the embedded problems generally do not possess a unique optimal solution. Rather, they admit a set of optimal solutions, all of which fail to satisfy strict complementarity.

\subsubsection{Experiments on instances satisfying strict complementarity}\label{ssec:SSC}
We conclude by testing our \IPM~on instances for which strict complementarity is known to hold. Since the strict complementarity condition is satisfied, no reparametrization is required, and thus $\rho = 1$ is sufficient. Nevertheless, for comparison purposes, we also include the cases $\rho = 2,3$.

\vspace{5px}
\noindent
\paragraph{\textbf{Random problems with dual non-degeneracy and strict complementarity.}}
In this section, we consider randomly generated problems that satisfy the dual non-degeneracy condition of~\cite{AHO97} (which implies the uniqueness of the primal optimal solution~\cite[Theorem~10]{AHO97}) and the strict complementarity condition, but they violate the primal non-degeneracy condition (see Table~\ref{tab:D}). Similar to the problem set considered in Section~\ref{ssec:GenP}, these instances are generated using the procedure described in~\cite[Algorithm~11]{MFAT24}. 

\begin{figure}[htbp]
\centering

\begin{minipage}[t]{0.48\textwidth}
  \centering
  \includegraphics[width=\linewidth]{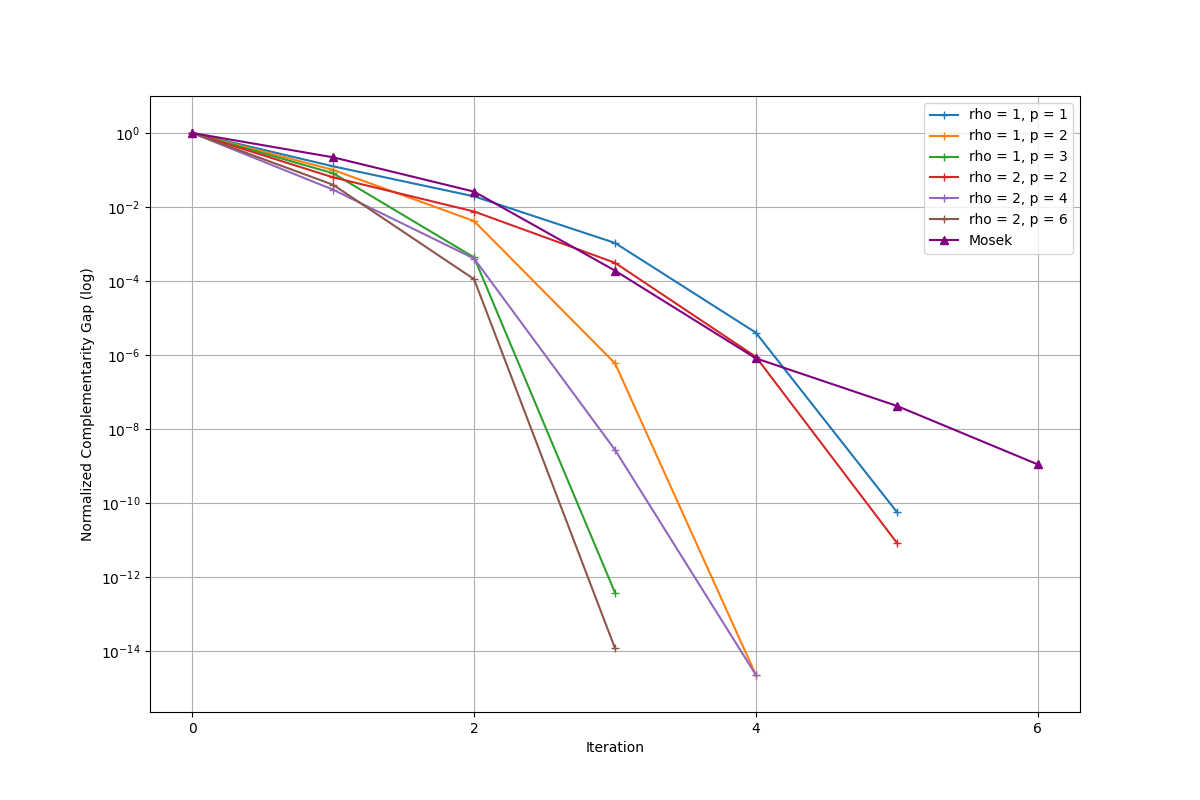}
  \par\vspace{0.3em}\small D-1: $\rho = $ 1 vs 2
\end{minipage}\hfill
\begin{minipage}[t]{0.48\textwidth}
  \centering
  \includegraphics[width=\linewidth]{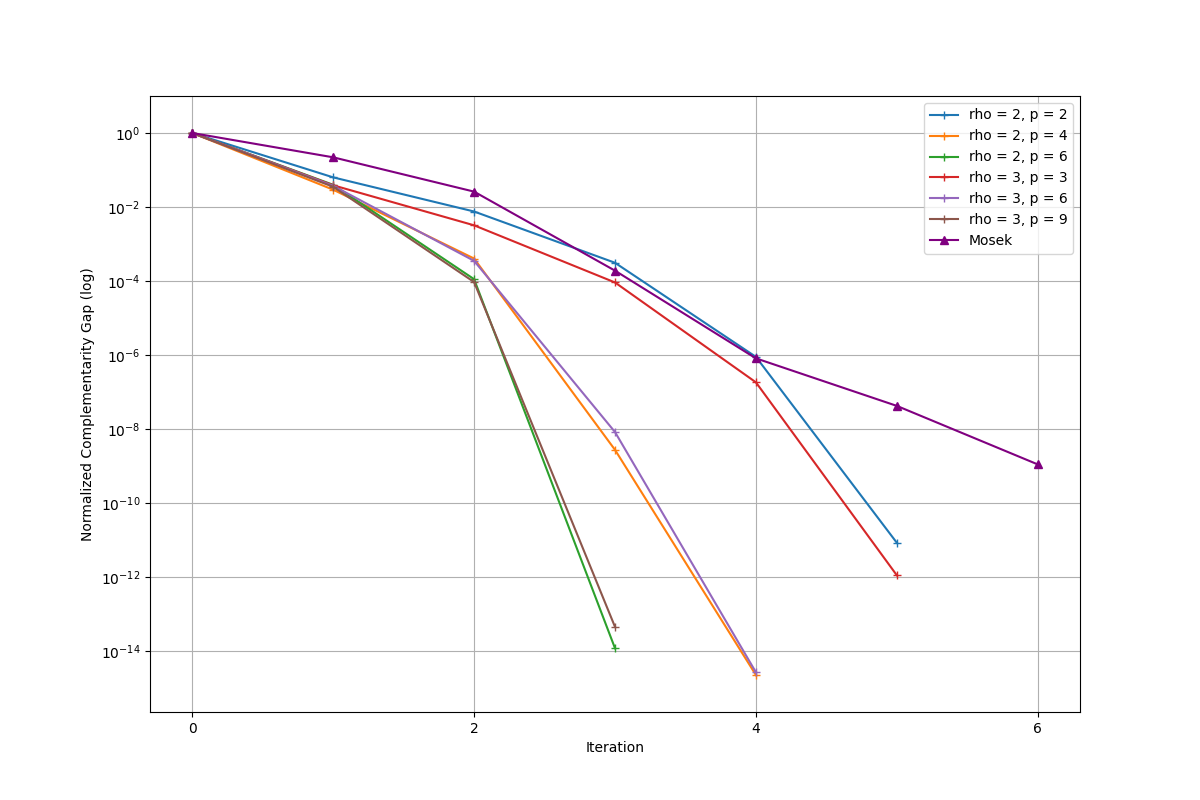}
  \par\vspace{0.3em}\small D-1: $\rho = $ 2 vs 3
\end{minipage}

\begin{minipage}[t]{0.48\textwidth}
  \centering
  \includegraphics[width=\linewidth]{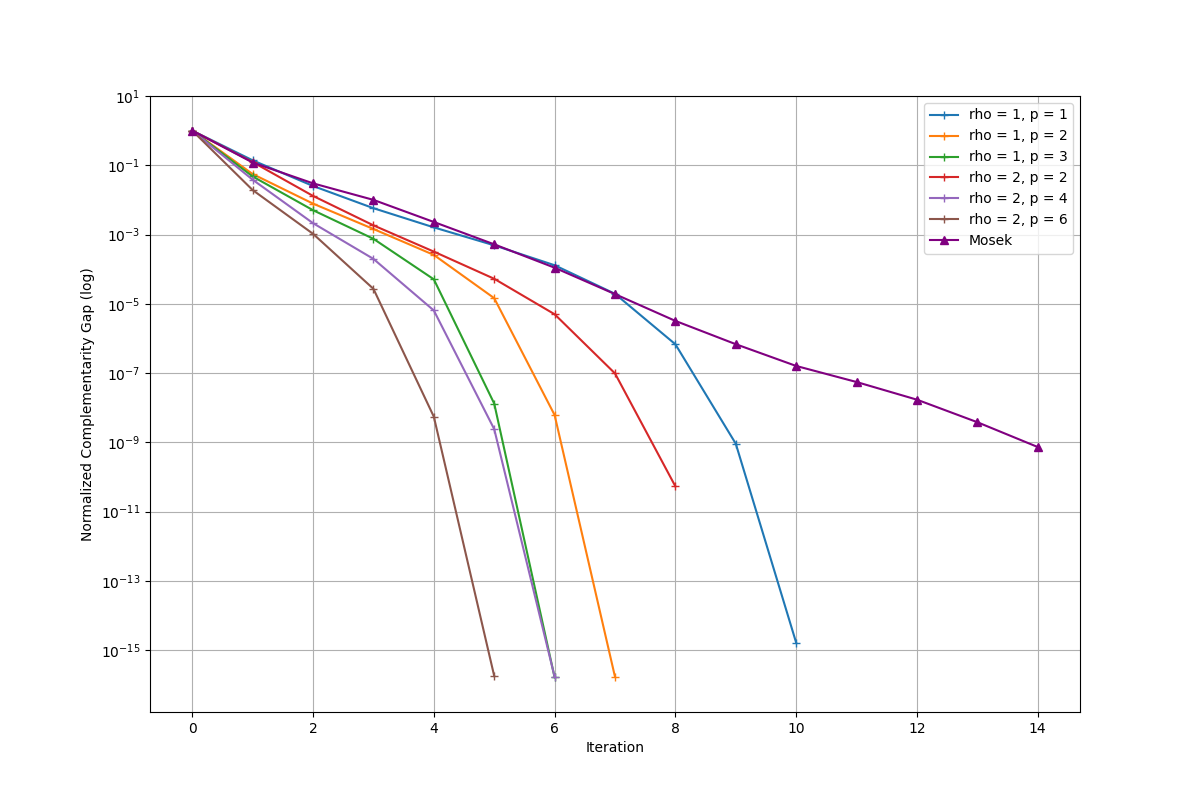}
  \par\vspace{0.3em}\small D-2: $\rho = $ 1 vs 2
\end{minipage}\hfill
\begin{minipage}[t]{0.48\textwidth}
  \centering
  \includegraphics[width=\linewidth]{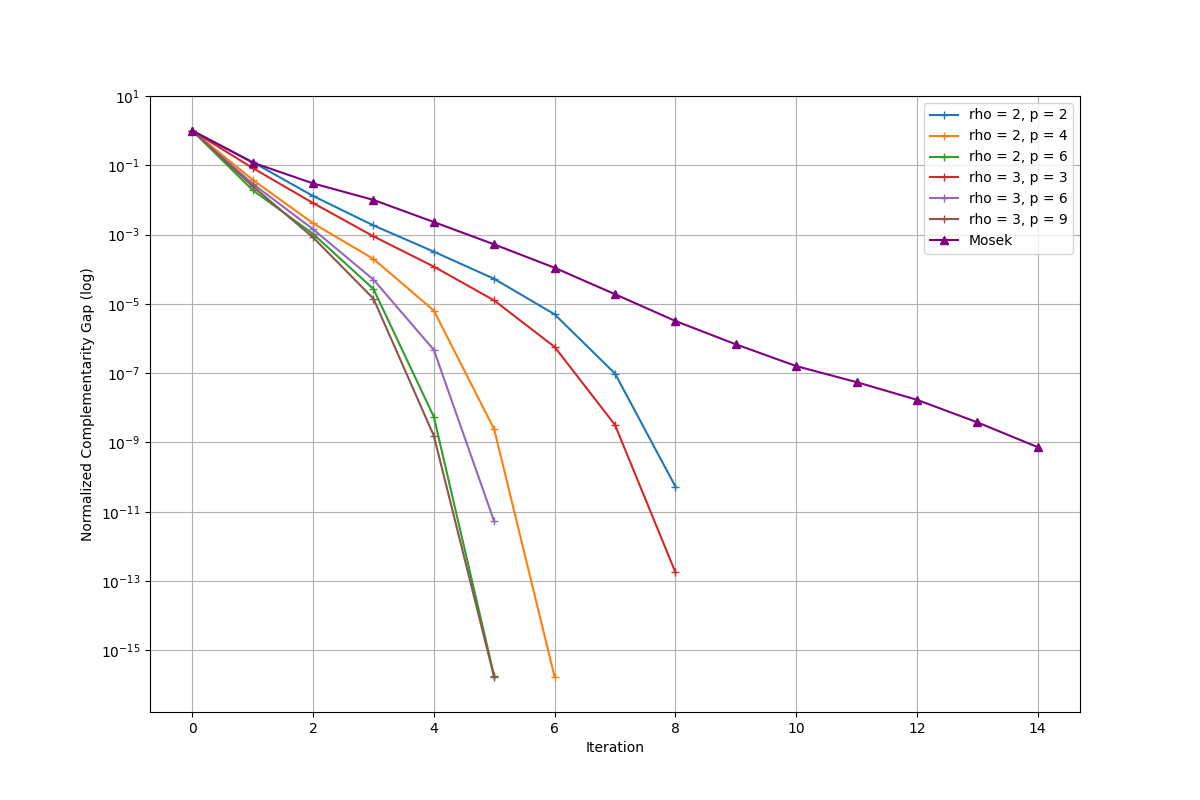}
  \par\vspace{0.3em}\small D-2: $\rho = $ 2 vs 3
\end{minipage}

\vspace{0.8em}

\begin{minipage}[t]{0.48\textwidth}
  \centering
  \includegraphics[width=\linewidth]{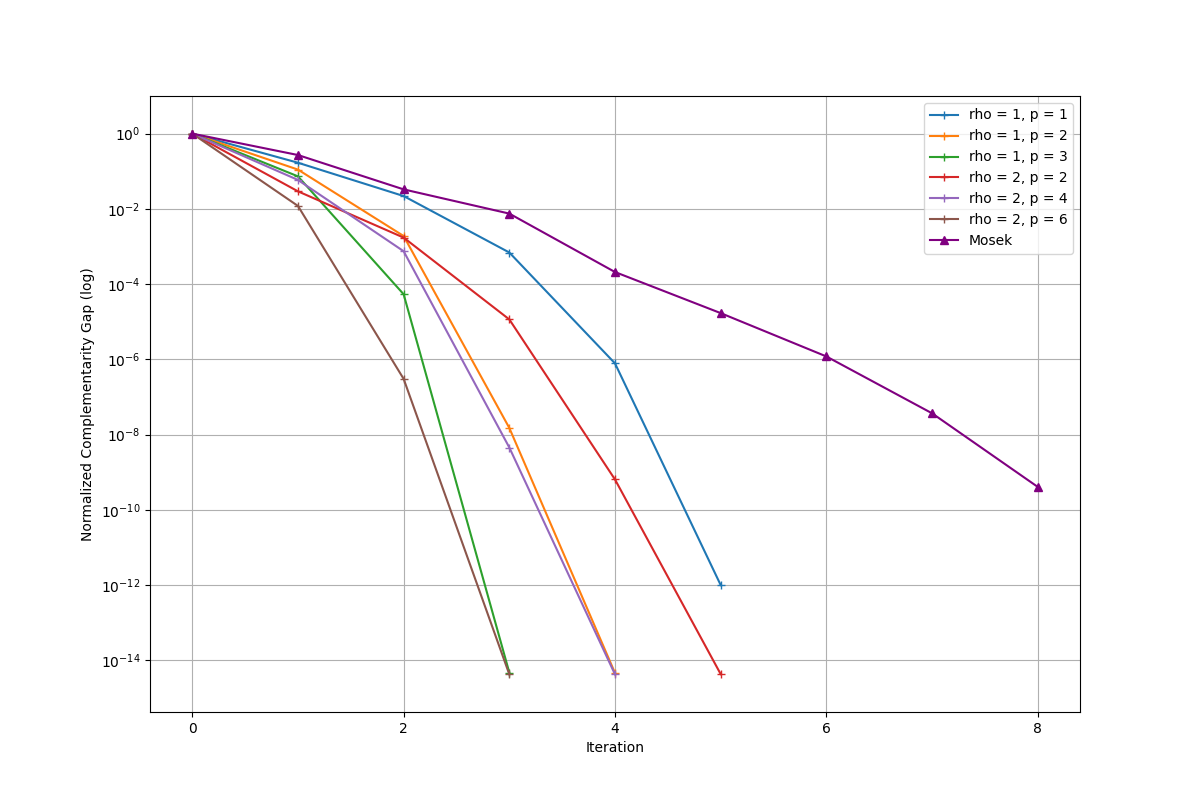}
  \par\vspace{0.3em}\small D-3: $\rho = $ 1 vs 2
\end{minipage}\hfill
\begin{minipage}[t]{0.48\textwidth}
  \centering
  \includegraphics[width=\linewidth]{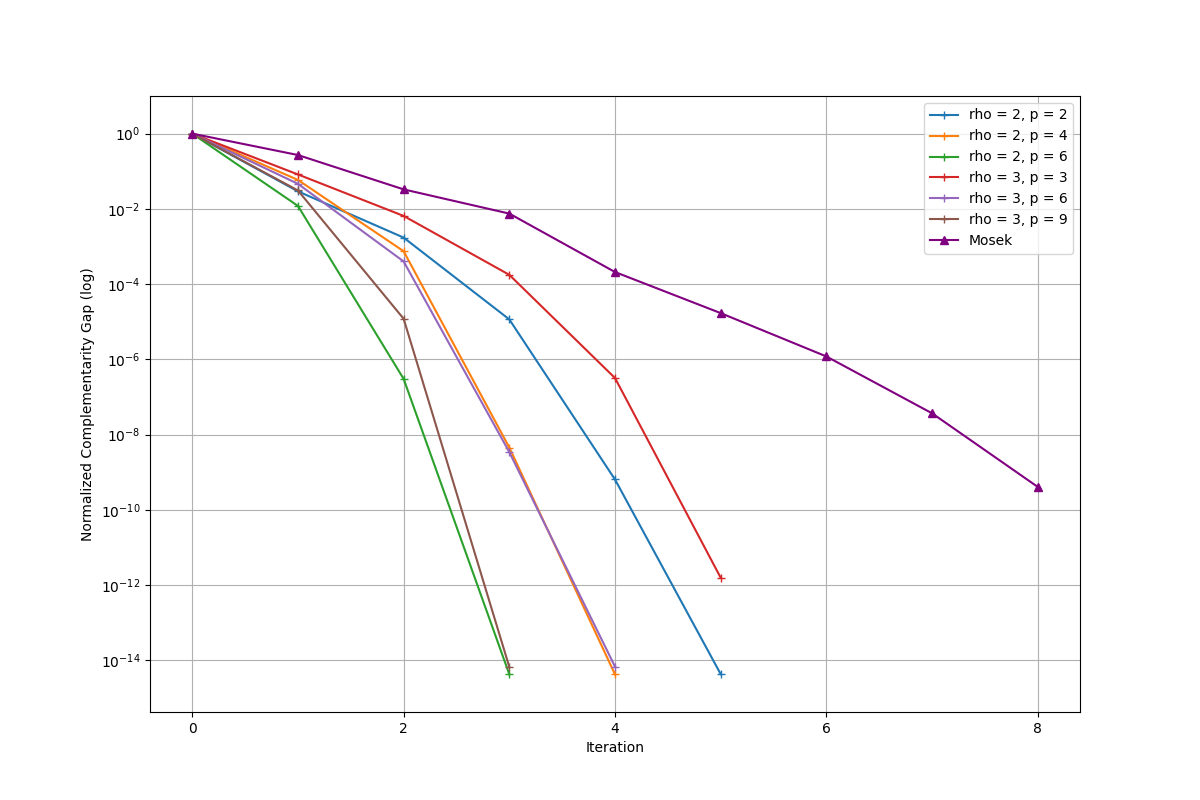}
  \par\vspace{0.3em}\small D-3: $\rho = $ 2 vs 3
\end{minipage}

\caption{Convergence of random instances with dual non-degeneracy and strict complementarity}
\label{fig:DConv1}
\end{figure}

\vspace{5px}
\noindent
The convergence histories of our \IPM~on three such instances are presented in Figure~\ref{fig:DConv1}. The super-linear convergence of predictor-corrector \IPM s under strict complementarity is well established (see Section~\ref{local_convergence_IPMs}), and our numerical results are consistent with this theory. In particular, our \IPM~converges substantially faster than in instances whose unique optimal solution does not satisfy strict complementarity. The convergence profiles exhibit a super-linear rate across all instances and for a wider range of choices of $\rho$. Furthermore, although the target complementarity gap was set to $10^{-10}$, the algorithm attained solutions with significantly smaller complementarity gaps in all tested instances.

\begin{table}[htbp]
\centering
\caption{Random instances with dual non-degeneracy and strict complementarity}
\renewcommand{\arraystretch}{1.25}
\begin{tabular}{ccccccc}
\hline
Instance & $m$ & $n$ & $|\mathcal{B}|$ & $|\mathcal{N}|$ & $\|X^*\|$ & $\|S^*\|$  \\ \hline
D-1        &  6  &  5  &  1  &  4  & 2.18E-01 & 1.13E+00 \\ \hline
D-2        &  8  &  8  &  1  &  7  & 6.48E-01 & 3.89E+00 \\ \hline
D-3        &  10 &  9  &  1  &  8  & 4.14E-01 & 1.11E+00 \\ \hline
\end{tabular}
\label{tab:D}
\end{table}

\vspace{5px}
\noindent
\paragraph{\textbf{Random problems with strict complementarity.}}
Next, we consider two benchmark instances: one from the literature~\cite{HLMZ21}:
\begin{equation}\label{Ha_bench}
\inf \left\{ 2x_{11} - 2x_{12} \mid \
        \begin{bmatrix}
            x_{11} & x_{12} \\
            x_{12} & (1-2x_{11})/3
\end{bmatrix}
\succeq 0\right\},
\end{equation}
and another generated using the \SDO\ generator of~\cite[Algorithm~11]{MFAT24}. We generate a small instance where $n=5, \ m=5, \ |\mathcal{B}| = 2, \ |\mathcal{N}| = 3, \ |\mathcal{T}| = 0$ (recall that $|\mathcal{T}| = 0$ implies the strict complementarity), which also satisfies the primal-dual non-degeneracy conditions of~\cite{AHO97}.

\begin{figure}
    \centering
    \includegraphics[width=.75\linewidth]{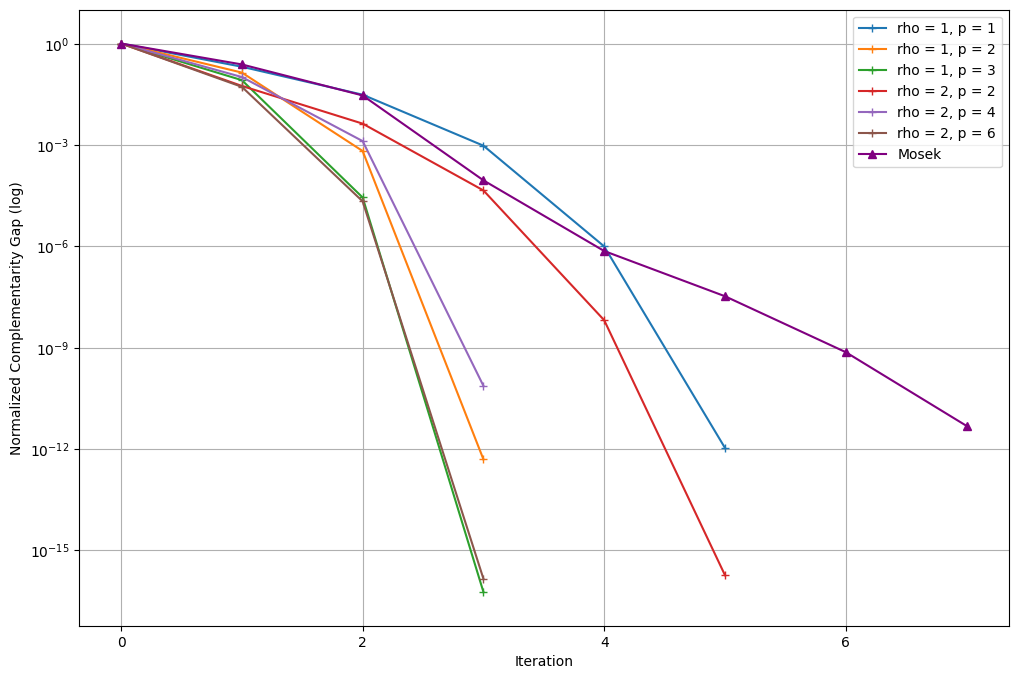}
    \caption{Convergence history on the benchmark problem~\eqref{Ha_bench}}
    \label{fig:HNS}
\end{figure}

\begin{figure}
    \centering
    \includegraphics[width=.75\linewidth]{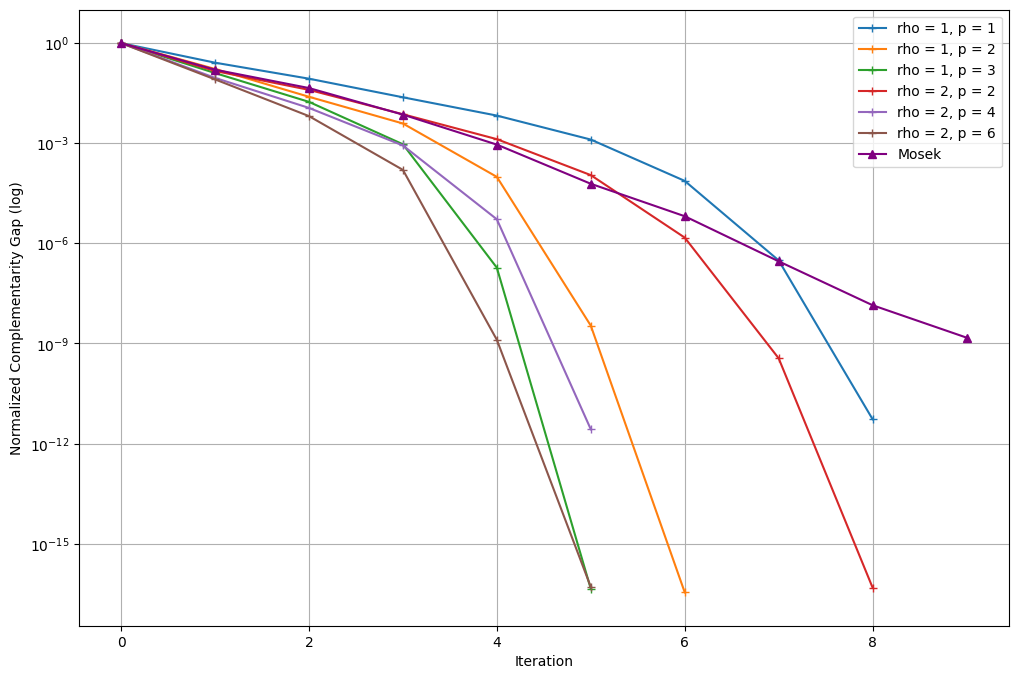}
    \caption{Convergence history on the random problem with strict complementarity}
    \label{fig:Random_strict}
\end{figure}

\vspace{5px}
\noindent 
Figures~\ref{fig:HNS} and~\ref{fig:Random_strict} present the numerical results for the benchmark problem~\eqref{Ha_bench} and the random problem. Across all experiments, we observe a super-linear rate of convergence for every choice of $p$. Moreover, as expected, increasing the order of the Taylor approximation leads to greater reductions in the complementarity gap.

\noindent

\subsection{Convergence of the sequence of complementarity gaps}\label{sec:convseq}
In this section, we examine the convergence behavior of the sequence of complementarity gaps in greater detail. Under the assumptions of Theorem~\ref{Algorithm_Convergence}, the
complementarity gaps satisfy
\begin{align*}
\mu_k=O\!\left(\mu_{k-1}^{(p+1)/\rho}\right).
\end{align*}
Hence, for any $\omega<(p+1)/\rho$, the ratio
\begin{align*}
\frac{\mu_k}{\mu_{k-1}^{\omega}}
\end{align*}
tends to zero in the local regime. We take $\omega=1.05$ throughout this section as a numerical diagnostic of convergence behavior.

\vspace{5px}
\noindent
As before, we begin with the 3-elliptope problem. Previously, we demonstrated the promising performance obtained with $\rho = 2$. For several choices of $p$, Table~\ref{tab:convhist-3ell-rho2} reports the convergence history of our \IPM.
\begin{table}[H]
\caption{Convergence history for the 3-elliptope problem (Figure~\ref{fig:3ell-rho1vs2Multip})}
\begin{tabular}{ccccccc}
\hline
~ & \multicolumn{2}{c}{$(\rho,p)=(2,2)$} & \multicolumn{2}{c}{$(\rho,p)=(2,4)$} & \multicolumn{2}{c}{$(\rho,p)=(2,8)$} \\ \cline{2-7} 
Iter. ($k$) & $\langle X_k,S_k\rangle/n$ & Ratio       & $\langle X_k,S_k\rangle/n$ & Ratio       & $\langle X_k,S_k\rangle/n$ & Ratio       \\ \hline
0 & 1.00E+00 & N/A & 1.00E+00 & N/A & 1.00E+00 & N/A \\
1 & 1.49E-01 & 1.49E-01 & 2.58E-02 & 2.58E-02 & 1.00E-03 & 1.00E-03 \\
2 & 1.14E-02 & 8.46E-02 & 1.08E-07 & 5.01E-06 & 7.51E-15 & 1.06E-11 \\
3 & 3.09E-04 & 3.38E-02 & 6.29E-15 & 1.30E-07 & - & - \\
4 & 1.73E-06 & 8.41E-03 & - & - & - & - \\
5 & 8.48E-10 & 9.51E-04 & - & - & - & - \\
6 & 1.17E-14 & 3.93E-05 & - & - & - & - \\
\hline
\end{tabular}
\label{tab:convhist-3ell-rho2}
\end{table} 

For all choices of $p$, the reported values of $\frac{\mu_k}{\mu_{k-1}^{\,1.05}}$ decrease rapidly toward zero. This provides numerical evidence of a $Q$-order exceeding $1.05$ over the reported iterates. Moreover, larger values of $p$ lead to a more rapid decrease in this ratio, consistent with the higher convergence orders predicted by Theorem~\ref{Algorithm_Convergence}.

\vspace{5px}
\noindent
Earlier, we noted that using $\rho=3$ for the 3-elliptope problem may reduce the number of iterations, although it is not expected to preserve super-linear convergence near optimality.  To investigate this behavior, Table~\ref{tab:convhist-3ell-rho3} presents the convergence history obtained with $\rho = 3$. As the iterates approach optimality, the ratio eventually ceases to decrease, so the reported iterates do not provide evidence of a $Q$-order exceeding $1.05$ in the asymptotic regime.

\begin{table}[H]
\caption{Convergence history for the 3-elliptope problem (Figure~\ref{fig:3ell-rho3})}
\begin{tabular}{ccccccc}
\hline
~ & \multicolumn{2}{c}{$({\rho},p)=(3,3)$} & \multicolumn{2}{c}{$({\rho},p)=(3,6)$} & \multicolumn{2}{c}{$({\rho},p)=(3,9)$} \\ \cline{2-7}
Iter. ($k$) & $\langle X_k,S_k\rangle/n$ & Ratio & $\langle X_k,S_k\rangle/n$ & Ratio & $\langle X_k,S_k\rangle/n$ & Ratio \\ \hline
0 & 1.00E+00 & N/A & 1.00E+00 & N/A & 1.00E+00 & N/A \\
1 & 1.09E-01 & 1.09E-01 & 1.89E-02 & 1.89E-02 & 1.77E-03 & 1.77E-03 \\
2 & 6.52E-03 & 6.72E-02 & 7.47E-06 & 4.81E-04 & 1.90E-07 & 1.47E-04 \\
3 & 1.47E-04 & 2.90E-02 & 2.71E-09 & 6.54E-04 & 2.02E-11 & 2.31E-04 \\
4 & 9.80E-07 & 1.04E-02 & 9.81E-13 & 9.72E-04 & - & - \\
5 & 3.00E-09 & 6.11E-03 & - & - & - & - \\
6 & 8.89E-12 & 7.92E-03 & - & - & - & - \\
\hline
\end{tabular}
\label{tab:convhist-3ell-rho3}
\end{table}

\vspace{5px}
\noindent
Next, we consider the random instance ND-5. Table~\ref{tab:convhist-ND-5} presents the
convergence history of our \IPM~for this instance. After an initial transient phase, the ratios eventually decrease for all three choices
of $p$, with a substantially faster decrease for the higher-order predictors. This behavior is consistent with the onset of the local super-linear regime.

\begin{table}[H]
\caption{Convergence history for the random instance ND-5 (Figure~\ref{fig:NDConv2})}
\begin{tabular}{ccccccc}
\hline
~ & \multicolumn{2}{c}{$({\rho},p)=(2,2)$} & \multicolumn{2}{c}{$({\rho},p)=(2,4)$} & \multicolumn{2}{c}{$({\rho},p)=(2,6)$} \\ \cline{2-7}
Iter. ($k$) & $\langle X_k,S_k\rangle/n$ & Ratio & $\langle X_k,S_k\rangle/n$ & Ratio & $\langle X_k,S_k\rangle/n$ & Ratio \\ \hline
0 & 1.00E+00 & N/A & 1.00E+00 & N/A & 1.00E+00 & N/A \\
1 & 6.97E-02 & 6.97E-02 & 5.94E-02 & 5.94E-02 & 2.78E-02 & 2.78E-02 \\
2 & 5.91E-03 & 9.69E-02 & 2.52E-03 & 4.89E-02 & 1.78E-03 & 7.64E-02 \\
3 & 1.60E-03 & 3.50E-01 & 7.50E-04 & 4.02E-01 & 2.72E-04 & 2.10E-01 \\
4 & 6.03E-04 & 5.20E-01 & 1.11E-04 & 2.13E-01 & 1.25E-06 & 6.94E-03 \\
5 & 1.48E-04 & 3.55E-01 & 1.44E-06 & 2.04E-02 & 1.16E-14 & 1.84E-08 \\
6 & 3.06E-05 & 3.22E-01 & 9.05E-12 & 1.23E-05 & - & - \\
7 & 4.06E-06 & 2.23E-01 & - & - & - & - \\
8 & 2.78E-07 & 1.28E-01 & - & - & - & - \\
9 & 7.06E-09 & 5.40E-02 & - & - & - & - \\
10 & 3.77E-11 & 1.37E-02 & - & - & - & - \\
\hline
\end{tabular}
\label{tab:convhist-ND-5}
\end{table}
\noindent
Finally, we consider one of the PDG problems and examine the convergence history of our \IPM~applied to its homogeneous self-dual embedding. Specifically, we revisit PDG Instance 3, whose embedded problem has dimensions $m_{HSD}=25$ and $n_{HSD}=30$. As in the previous experiments, we set $\rho=3$, which was observed to perform well for this instance. Table~\ref{tab:convhist-PDG-3} shows that, for all three choices of $p$, the reported ratios
eventually decrease toward zero. This provides numerical evidence of a $Q$-order exceeding $1.05$ over the reported iterates, with the decrease becoming more pronounced as $p$ increases.

\begin{table}[H]
\caption{Convergence history for the embedding of PDG Instance 3 (Figure~\ref{fig:PDGConv})}
\begin{tabular}{ccccccc}
\hline
~ & \multicolumn{2}{c}{$({\rho},p)=(3,3)$} & \multicolumn{2}{c}{$({\rho},p)=(3,6)$} & \multicolumn{2}{c}{$({\rho},p)=(3,9)$} \\ \cline{2-7}
Iter. ($k$) & $\langle X_k,S_k\rangle/n$ & Ratio & $\langle X_k,S_k\rangle/n$ & Ratio & $\langle X_k,S_k\rangle/n$ & Ratio \\ \hline
0 & 1.00E+00 & N/A & 1.00E+00 & N/A & 1.00E+00 & N/A \\
1 & 2.19E-01 & 2.19E-01 & 1.95E-01 & 1.95E-01 & 1.54E-01 & 1.54E-01 \\
2 & 6.06E-02 & 2.99E-01 & 3.44E-02 & 1.92E-01 & 2.18E-02 & 1.56E-01 \\
3 & 9.76E-03 & 1.85E-01 & 9.07E-04 & 3.12E-02 & 3.35E-05 & 1.87E-03 \\
4 & 7.66E-04 & 9.89E-02 & 1.59E-07 & 2.48E-04 & 4.20E-13 & 2.10E-08 \\
5 & 2.88E-05 & 5.39E-02 & 1.86E-12 & 2.57E-05 & - & - \\
6 & 4.94E-07 & 2.89E-02 & - & - & - & - \\
7 & 3.10E-09 & 1.30E-02 & - & - & - & - \\
8 & 6.09E-12 & 5.23E-03 & - & - & - & - \\
\hline
\end{tabular}
\label{tab:convhist-PDG-3}
\end{table}
\noindent
We also monitored the number of corrector steps and the norms of the
higher-order derivatives throughout these experiments. In most runs, a
single full centering step was sufficient to return the iterate to the
$10^{-1}$-neighborhood. The main systematic exception was the
3-elliptope with $\rho=3$. Since $W=I_n$ belongs to $\mathcal W_\beta$
and the ordinary central path is not analytic after the
reparametrization $\mu\mapsto\mu^3$, this value of $\rho$ is not
admissible on $\mathcal W_\beta$. In this case, as the iterations
progressed, two full centering steps were required. For ND-5, two
centering steps were needed only at isolated iterations, while one full
step sufficed otherwise.

\vspace{5px}
\noindent
The norms of the computed higher-order derivatives did not exhibit
systematic growth as the iterates approached optimality. Although they
were not necessarily monotone and could become relatively large at
intermediate iterations, they generally decreased in the local regime.
This behavior suggests that the scaling strategy of Section~\ref{modified_system} is
effective in controlling the higher-order systems, while quantitative
uniform bounds on these derivatives remain an interesting subject for
further study.

\section{Concluding remarks and future research}\label{sec:conclusion}
In this paper, we developed a higher-order feasible predictor-corrector \IPM~for \SDO~based on the \AHO~search directions. The method exploits the semi-algebraic structure of weighted central paths and their ramified analytic reparametrizations near the limiting point. This makes it possible to construct higher-order predictors even when strict complementarity fails. Under a uniform predictor regularity condition controlling these local models as the weight varies, we prove super-linear convergence of the normalized complementarity gap without assuming strict complementarity. We also establish super-linear convergence of the vanishing blocks of the primal-dual iterates, derive a quantitative reduction estimate throughout the predictor-regularity region, and obtain a doubly logarithmic outer iteration bound for the resulting predictor-corrector \IPM.

\vspace{5px}
\noindent
We evaluated the proposed algorithm on five problem sets with diverse characteristics, including four that violate the strict complementarity condition. For suitable choices of $(\rho,p)$, the numerical results exhibit behavior consistent with super-linear convergence across all test sets, including instances without strict complementarity.

\vspace{5px}
\noindent
\paragraph{\textbf{Local convergence and global iteration complexity}.}
We establish super-linear convergence of the complementarity gap and an explicit quantitative reduction estimate when $p\geq\rho$. It would be of considerable interest to further investigate the convergence behavior of Algorithm~\ref{alg:higher-order-IPM-imp} in the regime $p<\rho$. Another important question is whether the higher-order method itself admits a polynomial global iteration complexity bound in the standard input-size model. Our present global estimate depends on the predictor regularity constant $c$; establishing polynomially controlled bounds for the normalized higher-order \AHO~systems may provide a route toward such a result.

\vspace{5px}
\noindent
\paragraph{\textbf{Corrector complexity and higher-order conditioning}}
In all of our numerical experiments, the centering phase required only full \AHO~steps, and in most iterations a single full step was
sufficient. It would be interesting to identify conditions under which one full centering step is guaranteed locally, or more generally to
obtain a uniform bound on the number of corrector steps required after a higher-order predictor. The experiments also show that the norms of
higher-order derivatives can vary substantially with the iteration, the weight, and the derivative order. Developing quantitative uniform
bounds for these derivatives, and relating such bounds to the predictor regularity condition, could lead both to sharper convergence results and to more robust numerical implementations.

\vspace{5px}
\noindent
\paragraph{\textbf{Approximation of $\rho$}}
We did not address the problem of estimating an admissible reparametrization exponent $\rho$ for the family of weighted central paths. For the ordinary central path, an exact symbolic procedure for computing a feasible reparametrization exponent was developed in~\cite{BM24}. In our numerical experiments, we instead use heuristic candidate values motivated by~\cite{HLMZ21}. Developing reliable numerical techniques for estimating and certifying an admissible exponent for the weighted family remains an interesting direction for future research.

\section*{Acknowledgments} 
\noindent
The authors were supported by NSF grants CCF-2128527, CCF-1910441 and CCF-2128702 during the preparation of the initial version of this paper.

\bibliographystyle{abbrv}
\bibliography{MathSciNetBIB}
\end{document}